\input amstex
\documentstyle{amams} 
\document
\annalsline{153}{2001}
\startingpage{471}
\font\tenrm=cmr10

\catcode`\@=11
\font\twelvemsb=msbm10 scaled 1100

\font\ninemsb=msbm10 scaled 800
\newfam\msbfam
\textfont\msbfam=\twelvemsb  \scriptfont\msbfam=\ninemsb
  \scriptscriptfont\msbfam=\ninemsb
\def\msb@{\hexnumber@\msbfam}
\def\Bbb{\relax\ifmmode\let\next\Bbb@\else
 \def\next{\errmessage{Use \string\Bbb\space only in math
mode}}\fi\next}
\def\Bbb@#1{{\Bbb@@{#1}}}
\def\Bbb@@#1{\fam\msbfam#1}
\catcode`\@=12

 \catcode`\@=11
\font\twelveeuf=eufm10 scaled 1100
\font\teneuf=eufm10
\font\nineeuf=eufm7 scaled 1100
\newfam\euffam
\textfont\euffam=\twelveeuf  \scriptfont\euffam=\teneuf
  \scriptscriptfont\euffam=\nineeuf
\def\euf@{\hexnumber@\euffam}
\def\frak{\relax\ifmmode\let\next\frak@\else
 \def\next{\errmessage{Use \string\frak\space only in math
mode}}\fi\next}
\def\frak@#1{{\frak@@{#1}}}
\def\frak@@#1{\fam\euffam#1}
\catcode`\@=12


\define\s{\sigma}
\define\la{\lambda}
\define\al{\alpha}

\define\F{\Cal F}
\define\G{\Cal G}

\define\th{\theta}
\define\sing{\text{ sing }}
\define\pprime{{\prime\prime}}

\define\trace{\operatorname{trace}}
\define\grad{\operatorname{grad}}

\define\res{\operatorname{res}}
\define\Hess{\operatorname{Hess}}

\define\TR{\operatorname{TR}}
\define\KERD{\operatorname{KERD}}
\define\RES{\operatorname{RES}}
\define\RESD{\operatorname{RESD}}

\define\detp{{\operatorname{det}^{\prime}}}

\define\lan{\langle}
\define\ran{\rangle}

\define\reg{{\text{reg}}}
\redefine\sing{{\text{sing}}}

\define\M{\Cal M}

\redefine\D{\Cal D}
\define\ord{{\operatorname{ord}}}

\define\mer{{\text{mer}}}
\define\diff{{\operatorname{diff}}}
\define\Diff{{\operatorname{Diff}}}
\define\conf{{\operatorname{conf}}}
\define\Conf{{\operatorname{Conf}}}

\title{Critical metrics for the determinant of\\ the Laplacian in odd dimensions} 
\shorttitle{Determinant of the Laplacian}

\acknowledgements{The author was supported by the National Science Foundation and the Alfred P. Sloan Foundation.}
 \author{K. Okikiolu}
   \institutions{University of California San Diego, La Jolla, CA\\
{\eightpoint {\it E-mail address\/}: okikiolu@math.ucsd.edu}}
\bigbreak \centerline{\bf Abstract}
\medbreak

Let $M$ be a  closed compact  $n$-dimensional manifold with $n$ odd.
We calculate the first and second variations of the zeta-regularized
determinants $\detp \Delta$ and $\det L$ as the metric on $M$ varies,
where $\Delta$ denotes the Laplacian on functions 
and $L$ denotes the conformal Laplacian.   We see that
the behavior of these functionals depends on the dimension.
Indeed, every critical metric for $(-1)^{(n-1)/2}\detp\Delta$ or 
$(-1)^{(n-1)/2}|\det L|$ 
has finite index. Consequently  there are no local maxima
if $n=4m+1$ and no local minima if $n=4m+3$. We show that the
standard $3$-sphere is a local maximum for $\detp\Delta$ while 
the standard $(4m+3)$-sphere with $m=1,2,\dots,$ is a saddle point.
By contrast, for all odd $n$, the standard $n$-sphere is a local
extremal for $\det L$.  

An important tool in our work is the 
{\it canonical trace} on odd class operators in odd dimensions.
This trace is related to the determinant by the formula
$\det Q=\TR \log Q$, and we prove some basic results on how
to calculate this trace.

\def\sni#1{\smallbreak\noindent{#1}. }
\def\ssni#1{\vglue-1pt\noindent\hskip18pt {#1}. }

\bigbreak

 \centerline{\bf Contents}
 \sni{1} Introduction
\ssni{1.1} The determinant of the Laplacian
\ssni{1.2} The canonical trace
\ssni{1.3} Critical metrics
\sni{2} Regularization theorems
\ssni{2.1} Symbols
\ssni{2.2} Canonical splitting of operators
\ssni{2.3}  Canonical trace for products
\sni{3} Variation formulas
\ssni{3.1} First variation of det$'\Delta$
\ssni{3.2} First variation of det $L$
\ssni{3.3} Second variation of det$'\Delta$
\ssni{3.4} Second variation of det $L$
\sni{4} Critical points have finite index
\sni{5} Standard spheres which are extremal
\ssni{5.1} Fourier series for singular kernels
\ssni{5.2} The hessian of $F$ at $S^3$
\ssni{5.3} The hessian of log det $L$ at $S^n$
\sni{6} Saddle points for det$'\Delta$  

\section{Introduction}

In this paper we prove several results about critical
metrics for the determinant of the Laplacian.  One of our main analytic tools
is the {\it canonical trace} on odd class operators, and we also present here 
some basic results on how to compute this trace.  In Section 1.1 we give a brief 
introduction to the determinant
of the Laplacian, in Section 1.2 we state  results on the canonical trace
and in Section 1.3 we state the results on critical metrics for the
determinant of the Laplacian.
 
\demo{{\rm 1.1.}  The determinant of the Laplacian} 
Let $M$ be a smooth compact $n$-dimensional manifold
without boundary, let $g$ be a 
Riemannian metric  on $M$ and let $\Delta:C^\infty(M)\to C^\infty(M)$ be the 
(positive) Laplacian for the metric $g$.
The conformal Laplacian is 
$$
L\ =\ \Delta\ +\ \al \operatorname{S},\qquad\qquad\qquad \al=\frac{n-2}{4(n-1)},
$$
where $\operatorname{S}$  denotes the scalar curvature.
We will define the {\it zeta regularized determinants}
$\det L$ and $\detp \Delta$.

Let $B:C^\infty(M)\to C^\infty(M) $ be a second order elliptic differential operator
which is bounded below, with eigenvalues
$$
\la_1 \, \leq \,  \la_2 \, \leq \, \dots 
$$
counted by multiplicity.
Define
$$
Z(s)\ =\ Z(B, s)\ =\ \sum_{\la_j\neq 0} \la_j^{-s}
\tag 1.1.1
$$
where the branch of $\la^s$ is chosen so that $\la^s>0$ when $\la>0$
and $s\in \Bbb R$.
The sum in (1.1.1) converges for $\Re s>n/2$ by Weyl's law, and can 
be analytically continued to a meromorphic function on $\Bbb C$
which is regular at $s=0$; see [S].
When $B$ is invertible, define
$$
 \det B\ =\ e^{-Z^\prime(0)},
\tag 1.1.2
$$
which is formally equal to the product of the
eigenvalues of $L$. 
When $B$ is not invertible, define $\det B=0$ and 
$$
 \detp B\ =\ e^{-Z^\prime(0)}.
$$
 For some general theory of zeta-regularized determinants, see 
[BKF1], [BKF2],  [Fo1], [Fo2], [Fo3], [KV],
[Ok1], [Ok2].

The zeta regularized determinant of the Laplacian 
was introduced in [RS] in connection with topology,   
and the idea was then  taken up by physicists
in order to formally evaluate 
Gaussian functional integrals.  In particular A.~Polyakov applied 
$\detp\Delta$ on compact surfaces to string theory.  
A few years later, B. Osgood, R. Phillips and P. Sarnak studied the determinant
of the Laplacian on surfaces and investigated its applications to spectral
geometry, see [OPS1], [OPS2], [OPS3].\enddemo

\nonumproclaim{Theorem {\rm [OPS1]}}
If $M$ is a closed surface then of all metrics in a given conformal
class and of a given area{\rm ,}  the uniform metric has the maximum
determinant{\rm .}
\endproclaim

\demo{{R}emarks}
  1. All metrics on the $2$-sphere are conformally equivalent and 
so among all metrics on the $2$-sphere of area $4\pi$,
the standard $2$-sphere maximizes the determinant of the Laplacian.
This special case of the result was first proved by Onofri [On]. 

\smallbreak  2.  In [OPS2], the functional $\log\detp\Delta$ is  used to
prove that the space of isospectral metrics on a surface is compact
in the smooth topology.  
\enddemo

In higher dimensions, the functional $\det  L$ has been studied 
 by T. Branson, A. Chang, M. Gursky,  B. {\O}rsted, J. Qing, P. Yang and others; see 
[BO1], [BO2], [BCY], 
[CY], [CQ1], [CQ2], [CQ3], [Gur].

\nonumproclaim{Theorem {\rm [CY]}} For the following cases of $4$\/{\rm -}\/manifolds $M${\rm ,}
the standard metric is a global maximum
for $\det L$ among metrics in the same conformal class and with the same volume{\rm :} 
\smallbreak
 $S^4${\rm ,} $\Bbb C P^2${\rm ,} $S^2\times S^2${\rm ,} $\Bbb R^4/\Gamma$ for any lattice $\Gamma${\rm ,} 
$\Bbb H^2\times\Bbb H^2/\Gamma$ for any lattice $\Gamma${\rm ,}  $\Bbb C H^2/\Gamma$ for any lattice $\Gamma${\rm ,}
$\Sigma\times S^2$ with $\Sigma$ hyperbolic and those K{\rm \"{\it a}}hler\/{\rm -}\/Einstein surfaces which 
are not locally symmetric{\rm .}
 \endproclaim

Here, we collect some extremal properties of spheres for the functionals
$\detp\Delta$ and $\det L$.  The results can be found in 
 [On], [Ri], [BCY], [Br].
\medbreak

\phantom{down a bit}
\vglue-8pt

{\ninepoint
\vbox{\offinterlineskip
\hrule
\halign{&\vrule#&
  \strut \enspace\hfil#\hfil\enspace \cr
height8pt&\omit&&\omit&&\omit&&\omit&&\omit&\cr
&    standard \hfil    &&   extremal type &&  for&& among metrics fixing && proved by  &\cr
height8pt&\omit&&\omit&&\omit&&\omit&&\omit&\cr
\noalign{\hrule}
height8pt&\omit&&\omit&&\omit&&\omit&&\omit&\cr
& $S^2$   &&  global max  && $\det^\prime\Delta$ && area && Onofri &\cr
height8pt&\omit&&\omit&&\omit&&\omit&&\omit&\cr
& $S^3$  &&  local max  && $\det^\prime\Delta$  && volume + conformal class && Richardson &\cr
height8pt&\omit&&\omit&&\omit&&\omit&&\omit&\cr
& $S^4$  &&  global min  && $\det L$  && volume + conformal class && Branson-Chang-Yang & \cr
height8pt&\omit&&\omit&&\omit&&\omit&&\omit&\cr
& $S^6$  &&  global max  && $\det L$  && volume + conformal class && Branson &\cr
height8pt&\omit&&\omit&&\omit&&\omit&&\omit&\cr}
\hrule}

}

\demo{{R}emarks}
 1. 
Notice that whether $S^n$ is a maximum or a minimum
for the conformal Laplacian under conformal deformation
depends on the dimension.

\smallbreak 2. In odd dimensions, $\det L$ is constant on
each conformal class.

\smallbreak 3.  A. Chang and  J. Qing  studied $\det  L$ for
manifolds with boundary.   In particular in [CQ2] they consider the conformal
Laplacian on the unit ball $B^4$ in $\Bbb R^4$ with Robin boundary
conditions on $S^3$,  and show that for metrics
in the conformal class of the standard metric on $B^4$
with fixed volume,  the standard metric is extremal for $\detp  L$.
\enddemo

An important  tool in proving the $2$-dimensional results and
the results involving $\det L$ is an explicit 
formula for the variation of $\log\det L$ under
conformal variations  of the metric $g_0$ which  is 
computable from   the local data
$g_0$ and $\dot g(0)$; see [OPS1], [BO1], [Ok2].  On surfaces this is
the {\it Polyakov-Ray-Singer variation formula}. 
Suppose that  $g_0$ is a fixed metric on the
surface $M$.  Let $g=e^{2\phi}g_0$
be a conformal metric.  Here, $\phi$ is a smooth function on~$M$.
Write $\Delta_0$, $\Delta$ for the Laplacians in the metrics $g_0,g$ 
respectively, so that $\Delta=e^{-2\phi} \Delta_0$.  Then
$$
\multline
 \log\detp \Delta\ -\ \log\detp \Delta_0\ = \ 
-\frac{1}{6\pi} 
\left(  \frac12  \int_M |\nabla_0\phi|^2\, dA_0
\ +\ \int_M K_0 \phi\, dA_0\right)   \\
 \ +\ \log A \ -\ \log A_0
\endmultline
$$
where $\nabla_0$ is the gradient, $dA_0$ is the area 
element, $K_0$ is the Gaussian curvature  in the metric $g_0$
and $A$, respectively $A_0$, is the area of $M$ in the 
metric $g$, $g_0$. These quantities can be computed from $g_0$ and $\phi$.
In dimension $n>2$,
there is no such formula for $\detp\Delta$ under conformal variations,
or for $\det L$ or $\detp \Delta$ under general
variations, which is why these situations have been
studied less.  Next we state two results that deal with these situations.

\nonumproclaim{Theorem {\rm [Ri]}} Let $g_0$ be a metric on a compact closed 
$3$\/{\rm -}\/manifold{\rm .}  If $g_0$ is a critical point of $\detp\Delta$ for metrics in the  conformal class of $g_0$ with
the same  volume{\rm ,} and if  $-Z(1)\la_1\geq 5${\rm ,} then $g_0$ is a local maximum for $\detp\Delta$
among metrics in this class{\rm .}  The standard metric on the cubic torus is an example of such a metric $g_0${\rm .}
\endproclaim

\demo{{R}emark} The standard metric  $g_0$ on the $3$-sphere does not satisfy
the condition $-Z(1)\la_1\geq 5$, but Richardson dealt with this case separately.
\enddemo

\nonumproclaim{Theorem {\rm [Chi]}}
The determinant of the Laplacian on the space of
flat $3$\/{\rm -}\/tori of volume $1$ has a local maximum at the torus
corresponding to the face\/{\rm -}\/centered cubic lattice{\rm .}
\endproclaim

 For  a Riemannian surface, the functional,
$$
 F\ =\ \log \detp \Delta\ -\ \left( 1- \frac{\chi(M)}{6}\right) \log A,
$$
is homogeneous of degree zero,
where $\chi(M)$ is the Euler characteristic of $M$.
Osgood, Phillips and Sarnak  noticed that if the metric is constrained 
to vary within a conformal class then
the gradient flow for $F$ is equal to the Ricci flow.
$$
\frac{dg}{dt}\ =\ \grad F\ =\ -(K-\overline K)g
$$
where $K$ is the Gauss curvature and  $\overline K$ is the average curvature
$$
\overline K=\frac1A \int_M K\, dA.
$$
This flow was studied by R. Hamilton and B. Chow.

\nonumproclaim{Theorem {\rm [Ha], [Cho]}}
Given any initial metric on the closed compact surface $M${\rm ,}
the Ricci flow exists for all time and converges to a metric
of constant curvature{\rm .}
\endproclaim

\demo{{\rm 1.2.} The canonical trace} 
We summarize the theory of the canonical trace.
Let $M$ be a closed compact $n$-dimensional manifold with $n$  odd, and  
let $E$ be a complex $N$-dimensional vector bundle over $M$.
Write $\Bbb C(N)$ for the complex $N\times N$ matrices.
Let $Q:C^\infty(E) \to C^\infty(E) $
be a {\it polyhomogeneous} pseudodifferential
operator. Choice of local coordinates identifies
a local trivialization of $E$ with $\Omega\times \Bbb C^N$ where
$\Omega\subset\Bbb R^n$, and the symbol of $Q$ 
has an expansion of the form
$$
q(x,\xi)\ \sim\ q_d(x,\xi)\ +\ q_{d-1}(x,\xi)\ +\ \dots, \qquad\qquad |\xi|\to\infty,
\tag 1.2.1
$$
where $q$, $q_j$ are $\Bbb C(N)$-valued functions on $\Omega$ with 
$q_j (x,t\xi)=t^j  q_j (x,\xi)$ for $t>0$. 
The Schwartz kernel  $K(Q,x,y)$ of $Q$ is smooth away from  the diagonal.
If $d<-n$ then $Q$ is trace class, $K(Q,x,y)$ is continuous across the diagonal,
and the trace of $Q$ is given by the formula 
$$
\trace Q\ =\ \int_M \trace K(Q,x,x).
$$
When $d> -n$, $K(Q,x,y)$ blows up at the diagonal and
$K(Q,x,x)$ and $\trace Q$ are  undefined. Substitutes for these quantities 
can often be defined by embedding $Q$
in an analytic family of operators $Q(z)$ and then  analytically continuing  $K(Q(z),x,x)$
and $\trace Q(z)$.

\demo{{R}emark}
Write $\wedge$ for the bundle of smooth densities on $M$.
Away from the diagonal,  $K(Q,x,y)$ is a smooth section of $\pi_1^* E \otimes \pi_2^*(E^\ast\otimes \wedge)$
where\break $\pi_i:M\times M\to M$ is the projection onto the $i^{\rm th}$ factor.
When $d<-n$, $K(Q,x,x)$ is a smooth section of  $E\otimes E^*\otimes\wedge$, which is 
a bundle over $M$.    If we fix local coordinates $(x^1,\dots,x^n)$ on $M$ then we 
can identify the density $a(x)dx^1\dots dx^n$  with the function $a(x)$, 
thus trivializing $\wedge$.  If we also fix a local 
trivialization of $E$, then $K(Q,x,x)$ is identified with a $\Bbb C(N)$-valued function
of the coordinates.
\enddemo

Let $B:C^\infty(M)\to C^\infty(M)$ be a
positive self-adjoint elliptic  differential operator of order $\beta>0$.
 For $x\in M$ and $\Re z <\!\!<0$,
$$
K( QB^{z/\beta}, x,x)
$$
is well-defined.  It extends to a meromorphic function 
with simple poles. The  {\it Guillemin-Wodzicki residue density}
of $Q$ is the  section of $E\otimes E^*\otimes \wedge$ given by
$$
\RESD (Q,x)\ :=\ \res|_{z=0}^\mer \ K( QB^{z/\beta}, x,x)
$$
where the notation $\res|^\mer_{z_0}f(z)$ means the residue
of the analytic continuation of  $f(z)$ at $z=z_0$.
It is independent of the choice of $B$ and can 
be computed in local coordinates from the symbol:
$$
\RESD (Q,x)\ =\ \frac{1}{(2\pi)^n}
\left( \int_{\{|\xi|=1\}}  q_{-n}(x,\xi)\, d\s_{n-1}(\xi)\right) \, d\xi.
\tag 1.2.2
$$
Here, the sphere $|\xi|=1$ is defined with respect to any metric on $M$, and $d\s_{n-1}$ is
surface measure on the sphere;  see [Gui], [Ka], [Wo].
Define
$$
\RES Q\ :=\ \res|_{z=0}^\mer\  \trace QB^{z/\beta}\ =\ \int_M \trace \RESD (Q,x).
$$
In fact the family $QB^{z/\beta}$ can be replaced by any analytic family $Q(z)$ where 
$Q(z)$ has degree $d+z$ and $Q(0)=Q$.
 From (1.2.2) it is clear that for the following classes of operators, 
the local residue vanishes identically.
\medbreak
\item{(a)}  Trace class pseudodifferential operators (that is, operators of order less than $-n$)
\smallbreak \item{(b)}  Operators of nonintegral degree
\smallbreak\item{(c)}  Odd class operators (defined below).

\demo{Definition {\rm 1.2.1}}
Let $q\in C^\infty(\Bbb R^n\setminus 0, \Bbb C(N))$ be  homogeneous of 
degree $j$. Then $q$ has {\it regular parity} if
$$
\ q(-\xi)\ =\ (-1)^j q(\xi),  \tag 1.2.3
$$
and has {\it singular parity} if 
$$
q(-\xi)\ =\ (-1)^{j+1} q(\xi). \tag 1.2.4
$$
Now suppose that $q$ has an expansion
$$
q(\xi)\ \sim\ q_d(\xi)\ +\ q_{d-1}(\xi)\ +\ \dots, \quad\quad \xi\to\infty,
\quad \quad q_j(t\xi)=t^j q_j(\xi), \quad t>0.
\tag 1.2.5
$$
Then $q$ has regular, respectively  singular, parity in $\xi$ as $\xi\to\infty$ if all of the homogeneous functions
$q_j$ have regular, respectively singular,  parity.  Notice that polynomials have regular parity.
\enddemo

Suppose that the dimension $n$ is odd.
The {\it odd class} operators $Q$  are those classical pseudodifferential operators for which the symbol $q(x,\xi)$ has
regular parity in $\xi$ as $\xi\to\infty$.  Such operators $Q$ are also said to satisfy 
the {\it transmission property}.
The odd class operators on $E$ form an algebra containing the 
differential operators and inverses of invertible elliptic operators. 

In the cases (a), (b) and (c), it can be shown that 
$$
\KERD(Q,x)\  :=\ K( QB^{z/\beta},x,x)\bigl|_{z=0}^\mer
\tag 1.2.6
$$
is independent of $B$, and hence the {\it canonical trace}
$$
\TR Q\ :=\ \trace QB^{z/\beta}\bigl|_{z=0}^\mer \ =\ \int_M \trace \KERD(Q,x)
$$
is independent of $B$.  In fact the analytic family $QB^{z/\beta}$ in (1.2.6) can 
be replaced by any analytic family $Q(z)$ where $Q(z)$ has degree $d+z$, $Q(0)=Q$,  and 
the symbol $q(z)$ of $Q(z)$ has the expansion
$$
q(z)(x,\xi)\ \sim\ q_d(z)(x,\xi)\ +\ q_{d-1}(z)(x,\xi)\ +\ \dots,
$$
where $q_j(z)(x,\xi)$ is homogeneous in $\xi$ of degree $j+z$ and 
$$
q_j(z)(x,-\xi)\ =\ (-1)^j q_j(z)(x,\xi).
$$
The canonical trace was introduced by 
M. Kontseivich and S. Vishik who verified some basic properties; see  [KV].
\smallbreak \item{(A)}  
If $Q$ is trace class then  $\dsize{\KERD(Q,x)=K(Q,x,x)}$ and so
$\dsize{\TR Q=\trace Q}$.
\smallbreak \item{(B)}  If $n$ is odd and $\partial$ is a differential operator then
$\dsize{\KERD (\partial,x)=0}$ and so
$\dsize{\TR \partial=0}$.
\smallbreak \item{(C)} 
If $Q$, $R$ and $[Q,R]$ are in class (a), (b) or (c), then
$\dsize{\TR [Q,R]= 0}$.

 \medbreak
It is possible to calculate $\KERD(Q,x)$  from  $K(Q,x,y)$
without making an analytic continuation. In the case when $d$ is
nonintegral this was shown in [KV].  Here  we will show
how to do it in the more delicate case when $Q$ is odd class.
We make a splitting  $Q= Q_\reg+Q_\sing$ with
$Q_\reg$  smoothing, so that  the following result holds.

\nonumproclaim{Theorem 1.2.2}
$$
\KERD (Q,x)\ =\ K(Q_\reg,x,x),\qquad\qquad \TR Q=\TR Q_\reg
\tag 1.2.7
$$
If $\partial$ is a differential operator on $E$ then
$$
\KERD( \partial Q,x)\ =\  \partial_x K(Q_\reg,x, y)\,\bigl|_{y=x},
\qquad\qquad \TR \partial Q=\TR \partial Q_\reg.
$$
\endproclaim

The splitting $Q=Q_\reg+Q_\sing$ is made by examining the asymptotic behavior
of the Schwartz kernel $K(Q,x,y)$ of
$Q$ close to the diagonal, and separating the regular terms from the
singular terms.
The operator $Q_\sing$ is what we call {\it purely singular}, a concept which
we now define.

\demo{Definition {\rm 1.2.3}}
(a) Let $q(\xi)$ be a tempered distribution on $\Bbb R^n$ which is homogeneous
of degree $j$. Then $q$ has {\it regular parity} if it satisfies (1.2.3)
on $\Bbb R^n$ in the distributional sense, and has {\it singular parity}
if it satisfies (1.2.4) in the distributional sense.

\smallbreak
(b)  Suppose $S$ and $S_j$ are $\Bbb C(N)$-valued tempered distributions on $\Bbb R^n$ for $j=d,d+1,\dots$,
and suppose $S_j$ is homogeneous of degree $j$ on $\Bbb R^n$ and smooth away from the origin.  
We write 
$$
S(w)\ \sim\ S_d(w)\ +\ S_{d+1}(w)\ +\ \dots, \qquad\qquad w\to0,
$$
if 
$$
S(w)- S_d(w)-\dots-S_j(w)\ 
\cases   \in\  C^{j}(\Bbb R^n) \\
=\ O(|w|^{j+1}), \qquad w\to 0. \endcases
\tag 1.2.8
$$

\smallbreak (c)  If $S$ satisfies (1.2.8), 
we say that $S$ has {\it regular parity}, respectively {\it singular parity}, in $w$ as $w\to0$ if every
$S_j$ has regular parity, respectively singular parity in $w$.

\smallbreak (d)  When $n$ is odd, the odd class pseudodifferential operator $Q$  is {\it purely singular} if 
the Schwartz kernel in local coordinates has the form $K(Q,x,y)=S(x,x-y)$ where 
$S(x,w)$ has singular parity in $w$ as $w\to 0$.  In particular differential
operators have singular parity. 
\enddemo

It is easy to show that the Schwartz kernel of an odd class operator $Q$  can be 
split in local coordinates as 
$$
K(Q,x,y)\ =\ R(x,x-y)\ +\ S(x,x-y)
$$
where $R$ is smooth and $S(x,w)$ has singular parity as $w\to 0$.  The operator $Q_\sing$ is
formed by patching together the kernels $S(x,x-y)$ corresponding to different coordinate charts, and 
the result is a purely singular operator.  It will be shown that the canonical trace of any purely singular
operator vanishes, which leads to Theorem 1.2.2.  The details of  this argument are given in 
Section~2.2.  

\smallbreak {\it Example} 1.2.4.
 For the Laplacian on $S^3$, the eigenvalues are $k(k+2)$ with multiplicity $(k+1)^2$.  
Explicit calculation   gives
$$
\TR \Delta^{-1}\ =\ Z(1)\ =\ \sum_{k=1}^\infty \frac{ (k+1)^2}{(k(k+2))^s}\biggl|^\mer_{s=1}\ =\ -\frac34.
$$
On the other hand writing $d\mu_3$ for surface measure on
$S^3$ normalized to total volume $1$, and $r$ for the distance between points $x$ and $y$ on $S^3$, 
the Schwartz kernel  $K(\Delta^{-1},x,y)$   has the form
$G(r) d\mu_3(y)$ where $G$  satisfies the  differential equation
$$
\frac 1{\sin^2 r} \frac{d}{dr} \sin^2 r\frac{d}{dr} G(r)\ =\ 1,\qquad\qquad\qquad \int_{0}^\pi G(r)\sin^2 r\, dr=0.
$$
The solution is
$$
G(r)\ =\ \frac{(\pi-r)\cos r }{2\sin r}\ -\ \frac14.
$$
The regular part and singular parts for $r$ small are given by 
$$
G_\reg(r)\ =\ -\frac{r\cos r }{2\sin r}\ -\ \frac14, \qquad\qquad G_\sing(r)\ =\  \frac{\pi\cos r }{2\sin r},
$$
and so by Theorem 1.2.2,
$$
\TR \Delta^{-1}\ =\ \int_{S^3} G_\reg(0)\, d\mu_3\ =\ G_\reg (0)\ =\ -\frac34.
$$
We see that these two calculations agree.
\smallbreak

Next we state a small technical lemma which we need in order to describe our
second result on computing the canonical trace.

\nonumproclaim{Lemma 1.2.5} {\rm (a)} If $n$ is odd and  $q(\xi)\in C^\infty(\Bbb R^n\setminus 0)$ is a 
smooth homogeneous function on $\Bbb R^n\setminus 0$ of degree $d\in \Bbb Z$ and regular parity{\rm ,} then there is
a unique homogeneous distribution with regular parity on $\Bbb R^n$ which equals $q$ on 
$\Bbb R^n\setminus 0${\rm .} 
\endproclaim

An important consequence of this lemma is that it gives  a canonical way to multiply together the
kernels of two odd class operators.    

 For an integrable function $q$, define the Fourier and inverse Fourier transforms of $q$  by
$$
\hat q(w)\ =\  \int_{\Bbb R^n} q(\xi)\, e^{-iw\cdot\xi}\,d\xi,  \quad 
\check q(w)\ =\  \frac1{(2\pi)^n}\int_{\Bbb R^n} q(\xi)\, e^{iw\cdot\xi}\,d\xi\ =\  \frac1{(2\pi)^n}\hat
q(-w).
\tag 1.2.9
$$
This definition extends to tempered distributions on $\Bbb R^n$.
 For a symbol $q(x,\xi)$, define $\hat q(x,w)$ and $\check q(x,w)$ to be the
 Fourier and inverse Fourier transforms respectively in the second variable.
We state the following theorem for vector bundles since we will need this generality 
in applications.  However, it can be easily deduced from the scalar case.

\nonumproclaim{Theorem 1.2.6}
{\rm (a)} Suppose that $E_i\to M,\,i=0,1,2,3$ are  vector bundles over $M${\rm ,} that $A\in 
C^\infty(E_2\otimes E_1^*)$
and $Q:C^\infty(E_0)\to C^\infty (E_1)${\rm ,}
 $P:C^\infty (E_2)\to C^\infty (E_3)$ are 
odd class pseudodifferential operators{\rm .} Then $S:C^\infty (E_2\otimes E_1^*)\to C^\infty (E_3\otimes E_0^*\otimes\wedge)$
defined by
$$
S A(x)\ =\ \KERD( PAQ, x)
$$
is a classical pseudodifferential operator whose Schwartz kernel is 
$$
K(P,x,y)\otimes K(Q,y,x).
$$

{\rm (b)} If $\ord P=\al$ and $\ord Q=\beta$ and $\al,\beta> -n${\rm ,} and 
the principal symbols of $P$ and $Q$ are $p_\al$ and $q_\beta$
respectively{\rm ,} then the principal symbol of $S$ is 
$$
s(x,\xi)\ =\ (\check p_\al \otimes \check q_\beta^\prime)^\wedge(x,\xi)\in  \pi^*( E_3\otimes E_2^*\otimes E_1\otimes E_0^*),
\tag 1.2.10
$$
where $q^\prime_\beta(x,\xi)=q_\beta(x,-\xi)$ and  $\pi:T^*M\to M$ is the natural projection{\rm .}
\endproclaim

\demo{{R}emarks}
1. The definition of odd class operators acting between two different vector bundles over $M$
is identical to that  given in Definition 1.2.1 for operators acting on a single vector bundle.

2. To define the tensor product of distributions $K(P,x,y)\otimes K(Q,y,x)$, decompose $P=P_\reg+P_\sing$, $Q=Q_\reg+Q_\sing$.
Now the products
$$ \align
K(P_\reg,x,y)& \otimes K(Q_\reg,y,x),\\ K(P_\reg,x,y)&\otimes K(Q_\sing,y,x),\\
 K(P_\sing,x,y) &\otimes K(Q_\reg,y,x) \endalign
$$
are well-defined, being products of distributions with smooth functions.  The problematic term is 
$$
K(P_\sing,x,y)\otimes K(Q_\sing,y,x).
\tag 1.2.11
$$
This is well-defined as a smooth function away from the diagonal,
which in local coordinates has regular parity in $x-y$ as $y\to x$.
By Lemma 1.2.5, this extends uniquely over the diagonal to a regular parity distribution, which we
take as the definition of the product in 1.2.11. More explicitly, there is a polyhomogeneous expansion in local coordinates,
$$
\align
& K(P_\sing, x, y)= S(x,x-y),\quad\ S(x,w)\ \sim\ S_\al(x,w)+ S_{\al+1}(x,w)+ \dots , \\
& K(Q_\sing,y,x)=T(x,x-y), \quad\ T(x,w)\ \sim\ T_\beta(x,w)+ T_{\beta+1}(x,w)+ \dots,
\endalign
$$
which leads to an expansion for the product of the functions $S(x,w)$ and $T(x,w)$ {\it away from $w=0$},
$$
S(x,w)\otimes T(x,w)\ \sim\ P_{\al+\beta}(x,w)\ +\ P_{\al+\beta+1}(x,w)\ +\ \dots,
\quad  P_j=\sum_{k+\ell=j} S_k\otimes T_\ell,
$$
so that $P_j(x,w)$ is homogeneous of degree $j$ in $w$.
The difference
$$
S(x,w)\otimes T(x,w)\ -\ \sum_{j\leq -n} P_j(x,w)
$$
is an integrable function and each of the terms $P_j(x,w)$ is homogeneous of regular parity and hence extends
to a regular parity distribution in $w$ by Lemma 1.2.5.  The fact that
this extension is independent of coordinates  follows from the fact that the property of being a
regular or singular parity distribution is invariant under a change of coordinates, which is
proved in Section 2.2.

In Section 5.1, we present some additional results concerning the calculation
of the canonical trace for certain operators on spheres.
\enddemo

\demo{{\rm 1.3.} Critical metrics}
We begin with a result describing the behavior of $\detp\Delta$ and $\det L$ 
in the neighborhood of a general critical metric when the dimension $n$ is
odd. The result states roughly that every critical metric has finite index.
 
Write $\Cal M$ for the space of metrics on $M$.
The tangent space  $T_{g_0}\Cal M$ at the metric $g_0$ 
 is the space $C^\infty(S^2T^*M)$ of smooth symmetric $(0,2)$ tensor fields on $M$.
Some basic theory of the space $\Cal M$ can be found in [Eb].

 For the most part, we will work in a frame which is orthonormal with respect
to a fixed metric $g_0$. 
At each point $x\in M$, take an orthonormal basis
$e_1,\dots,e_n$ for the tangent space $T_xM$.  Unless otherwise stated,
the components of tensor fields are given with respect to this basis,
with all indices written as subscripts.  Covariant derivatives in
the directions $e_i,e_j$ etc. will be denoted by $D_i,D_j$ etc.
or by the use of subscripts.  For example, if $A\in T_{g_0} \M$, then
$$
A_{ij,k\ell}\ =\ D_\ell D_k A_{ij}\ =\ (D^2A)(e_i,e_j,e_k,e_\ell).
$$
In this notation, the metric $g_0$ is the $(0,2)$ tensor field $\delta_{ij}$.
We use the summation convention that any repeated index is summed 
over.  Write $d\mu_0$ for the normalized canonical volume element for the
metric $g_0$, with normalization  $\mu_0(M)=1$.
The space $T_{g_0}\Cal M$ has a scalar product
  $$
 \langle a, b\rangle\ =\  \int_M  a_{ij} b_{ij} \  d\mu_0.
 $$

The group of diffeomorphisms of $M$ acts on $\Cal M$
in a natural way; if\break $\psi:M\to M$ is a diffeomorphism 
then  $\psi(g_0)$  is the pull-back $(\psi^{-1})^*g_0$. 
We denote the orbit of $g_0$ by $\Diff(g_0)$.
Then $\detp \Delta$ is constant on this orbit because every metric
$g\in \Diff(g_0)$ is isometric to $g_0$ and
hence isospectral to $g_0$.  
Write  $\diff(g_0)$ for the
tangent space of $\Diff(g_0)$ at $g_0$.  This space contains
tensor fields on $M$ of the form  $X_{i,j}+X_{j,i}$ where $X$ is any $1$-form.
The perpendicular space $\diff(g_0)^\perp$  contains divergence-free  fields,
that is, fields of the form $A_{ij}$ with $A_{ij,j}=0$. 
Write $\Conf(g_0)$ for the space of metrics conformal to $g_0$.  As we have
already remarked, in odd dimensions $\det L$ is constant on $\Conf(g_0)$. Let
$\conf(g_0)$ be the tangent space to $\Conf(g_0)$ at $g_0$, which
contains tensor fields on $M$ of the form $\phi \delta_{ij}$
where $\phi$ is a smooth function on $M$.  It can be proved using elliptic
regularity theory  that we have splittings 
$$
\align
T_{g_0}\M\ &=\ \diff(g_0)\oplus \diff(g_0)^\perp ,\\ 
T_{g_0}\M\ &=\ \diff(g_0)\ +\ \conf(g_0)\ \oplus\ \left( \diff(g_0)+ \conf(g_0)\right)^\perp,
\endalign
$$
where for a subspace $V$ of $T_{g_0}\Cal M$, $V^\perp$ is the orthogonal complement of $V$ in $T_{g_0}\Cal M$.
Since $\detp\Delta$  scales
like $V^{2/n}$ where $V$ is the volume, it is computationally simpler
to work with the functional
$$
 F\ =\ \log\detp\Delta\ -\  \frac2n\log  V
$$
which is homogeneous of degree zero under scaling.
The metric $g_0$ is critical or maximal
 for $\detp\Delta$ under deformations
which preserve the volume if and only if it is critical,
respectively maximal for $F$
under all deformations.

\nonumproclaim{Theorem 1.3.1}
Suppose the dimension $n\geq 3$ of $M$ is odd and let $g_0$ be critical for $F${\rm .}
Then 
$$
\Hess F(X,Y)\ =\ 0,  \qquad \qquad \text{for all }X\in \diff(g_0), \ Y\in T_{g_0}\M.
$$
Moreover{\rm ,} $\Hess F$ restricted to 
the subspace $(\diff(g_0))^\perp$ of $T_{g_0}\Cal M$ has at most 
finitely many nonpositive eigenvalues if $n=1$ {\rm mod}~$4$
and at most finitely many nonnegative eigenvalues if $n=3$ mod~$4${\rm .}
Similarly{\rm ,}  if $g_0$ is critical for $\log|\det L|$, then
$$
\Hess \log|\det L|(X,Y)\ =\ 0,  \qquad \quad 
\text{for all }X\in \diff(g_0)+\conf(g_0), \ Y\in T_{g_0}\M.
$$
Moreover{\rm ,}
$\Hess \log|\det L|$
restricted to the subspace $(\diff(g_0)+\conf(g_0))^{\perp}$ has at most 
finitely many nonpositive eigenvalues if $n=1$ {\rm mod}~$4$
and at most finitely many nonnegative eigenvalues if $n=3$ {\rm mod}~$4${\rm .}
\endproclaim

This result has the following immediate corollary:

\nonumproclaim{{C}orollary 1.3.2} 
If $n=3$ $\mod 4${\rm ,} then no metric can be a local 
minimum for $\detp\Delta$  under deformations fixing the volume{\rm .}  
If $n=1$ $\mod 4${\rm ,} then no metric can be a local 
maximum for $\detp\Delta$ under deformations fixing the volume{\rm .}  
The same is true for the functional $|\det L|$ where it does not vanish{\rm .}
\endproclaim

Theorem 1.3.1 also indicates that the set of critical Riemannian manifolds for $\detp \Delta$
or $\det L$ is finite dimensional. 
We calculate explicitly what sort of critical point the standard sphere is and 
the results are summarized in the following table.

\nonumproclaim{Theorem 1.3.3} When $n$ is odd{\rm ,} the standard $n$\/{\rm -}\/sphere is critical
for $\det L$ under all smooth deformations{\rm ,} and critical for $\detp\Delta$ under
deformations fixing the volume{\rm .}  We describe the types of critical points in several cases 
in the following table{\rm .}
\endproclaim

\bigbreak

{\ninepoint
{\vbox{\offinterlineskip
\hrule
\halign{& \vrule#&
  \strut\quad\ \hfil#\hfil\quad \cr
height8pt&\omit&&\omit&&\omit&&\omit&\cr
& standard sphere\hfil    &&   extremal type &&  for&& among metrics fixing &\cr
height8pt&\omit&&\omit&&\omit&&\omit&\cr
\noalign{\hrule}
height8pt&\omit&&\omit&&\omit&&\omit&\cr
& $S^3$  &&  local max  && $\det^\prime\Delta$  && volume + smooth structure\rlap*  & \cr
height8pt&\omit&&\omit&&\omit&&\omit&\cr
& $S^{4m+3},\, m=1,2,\dots$  &&  saddle  && $\det^\prime \Delta$  && volume + conformal class &\cr
height8pt&\omit&&\omit&&\omit&&\omit&\cr
& $S^{4m+3}$    &&  local max  && $\det L$ && smooth structure \rlap{**} &\cr
height8pt&\omit&&\omit&&\omit&&\omit&\cr
& $S^{4m+1}$  &&  local min  && $\det L$  && smooth structure\rlap{**}  &\cr
height8pt&\omit&&\omit&&\omit&&\omit&\cr}
\hrule}
\noindent *extends Richardson's result outside the conformal class, \hfill\break
**only nontrivial across conformal classes.}} 
\bigbreak

In particular, notice that the standard $3$-sphere is  a local 
maximum for $\detp\Delta$, while the standard $4m+3$-sphere is 
a saddle for $m=1,2,\dots$.  
There is overwhelming evidence that for all odd $n>5$, the $n$-sphere is a saddle point for
$\detp\Delta$ under conformal deformations preserving the volume. 
We prove the case $n=4m+3$ in Section $6$. 
The case $n=4m+1$ is computationally more complicated, but we have checked 
$n=9,13,17$ with the help of Mathematica. In fact there is a clear pattern
 to how the Hessian of the determinant of the Laplacian behaves
for conformal deformations, and this is described at the end of this section. 

The local extrema in the above table are all strict in the sense of the following results.
\nonumproclaim{Theorem 1.3.4}
Let $g_0$ be the standard  metric on $S^3$  and
let $g(t)$ be a deformation of $g_0${\rm .}  If $g(t)$ fixes the 
volume  and $\dot g(0)\notin \diff(g_0)${\rm ,} then 
$\detp \Delta_{g(t)}$ has a strict local maximum at $t=0${\rm .}
\endproclaim

\nonumproclaim{Theorem 1.3.5}
Let $g_0$ be the standard  metric on $S^n$ where $n$ is odd and $n\geq 3${\rm ,} and
let $g(t)$ be a deformation of $g_0${\rm .}  If
$\dot g(0)\notin \diff(g_0)+\conf(g_0)$ then
$\det L_{g(t)}$ has a strict local maximum at $t=0$ if $n=3\,\mod 4${\rm ,}
and a strict local minimum if $n=1\,\mod 4${\rm .}
\endproclaim

We will prove Theorems 1.3.4 and 1.3.5 in this paper.   From Theorem 1.3.4, 
it follows by applying the Nash-Moser implicit function theorem that 
if $g(t)$ is a deformation of the standard metric $g_0$ on $S^3$ with 
$\dot g(0)\in \diff(g_0)$, then $\detp\Delta_{g(t)}$ has a local maximum at $t=0$.
The details will be given in a subsequent paper [Ok4], where  upper bounds on
$\detp\Delta$ for metrics close to $S^3$ will also be given, as well as the 
analogous results for $\det L$ close to standard spheres.

Our next results  characterize the  critical metrics for $\det L$ and
$\detp\Delta$.
Let $\mu_0$ be the canonical
volume element for $g_0$ normalized so that  $\mu_0(M)=1$.  Let $Q$ be $\Delta_0$ or $L_0$ and 
let $G(x,y)$ be  Green's function for $Q$; that is, 
$G(x,y)\,d\mu_0(y)$ is the Schwartz kernel for $Q^{-1}$.   Then $G(x,y)$ is
a smooth function on $M\times M$ except at the diagonal where
it blows up. Define $\G(x,y)$ to be the regular part of the Green function;
that is, $\G(x,y)\,d\mu_0(y)$ is the Schwartz kernel for $(Q^{-1})_\reg$.
The operator $Q^{-1}_\reg$ was defined in Section 1.2.

\nonumproclaim{Theorem 1.3.6 {\rm (first variation formula for $F$)}}
 Let $g_0$ be a metric on $M$ which has odd dimension
$n\geq 3${\rm .} Write $\G(x,y)$ for the regular part of the Green function for $\Delta_0${\rm .}
Let $g(t)$ be a deformation of $g_0${\rm ,} and write $\dot g(0)=2\phi \delta_{ij}+A_{ij}$
where $A_{ii}=0${\rm .}  Then
$$
\frac{d F(g(t))}{dt}\biggl|_{t=0}
\ =\ \int  \left( (n-2)\phi \delta_{ij}-A_{ij}\right) \left(D_{x,j} D_{y,i}\G(x,y)\bigl|_{y=x}
+\frac1n \delta_{ij}\right)\, d\mu. 
$$
\endproclaim

\demo{{R}emark} The notation $D_{x,i} D_{y,j} \G(x,y)\bigl|_{y=x}$ 
should be interpreted in the following sense.  The tangent and
cotangent spaces of $M\times M$ split canonically  as products and
this enables us to  take covariant derivatives in each variable
separately in a canonical way.  (In fact the covariant derivatives 
$D_{x,i}$ and $D_{y,j}$ commute
with each other.) The $(0,2)$ tensor field $ D_{x,i} D_{y,j} G(x,y)$ on $M\times M$
can be pulled back to $M$ by the diagonal map $x\to (x,x)$.
\enddemo

Theorem 1.3.6 leads to the following characterization of critical metrics for $F$.
\nonumproclaim{{C}orollary 1.3.7} When the dimension $n\geq 3$ of $M$ is odd{\rm ,} the following conditions 
{\rm (1.3.1), (1.3.2), (1.3.3)} and {\rm (1.3.4)} are equivalent\/{\rm :}  
$$ 
\align
& \text{The metric }g_0\text{ is  critical for }\detp \Delta\text{ under
all}\tag 1.3.1 \\  
& \text{deformations of }g_0  \text{ which fix the volume.}
\endalign
$$
$$
D_{y,i} D_{x,j} \G(x,y)\bigl|_{y=x}    \ =\ -\frac{ \delta_{ij}}n, 
\tag 1.3.2
$$
$$
\align
& \text{{\rm (a)}}\qquad \G(x,x)\  \text{ is independent of }x  \tag 1.3.3 \\
\text{and} &    \\
& \text{{\rm (b)}} \qquad (D_{x,i} D_{x,j} 
 \G(x,y))\bigl|_{y=x}\ =\ \frac{\delta_{ij}}n. 
\endalign
$$
$$
\align
& \text{\rm (a)}\qquad \G(x,x)\  \text{ is independent of }x    \tag 1.3.4\\
\text{and} &   \\
& \text{\rm (b)} \qquad (D_{x,i} D_{x,j} 
 \G(x,y))\bigl|_{y=x}\ =\ c\delta_{ij},\qquad  c\text{ constant}.
\endalign
$$
Here{\rm ,} $\G(x,y)$ is the regular part of the Green function of $\Delta_{g_0}${\rm .}
\endproclaim

\demo{{R}emark} The condition  (1.3.4)(a) holds if and only if $g_0$ is critical
for $\detp\Delta$ under {\it conformal} deformations which fix the volume.
This was proved by  K. Richardson in [Ri].
 The constant $\G(x,x)$  equals $Z(1)$ where $Z(s)$ is the
zeta function for $\Delta_0$.
\enddemo

We have similar results for the conformal Laplacian.

\nonumproclaim{Theorem 1.3.8}
 Suppose that the dimension $n\geq 3$ of 
$M$ is odd{\rm .}  Write $G(x,y)$ for  Green\/{\rm '}\/s function for
 $L${\rm ,} and $\G(x,y)$ for the regular part of $G(x,y)${\rm .}
Write $R_{ij}$ for the Ricci curvature tensor{\rm .}
Let $g_0$ be a metric on $M$ for which $L_{g_0}$ is invertible
and let $g(t)$ be a deformation of $g_0${\rm .}  Decompose $\dot g_{ij}(0)=2\phi\delta_{ij}+A_{ij}$
where $A_{ii}=0${\rm .}  Now{\rm ,} 
$$
\multline
 \frac{d\log\det L}{dt}  =\ 
 \int \dot g_{ij}(0)\Bigg( \left( D_{x,j}D_{x,i} -\al R_{ij}\right)     \G(x,y)\bigl|_{y=x}\\
-\frac{n}{4(n-1)} \left( D_iD_j+\frac{\delta_{ij}}n \Delta \right)\G(x,x) \Bigg)\, d\mu_0 \endmultline
$$
$$\multline
\phantom{\frac{d\log\det L}{dt} } =\ \int A_{ij}(0)\Bigg( \left( D_{x,j}D_{x,i} -\al R_{ij}\right)     \G(x,y)\bigl|_{y=x}\\
   -\frac{n}{4(n-1)} D_iD_j(\G(x,x)) \Bigg)\, d\mu_0. \endmultline
$$
\endproclaim

This leads to the following characterization of critical metrics for $\log\det L$.

\nonumproclaim{{C}orollary 1.3.9}
Suppose that the dimension $n\geq 3$ of 
$M$ is odd{\rm .} 
The following conditions 
{\rm (1.3.5), (1.3.6)} and {\rm (1.3.7)}  are equivalent{\rm .}
$$ 
\align
& \text{The metric }g_0\text{ is  critical for }\det L\text{ under
all}\tag 1.3.5 \\  
& \text{deformations of }g_0.  
\endalign
$$
$$
\multline
  \left(-D_{x,i}D_{x,j}+\frac{n-2}{4(n-1)}\, 
R_{ij}(x)\right) \G(x,y)\, \biggl|_{y=x} \\
 +\  \frac{n}{4(n-1)}D_i D_j (\G(x,x))\ =\ c_0 \delta_{ij}
\quad \text{ with } c_0 \text {  constant}.
\endmultline  \tag 1.3.6 
$$
$$\multline
\left(-D_{x,i}D_{x,j}+\frac{n-2}{4(n-1)}\, 
R_{ij}(x)\right) \G(x,y)\, \biggl|_{y=x}\enspace
\\
=\  -\frac{1}{4(n-1)}\left( nD_i D_j +\delta_{ ij}\Delta_0\right)(\G(x,x)),
\endmultline \tag 1.3.7
$$
where $R_{ij}$ is the Ricci curvature tensor
and  $\G(x,y)$ is the regular part of the Green function of $L_{g_0}${\rm .}
\endproclaim

\demo{{R}emark}
It is clear that the standard metric on $S^n$ satisfies the conditions
(1.3.4) and (1.3.6), and hence is critical for $\detp\Delta$
and $\det L$.
\enddemo

To prove Theorem 1.3.8, we have
$$
\frac{d}{dt}\log\det L 
\ =\ \TR \dot L L^{-1}\ =\ \int_M \dot L _x \G(x,y)\bigl|_{y=x}\, dx,
\tag 1.3.8
$$
where $\G$ is the regular part of the Green function for $L$ and the first equality is established
in [KV] and the second equality follows from Theorem 1.2.2.
Computing the right-hand side of (1.3.8) gives Theorem 1.3.8.
A similar calculation   holds for the ordinary Laplacian and 
all the details are  carried out
in Sections 2.1 and 2.2.

Now we consider the second variation of the functionals $F$ and $\log\det L$.  We have
$$
\frac{d^2}{dt^2}\log\det L 
\ =\ \frac{d}{dt}\TR \dot L L^{-1}\ =\ \TR \ddot L L^{-1}\,+\, \TR \dot L L^{-1} \dot L L^{-1}.
\tag 1.3.9
$$  
By Theorem 1.2.2, 
$$
\TR \ddot L L^{-1}\ =\ \int_M \ddot L _x \G(x,y)\bigl|_{y=x}\, dx.
$$
The term $\TR \dot L L^{-1} \dot L L^{-1}$ can be somewhat simplified.  The functional $F$
can be treated similarly, leading to the following second variation formulas.

\nonumproclaim{Theorem 1.3.10}
Suppose $\phi\in C^\infty(M)${\rm ,}  let $g_0$ be a critical metric for
$F$ and let $A$ be a symmetric  $(0,2)$\/{\rm -}\/tensor field on $M$ which is 
pointwise trace\/{\rm -}\/free with respect to $g_0${\rm ;} that is{\rm ,}  $A_{ii}=0${\rm .} Write 
$$
\bar\phi\ =\ \phi-\int_M\phi\, d\mu_0, \qquad\qquad \bar{\Delta}_0^{-1}\ =\ 
\Delta_0^{-1}-Z(1)\Pi_0,
$$
where $Z(s)$ is the zeta function for $\Delta_{g_0}${\rm .}
Then 
$$
\multline
 \Hess F(2\phi g_0+A,2\phi g_0+A)
\ =\  - \frac1n\int_M A_{ij} A_{ij}\,d\mu_0 
\ +\ \frac{(n+2)(n-2)}2  \int_M \bar\phi^2\, d\mu_0  \\ 
 \qquad  -\ \TR \left(  \frac{n-2}2  \phi_{ii}
\ +\ D_i A_{ij} D_j \right) \bar{\Delta}_0^{-1}
\left(  \frac{n-2}2  \phi_{kk}
\ +\ D_k A_{k\ell} D_\ell \right) \bar{\Delta}_0^{-1}.
\endmultline  
\tag 1.3.10
$$
\endproclaim

\nonumproclaim{Theorem 1.3.11} Suppose $g_0$ is critical for $\log\det L$ and suppose $A\in T_{g_0}\M$ has
$A_{ii}=0$ and $A_{ij,j}=0${\rm .}  Then
$$
\align
 \Hess \log\det L(A,A)\ =& \  
\frac{1}{4(n-1)}\TR  (A_{ij}A_{ij})_{kk} L^{-1} 
\tag 1.3.11\\
& +\   \al\TR (A_{ij,k} A_{ik,j}- \tfrac12 A_{ij,k} A_{ij,k}) L^{-1} 
\\ &  -\ \TR  A_{ij}(D_i D_j-\al R_{ij}) L^{-1} A_{k\ell}(D_k D_\ell-\al R_{k\ell}) L^{-1}.
\endalign
$$
\endproclaim

Explicit evaluation of (1.3.10) at the standard $3$-sphere  leads to
Theorem 1.3.4.  Explicit evaluation of (1.3.11) at the standard $n$-sphere
with $n$ odd leads to Theorem 1.3.5. We now explain in more detail
what happens to the functional $F$ close to the standard $n$-sphere for odd $n>3$.

Let $g_0$ denote the standard metric on $S^n$.  From Theorem 1.3.10, 
$$ \multline \Hess F(2\phi g_0, 2\phi g_0)
\ = \  \frac{(n+2)(n-2)}{2}  \int_{S^n} \bar\phi^2\, d\mu_n 
\\
 -\ \frac{(n-2)^2}4 \TR  \phi_{ii} \bar{\Delta}_0^{-1}  \phi_{kk}
\bar{\Delta}_0^{-1}.
\endmultline \tag1.3.12
$$
 From Theorem 1.2.6,
$$
 \TR  \phi_{ii} \bar{\Delta}_0^{-1}  \phi_{kk}
\bar{\Delta}_0^{-1}\ =\       \int_{S^n} \phi_{ii} T \phi_{kk}\, d\mu_n
$$
where $T$ is a pseudodifferential operator with
$$
K(T,x,y)\ =\ K(\bar{\Delta}_0^{-1} ,x,y)K(\bar{\Delta}_0^{-1} ,y,x).
$$
This kernel is invariant under  transformations of $S^n\times S^n$
which preserve the  distance $r$ between 
$x$ and $y$ and hence $T$ is a function of the Laplacian, and
$$
\Hess F(2\phi g_0, 2\phi g_0)\ =\ \int_{S^n} \phi Q\phi \, d\mu_n
$$
where $Q$ is a pseudodifferential operator which is a function of the Laplacian.

Write $\Cal H_k$ for the space of spherical harmonics of degree $k$ on
$S^n$.  This is the space of  homogeneous harmonic, degree $k$ polynomials 
on $\Bbb R^{n+1}$ restricted to $S^n$, and is an eigenspace for
the Laplacian on $S^n$ and hence an eigenspace  for  $Q$.
$\Cal H_0$ is the space of constants, and  $Q$ vanishes on $\Cal H_0$
since $F$ is scale invariant.
$\Cal H_1$ is the space of linear functions on $\Bbb R^n$ restricted to $S^n$ and
equals  $\conf(g_0)\cap \diff(g_0)$, so $Q$ vanishes on $\Cal H_1$.  
By Theorem 1.3.1, we see that when $n=3$ $\mod 4$, 
$Q$ is strictly negative on a space
of finite codimension; so for $k$ sufficiently large, $Q$ is strictly negative on $\Cal H_k$.

\phantom{growl}

\nonumproclaim{Lemma 1.3.12}
If $n>3$ is odd and $\phi\in\Cal H_{2}$ then 
$$
\align
& \Hess F(2\phi g_0, 2\phi g_0)
\ =\ \frac{ 2(n-2)\left(-2 + (n-2)(n+1)\sum_{j=2}^{n-1}\frac1j\right)}
{(n-3)n}\int_{S^n} \phi^2\, d\mu_n.
\endalign
$$
\endproclaim

Since the right-hand side is clearly positive for $n>3$, we see that if $n=3$ ${\rm mod}\ 4$, then $Q$
has both positive and negative eigenvalues and $S^n$ is a saddle point for $F$.  
Concerning the case $n=1$ $\mod 4$, computer calculations indicate 
that $S^5$ is a local minimum for $F$ among metrics in the same conformal class, 
while $S^n$ is a saddle point for all odd $n>5$. Indeed, it appears that 
the index of this saddle point has a simple description, since 
for low values of $k$, the sign of $\Hess F$ on $\Cal H_k$
appears to have a $4$-fold periodicity in $k$. 

\phantom{growl}

\nonumproclaim{Conjecture 1.3.13}
$$
\align
 1<k<n-1\quad \Rightarrow \quad & \cases \Hess F <0 \text{ on }\Cal H_k
& \text{ if } k=1\ \mod 4, \\
 \Hess F >0 \text{ on }\Cal H_k
& \text{ if } k=2,3,0\ \mod 4, \endcases  \\
\noalign{\vskip6pt}
 k\geq n-1\quad \Rightarrow \quad & \cases \Hess F <0 \text{ on }\Cal H_k
& \text{ if } n=3\ \mod 4, \\
 \Hess F >0 \text{ on }\Cal H_k
& \text{ if } n=1\ \mod 4. \endcases  \\
\endalign
$$
\endproclaim

\phantom{growl}

This conjecture is confirmed by computer calculations using Mathematica for 
the case that $n$ is odd,  $3\leq n\leq 17$
and $k\leq 20$.   

\phantom{growl}

The outline of the paper is as follows.
Section 2 contains the proofs of the  results on the canonical trace 
described in Section 1.2, including proofs of Theorems 1.2.2 and 1.2.6.

In Sections 3.1 and 3.2 we compute  first variation formulas
for $F$ and $\log \det L$ and characterize the critical metrics; that is, we
prove Theorems 1.3.6--1.3.9. 
In Sections 3.3 and 3.4 we compute the second variation formulas
for  $F$ and $\log \det L$ given in (1.3.10) and (1.3.11).

In Section 4 we apply Theorem 1.2.6 to
analyze those terms in (1.3.10) and (1.3.11) which involve the canonical trace.  
This leads to the proof that critical points of $\detp\Delta$ and $\det L$
have finite index as described in  Theorem 1.3.1. 

In Section 5 we compute (1.3.10) and (1.3.11) explicitly at the 
standard spheres in certain dimensions to verify the  extremal properties 
of Theorems 1.3.3 and 1.3.4, while 
in Section 6, we compute (1.3.10) explicitly at the standard spheres
for conformal deformations corresponding to spherical harmonics in $\Cal H_2$, 
to prove Lemma 1.3.12.

Theorem 1.3.1 and the $3$-dimensional case of Theorem 1.3.3 were 
circulated in a preprint in 1997. The proofs appearing in
this paper are slightly simpler than the original ones.
In a sequel to this paper [Ok3], we will describe the behavior of critical
metrics for special values of the zeta function.  We have observed a marked 
difference in behavior between  
$Z(k)$ when it is summable, and when it is not.

I would like to thank Peter Sarnak for suggesting that I investigate the
determinant of the Laplacian on $3$-manifolds.  I would also like
to thank  Alice Chang,  Richard Hamilton, David Jerison, Ken Richardson,
Gang Tian,  Nolan Wallach and  Paul Yang for  useful conversations.
 
\section{Regularization theorems}

In this section we prove the results stated and outlined in Section 1.2.
\demo{{\rm 2.1.} Symbols} 
 For $\Omega\subset \Bbb R^n$,
write $\Cal D(\Omega)$ for the space of smooth compactly supported functions on
$\Omega$ and  $\Cal D^\prime(\Omega)$ for the space of distributions on $\Omega$.
\enddemo

\demo{Proof of Lemma {\rm 1.2.5}}
Consider the function $q(\xi)|\xi|^{z}$.  This is homogeneous in $\xi$ of degree
$z+d$ and hence   when $\Re z+d>-n$ it is locally integrable and gives  a homogeneous distribution 
$q_z(\xi)$ of degree $z+d$.
 For $f$ a Schwartz function, set
$$
f_N(\xi)\ =\ f(\xi)-\sum_{|\al|<N} \frac{\partial_\xi^\al f(0)}{\al!} \xi^\al,
$$
so that $f_N(\xi)=O(|\xi|^N)$.  Then  $\dsize{ \langle  q_z(\xi)\,,\, f(\xi)\rangle }$ equals
$$
\multline \int_{|\xi|>1} q(\xi)|\xi|^{z} f(\xi)\, d\xi
\\
 +\  \int_{|\xi|<1} q(\xi)|\xi|^{z} f_N(\xi)\, d\xi  
\ +\ \int_{|\xi|<1} q(\xi)|\xi|^{z}  \sum_{|\al|<N} \frac{\partial_\xi^\al f(0)}{\al!} \xi^\al   \, d\xi\endmultline
$$
$$\multline =\ \int_{|\xi|>1} q(\xi)|\xi|^{z} f(\xi)\, d\xi
\, +\, \int_{|\xi|<1} q(\xi)|\xi|^{z} f_N(\xi)\, d\xi\\
  +\, \sum_{|\al|<N}\frac{\partial_\xi^\al f(0)}{\al!(z+d+|\al|+n)}
\int_{|\th|=1}q(\th)\th^\al\, d\s_{n-1}(\th),
\endmultline
$$  
where $d\s_{n-1}$ is surface measure on the unit sphere in $\Bbb R^n$. 
Since the final expression defines a meromorphic function of $z$ when $\Re (z+d+N)>-n$ we conclude that
$\langle q_z(\xi)\,,\, f(\xi) \rangle$ extends to a meromorphic function on the whole complex
plane.   Away from the poles, $q(z)$   is a homogeneous distribution of degree $z+d$ and 
in the distributional sense, 
$$
q_z(-\xi)\ =\ (-1)^d q_z(\xi).
$$
Poles can only occur when $z+d$ lies in the set
$$
  -n, \, -n-1, \, -n-2, \dots,
$$
and the residue of $q_z$ at $z=-n-d-j$ is
$$
\sum_{|\al|=j} \frac{ (-1)^{|\alpha|}\partial_\xi^\al \delta_0}{2\al!}
\int_{|\th|=1}q(\th)\th^\al\, d\s(\th),
$$
where $\delta_0$ is the Dirac delta function at $0$.  In particular if $q$ has
regular parity there is no pole at $z=0$ and
$q_0$ is a regular parity homogeneous distribution 
of degree $d$ which equals $q$ away from the origin.

 For the uniqueness, if $q^\prime$ and $q^{\pprime}$ are two homogeneous distributions of regular parity
which equal  $q$ away from zero,
then $q^\prime-q^{\pprime}$ is supported at zero, hence
equals  a differential operator applied to the delta function at zero.  Since $n$ is odd we conclude
that  $q^\prime-q^{\pprime}$ has singular parity, but since it 
is the difference of regular parity distributions it also has regular parity and hence equals zero.
\enddemo 

In considering pseudodifferential operators, it will be convenient to 
allow homogeneous symbols as well as smooth symbols, so we will
define classes of symbols  $q(x,\xi)$ which are smooth functions for $\xi\neq 0$ and  classical conormal
distributions with respect to the set $\xi=0$.

\demo{Definition {\rm 2.1.1}} Write $S^d(\Omega)$ for the space of {\it smooth symbols on $\Omega$ of
order $d$}, that is, the space of functions
$q\in C^\infty (\Omega\times\Bbb R^n)$  satisfying the estimates
$$
\sup_{x\in K, \, |\xi|\geq 1} |\xi|^{-d+|\beta|}
|\partial_x^\al \partial_\xi^\beta q(x,\xi)|\ <\ \infty
\quad  \text{ for each }\al,\,\beta,\text{ and compact } K\subset \Omega.
\tag 2.1.1
$$
\enddemo

\demo{Definition {\rm 2.1.2}} Write $SD^d(\Omega)$ for the space of {\it order $d$ 
distributional symbols on $\Omega\times\Bbb R^n$ which are conormal with respect to $\Omega\times 0$}, that is,  
the space of distributions $q\in \Cal D^\prime(\Omega\times\Bbb R^n)$ such that 
$q(x,\xi)$ is smooth away from $\xi=0$, satisfies (2.2.1), and 
there exists   $N\in \Bbb N$  with 
$$
\xi^{\gamma} \partial_x^\al \partial_\xi^\beta  q(x,\xi)\ \in\ C(\Omega\times\Bbb R^n)
\qquad \text{ for all }\al,\beta,\gamma\quad\text{with } |\gamma|-|\beta|>N.
$$
Write $SD^d(\Omega,\Bbb C(N))$ for the space of $\Bbb C(N)$-valued distributions on $\Omega\times\Bbb R^n$
whose matrix entries are in $SD^d(\Omega)$.
\enddemo

\nonumproclaim{Lemma 2.1.3}  If $q\in SD^d(\Omega)$ where $d< -n${\rm ,}
then   $$\hat q, \quad \check q\in C^{-n-d-1}(\Omega\times\Bbb R^n).$$
\endproclaim

Now we introduce pseudodifferential operators.
Let $\pi:E\to M$ be an $N$-dimensional complex vector bundle over the closed compact
$n$-dimensional manifold $M$ where $n$ is odd.   We will generally work in a coordinate
chart $\Omega$ of $M$ over which $E$ is trivial, so 
$\pi^{-1}(\Omega)$ is identified with $U\times\Bbb C^N$, where $U$ is  a subset of  $\Bbb R^n$.

$Q:C^\infty(E)\to C^\infty(E)$ is a 
pseudodifferential operator of order $d$, if for each  local coordinate chart $\Omega$ 
over which $E$ is trivial
there exists a symbol $q\in S^d(\Omega,\Bbb C(N))$ such that whenever 
$f\in C^\infty(E)$ is supported  $\Omega$, in coordinates
$$
Qf(x)\ =\ \frac1{(2\pi)^n}\int_{\Bbb R^n}\int_{\Bbb R^n} q(x,\xi) e^{i(x-y)\cdot \xi} f(y)\, dy\, d\xi.
$$

Now formally
$$\align
Qf(x)& =\ \int_{\Bbb R^n}\left( \left( \frac1{2\pi}\right)^n
\int_{\Bbb R^n} q(x,\xi)\, e^{i(x-y)\cdot\xi}\,d\xi
\right) f(y)\, dy \\
& =\ \int_{\Bbb R^n} \check q(x,x-y)f(y)\, dy.\endalign
$$
Interpreted distributionally this shows that the Schwartz kernel $K(Q,x,y)$ 
equals the distribution $\check q(x,x-y)$. 

\demo{{\rm 2.2.} Canonical splitting of operators} 
Suppose throughout this section that the dimension $n$ is odd.  Let
$\pi:E\to M$ be an $N$-dimensional complex vector bundle over the closed compact
$n$-dimensional manifold $M$ and let\break $Q:C^\infty(E)\to C^\infty(E)$ be an odd class 
pseudodifferential operator of order $d$ with symbol expansion
$$
q(x,\xi)\ \sim \ q_d(x,\xi)\ +\ q_{d-1}(x,\xi)\ +\ \dots,  \qquad\qquad \xi\to\infty,
\tag 2.2.1
$$
where $q_j(x,\xi)$ is a homogeneous distribution of degree $j$ in $\xi$ with regular parity.
Then
$$
\check q(x,w)\  -\ \sum_{j\geq -n-M} \check q_j(x,w)\ \in\  C^{M}(\Omega\times\Bbb R^n, \Bbb C(N)).
$$
It is convenient at this point to label functions of $w$ by their homogeneity 
by setting
$$
S_j(x,w)\ =\  \check q_{-n-j}(x,w).
$$
It is well-known that 
given an arbitrary power series at a point, there is a smooth function $f$ whose Taylor series
equals that power series.  The function can be constructed by multiplying each term in the
power series by a suitable cut-off function.  The same construction gives  a 
distribution $S(x,w)$ which is
smooth for $w\neq 0$ and satisfies
$$
S(x,w)\ \sim\ \sum_{j=-n-\al}^\infty S_j(x,w),
$$
in the sense of (1.2.8).
Then $R(x,w):=\check q(x,w)-S(x,w)$  is smooth.  Summarizing, we have
$$
K(Q,x,y)\ =\ S(x,x-y)\ +\ R(x,x-y).
\tag 2.2.2
$$
where $R(x,w)$ is smooth and  $S(x,w)$ has singular parity in $w$ as $w\to0$.

Now we form  operators $Q_\reg$ and $Q_\sing$ with $Q=Q_\reg+Q_\sing$
by patching together the kernels $R(x,x-y)$ and $S(x,x-y)$.
Choose a finite partition of unity $\phi_j$ for  $M$ such that the supports of $\phi_j$ and  $\phi_\ell$, 
are either disjoint or contained in a coordinate patch $\Omega$ over which $E$ is trivial.
If the supports are disjoint then
$\phi_j Q \phi_\ell$ is smoothing and we set
$$  
Q_\reg^ {j\ell} \ =\ \phi_j  Q \phi_\ell,
\qquad\qquad\qquad Q_\sing^ {j\ell} \ =\ 0.
$$
If not, the supports are contained in a coordinate patch $\Omega$ and we set
$$ \align
Q_\reg^{j\ell} f(x)& =\ \int_{\Omega} \phi_j(x) R(x,x-y)\phi_\ell(y) f(y)\, dy,
\\
Q_\sing^{j\ell} f(x)& =\ \int_{\Omega} \phi_j(x) S(x,x-y)\phi_\ell(y) f(y)\, dy, \endalign
$$
where $R$ and $S$ are as defined above.  Now we set
$$
Q_\reg \ =\ \sum_{j,\ell} Q_\reg^{j\ell},
\qquad\qquad Q_\sing \ =\ \sum_{j,\ell} Q_\sing^{j\ell}.
$$

See Definition 1.2.3  for  the definition of purely singular operators.

\nonumproclaim{Lemma 2.2.1} {\rm (a)}
The class of purely singular operators contains the differential operators and
is closed under composition with differential operators{\rm .}
\smallbreak
\item{\rm (b)} The adjoint of a  purely singular operator is purely singular{\rm .}
\smallbreak
\item{\rm (c)}  The operator $Q_\sing$ constructed above is purely singular{\rm .}
\smallbreak
\item{\rm (d)} If $Q$ is purely singular  then $\KERD(Q,x)=0$ so that $\TR Q=0${\rm .}
\endproclaim 

\nonumproclaim{{C}orollary 2.2.2} If $\partial$ is a differential operator then
$$
\align
 \TR Q& = \TR Q_\reg + \TR Q_\sing\ =\ \TR Q_\reg,  \\
 \TR \partial Q& =  \TR \partial Q_\reg\ +\ \TR \partial Q_\sing \ =\ \TR \partial Q_\reg.
\endalign
$$
\endproclaim

\demo{{R}emark}  $Q_\reg$ and $Q_\sing$ are only well-defined 
modulo smoothing operators whose integral kernels vanish to infinite order at the diagonal.
\enddemo

\demo{Proof of Lemma {\rm 2.2.1}}
(a)  First we notice that the identity operator $I$ is purely singular, since 
 $K(I,x,y)=\delta(x-y)$, and the delta function $\delta(w)$ is even in $w$ and 
has homogeneity $-n$.   Next notice that 
if $S(x,w)$ is purely singular as $w\to 0$   and $\psi(x,w)$ is 
smooth, then $S(x,w)\psi(x,w)$ is purely singular as $w\to 0$. Indeed, 
$$ \align
S(x,w)& \sim  \sum_{j=d}^\infty S_j(x,w),  \\  \psi(x,w)& \sim  \sum_{j=0}^\infty
\psi_j(x,w) \ \text{ with } \  \psi_j(x,w)=
\sum_{|\al|=j} \frac{\partial_w^\al \psi(x,0)}{\al!} w^\al,\endalign
$$
where $S_j(x,w)$ is homogeneous of degree $j$ in $w$ with singular parity.
Clearly $S_k(x,w) \psi_j(x,w)$ is homogeneous with singular parity in $w$, and 
$$
S(x,w)\psi(x,w)\ \sim\ \sum_{j=d}^\infty \sum_{k+\ell=j} S_k(x,w)\psi_\ell(x,w).
\tag 2.2.3
$$

We  need to show that if 
$$
S(x,w)\ \sim\ \sum_{j=d}^\infty S_j(x,w)
$$
where $S_j(x,w)$ is homogeneous of order $j$ in $w$ with
singular parity, then for a differential operator 
$P(x,D_x)=\sum_{|\beta|\leq N} a_\beta(x)\partial_x^\beta$,
$$
P(x,D_x) S_j(x,x-w)
\ \sim\ \sum_{j=d-N}^\infty \tilde S_j(x,x-y)
\tag 2.2.4
$$
where $S_j(x,w)$ is homogeneous of order $j$ in $w$ with
singular parity.  We have already seen in (2.2.3) that (2.2.4) holds when
$P(x,D_x)$ is multiplication\break by a smooth function.
By  induction it suffices to prove (2.2.4) for 
$P(x,D_x)\break=\partial_{x_k}$, and this case follows easily by application of
the chain rule to each term $S_j(x,x-y)$.
We leave the details for the reader.

We omit the proof of (b), since it is not used in what follows.

To show (c), we want to see that  $Q^{j\ell}_\sing$ 
has singular parity.
Now if the supports of $\phi_j$ and $\phi_\ell$ are not disjoint, then
in some coordinate patch $\Omega$, 
$$
K(Q^{j\ell}_\sing,x,y)\ =\ \phi_j(x) S(x,x-y) \phi_\ell(y)\ =\ S^{j\ell}(x,x-y),
$$
where $S(x,w)$ has singular parity in $w$ as $w\to 0$.  By 
part (a), $S^{j\ell}(x,w)$ also has singular parity in $w$ as $w\to 0$.
We will now  show  that $S(x,w)$ still has singular parity 
under a change of coordinates.
To simplify the presentation we  fix $x$.  It is a trivial modification
to include the dependence on $x$. Suppose  that $\Omega\subset\Bbb R^n$ is open and
contains zero, and suppose that
$S\in \Cal D^\prime(\Bbb R^n)$ is compactly supported in $\Omega$, smooth away from $0$, 
and 
$$
S(w)\ \sim\ \sum_{j\geq d} S_j(w),  \qquad w\to 0,
$$
where $S_j(w)$ is a homogeneous distribution of degree $j$ with singular parity.
Let $F:U\to\Omega$ be a diffeomorphism with $F(0)=0$, and set $w=F(W)$.
We wish to show that
$$
S(F(W))\ \sim\ \sum_{j\geq d} \tilde S_{j}(W), \qquad W\to 0,
\tag 2.2.5
$$
where $\tilde S_j(W)$  is a homogeneous distribution of degree $j$ with singular parity in $W$.
We start by showing that this is true  away from $W=0$, which in fact proves the 
result in the case when $d>-n$.  We will assume that there is only one term in
the expansion of $S$, so that $S$ is homogeneous of degree $j$ for $w$ small. The general
case follows by summing over $j$.  Following a similar calculation in [KV], 
by Taylor's theorem, we have
$$
\align
 F(W)\ & =\ F^\prime(0)W+ \sum_{2\leq |\al|<N} 
\frac{\partial^\al F(0)}{\al!}
W^\al + O(|W|^{N}),\\
S(w_0+h)\ &=\ \sum_{ |\beta|< N}
\frac{\partial^\beta S(w_0)}{\beta!} h^\beta
\ +\ O\left(|h|^{N} |w_0|^{j-N}\right).
\endalign
$$
Hence, setting $w_0=F^\prime(0)W$ and $h=F(W)-F^\prime(0)W$,
$$
S(F(W))  =  \sum_{|\beta|< N} 
\frac{\partial ^\beta S(F^\prime(0)(W))}{\beta!}
\left( \sum_{2\leq |\al|<N} \frac{\partial^\al
 F(0)}{\al!} W^{\al}\right)^\beta
  +  O\left(|W|^{N+j-1}\right).
\tag 2.2.6
$$
Moreover, if the equation is differentiated in  $W$ any number
of times, the error $O(|W|^{N+j-1})$ is reduced in order by the number of
$W$-derivatives.  By grouping terms on the 
right-hand side of (2.2.6) which
are homogeneous of  degree $\ell$, we get (2.2.5).

In the case when $d\leq -n$, write
$$
S(w)\ =\  (1+\Delta_w)^N S_0(w)
$$
where $S_0$  has singular
parity as $w\to 0$, is integrable and smooth away from zero, and $N$ is a positive integer.  
This is easily accomplished by
using the Fourier transform.  A change of variables gives
$$
S(F(W))\ =\ P(W, D_W) S_0(F(W)),
$$
where 
$$
P(W, D_W)  =  \left(1\ -\ \sum_{j,k, \ell} A_{kj} (W) \frac{\partial}{\partial {W_k}}
 A_{\ell j} (W)\frac{\partial}{\partial {W_\ell}}\right)^N,
\enspace   A(W)  =    (F^\prime(W))^{-1}.
$$
We have already shown that $S_0(F(W))$ has singular parity in $W$ as $W\to 0$. 
Hence so does $S(F(W))$ by part (a).
This completes the proof of (c).  The same argument shows that the property of
being a regular parity 
function or distribution
is invariant under a change of coordinates.

 Finally, we prove (d). 
Suppose $Q$ is a purely singular odd class operator of order $d$.  Let $x\in M$, choose 
a coordinate chart $\Omega$ containing $x$ and let 
$q\in S^d(\Omega,\Bbb C(N))$ be the symbol of $Q$.
Let $\phi$ be a smooth positive even function on $(-\infty,\infty)$ with 
$$
\phi(r)= r, \quad r>1,  \quad\qquad \phi(r)=1, \quad |r|<1/2.
$$
Then for $x\in \Omega$,
define an analytic family of operators $Q(z)$ so that on a neighborhood of $x$,  $Q(z)$  has symbol
$$
q_z(x,\xi)\ =\ q(x,\xi) \phi(|\xi|)^z.
$$
Then  
$$
\KERD (Q,x)\ =\ \frac1{(2\pi)^n}\int_{\Bbb R^n} q(x,\xi) \phi(|\xi|)^{z} \, d\xi\biggl|^\mer_{z=0}.
$$
We wish to show that this vanishes.  Suppose $q$ has the expansion (2.2.1).
We will compute
$$
\frac1{(2\pi)^n} \int_{\Bbb R^n} q_j (\xi) (\phi(|\xi|))^{z}\, d\xi\biggl|^\mer_{z=0}.
\tag 2.2.7
$$
If $j>-n$ and $\Re z< -j-n $ then 
$$ \align
\int_{|\xi|>1} q_j (\xi) (\phi(|\xi|))^{z}\, d\xi
&=\ \int_1^\infty r^{j+z+n-1}\, dr\, \int_{|\th|=1} q_j(\th)\, ds(\th)
\\
 & =\ -\frac1{j+z+n} \int_{|\th|=1} q_j(\th)\, ds(\th).\endalign
$$
When we analytically continue to $z=0$ this equals 
$$
-\frac1{j+n} \int_{|\th|=1} q_j(x,\th)\, ds(\th).
$$
On the other hand 
$$
\int_{|\xi|<1} q_j (x,\xi) (\phi(|\xi|))^{z}\, d\xi\biggl|_{z=0}
\ =\ \int_{|\xi|<1} q_j (x,\xi)\, d\xi
\ =\ \frac1{j+n} \int_{|\th|=1} q_j(x,\th)\, ds(\th).
$$
Hence (2.2.7) equals zero.
If $j=-n$ then the distribution $q_j(x,\xi)(\phi(|\xi|))^{z}$ is odd in $\xi$ and
hence integrates to zero when $\Re z <\!< 0$. Hence again, (2.2.7)
vanishes. Now set 
$$
q^{(n)}(x,\xi)\ =\ q(x,\xi)\ -\ \sum_{j\geq -n} q_j(x,\xi)\ \in \ S^{-n-1}(\Omega,\Bbb C^{N\times N}).
$$
Then $\int_M q^{(n)}(x,\xi)\, d\xi=0$.  Indeed, $\check q(x,w)$
has singular parity in $w$ as $w\to 0$, hence so does $\check q^{(n)}(x,w)$, and
$$ \align
\frac1{(2\pi)^2n}\int_{\Bbb R^n}
q^{(n)}(x,\xi)e^{i\xi\cdot w} \, d\xi& =\ \check q^{(n)}(x,w)\\
& =\ \check q(x,w)\ -\ \sum_{j\geq -n} \check q_j(x,w)
\ \in \ C(\Omega\times\Bbb R^n),\endalign
$$
which vanishes as $w\to 0$. \enddemo

2.3. {\it Canonical trace for products}.

\demo{{P}roof of Theorem {\rm 1.2.6}}
Suppose that $P=\phi P\psi$ where $\phi$ and $\psi$ are functions on $M$ with disjoint supports.
Then $P$ and $PAQ$ are smoothing and 
$$\align
TA(x)\ =\ \KERD(PAQ,x)& =\  K(PAQ,x,x)\tag 2.3.1\\
& =\ \int_M K(P,x,y) A(y) K(Q,y,x)\, dy.
\endalign
$$
Hence $T$ is a smoothing operator  with integral kernel
$ K(P,x,y) \otimes K(Q,y,x)$. The same result holds if $Q=\phi Q\psi$.
By taking a suitable partition of 
unity we just need to prove Theorem 1.2.6 in the
case when there exists a coordinate chart $\Omega$ 
over which  the vector bundles $E_i$ are trivial, and
$$
\text{there exists }\phi\text{ supported in }\Omega\text{ with }
P=\phi P\phi, \quad Q=\phi Q\phi, \quad A=\phi A \phi.
\tag 2.3.2
$$ 
Now writing $P,A,Q$ as matrices
we have
$$
(PAQ)_{i\ell}\ =\ \sum_{j,k} P_{ij} A_{jk} Q_{k\ell}.
$$
If we prove the scalar case of the theorem we see that the general case follows 
immediately by considering the matrix entries.

We assume then that $P,Q:C^\infty(M)\to C^\infty(M)$ are odd-class pseudodifferential 
operators,
$A\in C^\infty(M)$ and we restrict attention to a coordinate chart
on which (2.3.2) holds.  Write $Q^t$ for the transpose of $Q$ with respect to the
standard density $dx$ on the coordinate chart. We consider separately the following cases:
\medbreak
\item{(1)} $P$ and $Q$ are smoothing.
\smallbreak\item{(2)} $P$ is purely singular and $Q$ is smoothing.
\smallbreak\item{(3)} $Q^t$ is purely singular and $P$ is smoothing.
\smallbreak\item{(4)} $P$ and $Q^t$ are purely singular with $\ord P+\ord Q<-n$.
\smallbreak\item{(5)} $P$ and $Q^t$ are purely singular and 
in local coordinates, $K(P,x,x-w)$ and
$K(Q^t,x,x-w)$ are homogeneous in $w$ (for $w$ small).
\medbreak\noindent 
The general case follows since
$$
P\ =\ P_\sing+P_\reg, \qquad\qquad Q^t=Q^t_{\hphantom{t} \sing}+Q^t_{\hphantom{t} \reg},
$$
and 
$$
P_\sing\ =\  P_1\ +\ \dots\ +\ P_j\ +\ P^{(j)},\qquad\qquad
Q_{\hphantom{t} \sing}^t\ =\  Q_1\ +\ \dots\ +\ Q_j\ +\ Q^{(j)}
$$
where $P^{(j)}, Q^{(j)}$ have low order and $K(P_k,x,x-w)$,
$K(Q_k,x,x-w)$ are homogeneous in $w$.

\medbreak
{\it Case} (1).
If $P$ and $Q$ are both smoothing then (2.3.1) holds and $T$  is a smoothing operator.

\medbreak
{\it Case} (2).
Similarly, when $P$ is purely singular and $Q$ is smoothing,
$PAQ$ is smoothing and 
$$\align
TA(x)& =\  K(PAQ,x,x)\\
& =\ \lan K(P,x,\cdot ),  A(\cdot) K(Q,\cdot,x)\ran\  =\ \lan K(P,x,\cdot ) K(Q,\cdot,x), A(\cdot)\ran.\endalign
$$
The product $ K(P,x,x-w ) K(Q,x-w,x)$ has singular parity in $w$ as $w\to0$, and 
the lowest homogeneity occurring in $w$ is at least $-n-\ord P$. Hence $T$ is an odd class pseudodifferential 
operator with $\ord T\leq \ord P$.

\medbreak
{\it Case} (3). follows in the same way as Case (2), or indeed can be deduced from Case (2) because
$$\align
TA(x)& =\ \KERD(PAQ,x)\\ & =\ K(PAQ,x,x)\ =\ K( Q^t A P^t,x,x)
\ =\ \KERD(Q^t A P^t,x).\endalign
$$
Here, $T$ is an odd class pseudodifferential operator with   $\ord T\leq\ord Q$.

\medbreak
{\it Case} (4). When $P$ and $Q$ are both purely singular and
\hbox{$\ord P+\ord Q <-n$,} the operator $PAQ$ is trace class and the function $K(P,x,y)K(Q,y,x)$
is integrable.  A simple continuity argument shows that (2.3.1) is still valid.
The kernel  $K(P,x,y )K(Q,y,x)$
has an expansion in homogeneous terms in $y-x$ in local coordinates and so $T$ is a classical
pseudodifferential operator.

\medbreak
{\it Case} (5).  Using  local coordinates to work in $\Bbb R^n$ and reducing the support of $A$ if necessary,
we assume that  $K(P,x,x-w)$ and $K(Q^t,x,x-w)$ are homogeneous distributions 
in $w$ for $w\in \Bbb R^n$.  Let $p(x,\xi)$ and $q^t(x,\xi)$ be the symbols of $P$ and $Q^t$ respectively,
which are homogeneous in $\xi$.  In a neighborhood of some point $x$, 
define the meromorphic families $P(z)$ and $Q(z)$ so that close to $x$ the symbols 
of $P(z)$ and $Q^t(z)$ are $p(z,x,\xi)$
and $q^t(z,x,\xi)$ respectively, where
$$
p(z,x,\xi)\ =\ p(x,\xi)|\xi|^{z/2}, \qquad\qquad q^t(z,x,\xi)\ =\ q^t(x,\xi)|\xi|^{z/2}. 
$$
As   shown in Section 2.1, this defines distributions which vary meromorphically with $z$ and 
have poles  lying in the set $\Bbb Z+1/2$.
We will show that 
$$
\KERD(PAQ,x)\ =\ K(P(z) A Q(z),x,x)\bigl|^\mer_{z=0}.
\tag 2.3.3
$$
Indeed,  let $\phi\in C^\infty(\Bbb R)$ be positive and even with
$$
\phi(r)\ =\ r, \qquad r>1, \qquad\qquad\qquad \phi(r)=1,\qquad |r|<1/2.
$$
Define $P_{1}(z),  Q^t_{1}(z)$ with symbols  $p_{1}(z,x,\xi), q^t_{1}(z,x,\xi)$ by
$$
p_{1}(z,x,\xi)\ =\ p(x,\xi)(\phi(|\xi|))^{z/2}, \qquad\qquad
q^t_{1}(z,x,\xi)\ =\ q^t(x,\xi)(\phi(|\xi|))^{z/2}.
$$
Then $P_1(z)$ and $Q_1(z)$  vary analytically in $z$ and so 
$$
\KERD(PAQ,x)\ =\ K(P_1(z) A Q_1(z),x,x)\bigl|^\mer_{z=0}.
$$
However,  $P(z)-P_1(z)$ and $Q(z)-Q_1(z)$ are smoothing operators depending
meromorphically on $z$ with poles in $\Bbb Z+1/2$, and $P(0)-P_1(0)=0$,\break  $Q(0)-Q_1(0)=0$.
Hence
$$
K((P(z)- P_1(z)) A Q_1(z),x,x)\bigl|^\mer_{z=0}\ =\ K((P(0)- P_1(0)) A Q_1(0),x,x)\ =\ 0,
$$
and similarly 
$$ K(P_1(z) A (Q(z)-Q_1(z)),x,x)\bigl|^\mer_{z=0}=0$$ and 
$$ K((P(z)- P_1(z)) A (Q(z)-Q_1(z)),x,x)\bigl|^\mer_{z=0}=0$$ so that  (2.3.3) follows.
Now for $\Re z <\!<0$, 
$$\align
K(P(z) A Q(z),x,x)& =\ \int_{\Bbb R^n} K(P(z),x,y)A(y) K(Q(z),y,x)\\ \noalign{\vskip5pt}
& =\ \int_{\Bbb R^n} K(P(z),x,y) K(Q^t(z),x,y)A(y). \endalign
$$
But $K(P(z),x,x-w) K(Q^t(z),x,x-w)$ is homogeneous in $w$ of degree\break  $-2n-\ord P-\ord Q-z$ and is
an even function of $w$. 
From the proof of Lemma 1.2.5 given in Section 2.1, we see that as a distribution, this product 
meromorphically continues and gives  a regular parity distribution at $z=0$, which is the kernel of 
a pseudodifferential operator of order $n+\ord P+\ord Q$.
This completes the proof of (5).   The symbol of $T$ in this case
is  $(\check p \check q^t)^\wedge(x,\xi)$, but $q^t(x,\xi)=q(x,-\xi)$ where $q$ is the principal symbol of $Q$.
\enddemo

\section{Variation formulas}

Throughout this section, $g_0$ is a metric on the $n$-dimensional manifold $M$ where 
$n$ is odd. $V_0$ is the total volume of $M$ in the metric $g_0$ and $d\mu_0$ is the 
canonical volume element  $g_0$, scaled to have total volume one. 
$\Pi_0$ is the $d\mu_0$-orthogonal projection 
of $L^2(M)$ onto the constants, and
$\Delta_0=\Delta_{g_0}$.   
The components of tensor fields are given with respect to a local orthonormal frame
for $g_0$, and $D_i$ denotes $g_0$-covariant differentiation in the direction of the 
$i^{\rm th}$ frame vector.  
Consider a deformation $g(t)$ of $g_0$ and write
$$ \matrix 
\noalign{\vskip4pt}
 g(t)= e^{2f(t)} B(t),\qquad 
\text{where}\qquad\hfill\det B_{ij}&\hskip-6pt =\, 1, \hfill\qquad \hfill f&\hskip-6pt =\, {\displaystyle\frac1{2n}}\log \det g_{ij}, 
\hfill
\cr \noalign{\vskip3pt}
 \hfill A_{ij}&\hskip-6pt =\, \dot B_{ij}(0),\hfill\qquad\hfill \phi&\hskip-6pt =\, \dot f(0).\hfill\endmatrix
\tag 3.0.1
$$
Write $V$ for the total volume of $M$ in the metric $g$ and $\Pi$ for the orthogonal 
projection onto the constants in the metric $g$.

\demo{{\rm 3.1.} First variation of $\detp\Delta$}  Our goal  is to compute the first variation of  the homogenized 
log determinant functional
$$
 F\ =\ \log \detp\Delta \ -\ \frac2n \log V,
$$
and use this to characterize the critical metrics.

Write 
$G(x,y)$ for the Green function for $\Delta_0$ and $\G(x,y)$ for the regular
part of the Green function, so that  $\G(x,y)\, d\mu_0$ is the Schwartz kernel for
$(\Delta_0^{-1})_\reg$.
With the notation of (3.0.1), set 
$$
E=e^{(n-2)f}B^{-1}.
$$
Then
$$
\matrix
  \Pi&\hskip-6pt =\,  (\Pi_0\,  e^{2f})^{-1} \Pi_0 e^{2f}, \hfill&\hfill \qquad \dot\Pi(0) &\hskip-6pt =\, 2\Pi_0\,
(\phi-(\Pi_0\phi)),\hfill\\
\noalign{\vskip5pt}
  \Delta &\hskip-6pt =\,  -e^{-nf} D_i E_{ij} D_j,\hfill&\hfill\qquad
\dot \Delta(0)&\hskip-6pt =\,  -n\phi\Delta_0-D_i \dot E_{ij}(0) D_j.\hfill
\endmatrix
$$
\enddemo

\demo{Proof of Theorem {\rm 1.3.6}}
$$
\align
 \frac{d}{dt}\log\det\Delta& =  \TR (\dot \Delta+\dot \Pi)(\Delta+\Pi)^{-1} 
  =  \TR \dot \Delta(\Delta^{-1}+\Pi)  +  \TR \dot \Pi (\Delta^{-1}+\Pi)\\
&  =  \TR \dot\Delta \Delta^{-1}  +  \TR \dot \Pi  \Pi 
  =  \TR \dot\Delta \Delta^{-1}. \\
\endalign
$$
The first equality follows from [Ok1], and the last from
the fact that  $\TR \dot \Pi\Pi =\frac12\trace (\Pi^2)^\cdot=\frac12(\trace \Pi)^\cdot=0$.
Now $V=(\Pi_0 e^{nf})V_0$, so that 
$$
\align
 \qquad \frac{dF}{dt}\biggl|_{t=0}\ & =\ \TR \dot\Delta_0\Delta_0^{-1}\ -\ 2(\Pi_0 \phi) \tag 3.1.1 \\
&=\ \TR \left(- n \phi \Delta_0- D_i \dot E_{ij}(0) D_j\right)
\Delta_0^{-1}\ -\ 2(\Pi_0\phi)\\ 
& =\ (n-2) \TR \phi\Pi_0
- \TR D_i \dot E_{ij}(0) D_j \Delta_0^{-1}  \\
&=\ (n-2) \int_M \phi\, d\mu_0
\ +\ \int_M  \dot E_{ij}(0)(x)D_{y,i} D_{x,j} \G(x,y)\bigl|_{y=x}\, d\mu_0  \\
&=\ \int_M  \dot E_{ij}(0)(x)\left(D_{y,i} D_{x,j} 
\G(x,y)\bigl|_{y=x}+\frac1n\delta_{ij}
\right)\, d\mu_0. 
\endalign
$$
Here, we have integrated by parts using the fact that for a 
tensor field $H(x,y)$,
$$
D_i \left(H(x,y)\bigl|_{y=x}\right) \ =\ 
D_{x,i}  H(x,y)\bigl|_{y=x}\ +\ D_{y,i} H(x,y)\bigl|_{y=x}.
\tag 3.1.2
$$ 
We have also used the fact that
$$
\dot E_{ij}(0)\ =\ (n-2)\phi\delta_{ij}\ -\ A_{ij}(0),
$$
so  that 
\smallbreak
\hfill ${\displaystyle
(n-2)\phi\ =\ \frac{\dot E_{kk}}n. 
}$\hfill
\enddemo
 
\phantom{downtown}

\demo{Proof of Corollary {\rm 1.3.7}} From Theorem 1.3.6,
we see that the metric $g_0$ is critical for $F$ 
if and only if (1.3.2) holds;  that is,
$$
D_{y,i} D_{x,j} \G(x,y))\bigl|_{y=x}\ +\ \frac{\delta_{ij}}n\ =\ 0.
\tag 3.1.3
$$
To see that this implies (1.3.3), we first note that
$$
\align
&\tag 3.1.4\\
\noalign{\vskip-18pt}
 (\Delta_{0,x} \G(x,y))\bigl|_{y=x}& =\ -1.  \\
\noalign{\vskip10pt}
&\tag 3.1.5\\
\noalign{\vskip-20pt}
\qquad \enspace
 D_{y,i} D_{x,j} \G(x,y))\bigl|_{y=x} & =\ \frac12 D_i D_j (\G(x,x))
-( D_{x,i} D_{x,j} \G(x,y)) \bigl|_{y=x}.  
  \\
\noalign{\vskip10pt}
&\tag 3.1.6\\
\noalign{\vskip-20pt}
 D_{y,i} D_{x,i} \G(x,y))\bigl|_{y=x} &=\ -\frac12 \Delta_0 (\G(x,x))-1.
  \\
\endalign
$$
Indeed  (3.1.4) follows since
$$
\Delta_{0,x} G(x,y)\ =\ -1,\qquad\qquad\qquad y\neq x.
$$
To see (3.1.5) note that 
$$
D_x \G(x,y)\bigl|_{y=x}\ =\ D_y \G(x,y)\bigl|_{y=x},
$$
which follows since $\G(x,y)$  with $\G(y,x)=\G(x,y)$.
By the chain rule,
$$
D (\G(x,x))\ =\ D_x \G(x,y)\bigl|_{y=x}\ +\ D_y 
\G(x,y)\bigl|_{y=x}.
$$
Putting these two facts together,
$$
D_x \G(x,y)\bigl|_{y=x}\ =\ D_y \G(x,y)\bigl|_{y=x}
\ =\ \frac12D (\G(x,x)).
\tag 3.1.7
$$
Applying the chain rule to (3.1.7) gives (3.1.5), and setting $i=j$ in (3.1.5)
gives (3.1.6).  The condition (3.1.3) now 
becomes
$$
 \frac12 D_i D_j (\G(x,x))
- (D_{x,i} D_{x,j} \G(x,y))\bigl|_{y=x}\ +\ \frac{\delta_{ij}}n\ =\ 0.
\tag 3.1.8
$$
Setting $i=j$ and using (3.1.4) we get
$$
\frac12\Delta_0(\G(x,x))\ =\ 0.
$$
Now, $\G(x,x)$ is constant, and  (3.1.8) becomes
$$
-D_{x,i} D_{x,j} \G(x,y))\bigl|_{y=x}\ +\  \frac{\delta_{ij}}n\ =\ 0,
$$
which is (1.3.3)(b).  Clearly (1.3.3) implies (1.3.4).  Now assume that
(1.3.4) holds; that is, $\G(x,x)$ is constant and
$$
D_{x,i} D_{x,j} \G(x,y))\bigl|_{y=x}\ =\  c\delta_{ij}.
$$
Setting $i=j$ and summing we get $c=1/n$ which implies (3.1.8) which implies (3.1.3).
\enddemo

 3.2. {\it First variation of $\det L$}.
In this section we calculate the first variation of $\log\det L$
and hence characterize the critical metrics for $\det L$.

\nonumproclaim{Lemma 3.2.1}  Let $g(t)$ be a deformation of the metric
$g_0$ and set $h=\dot g(0)${\rm .}  Then
$$
\frac{dS}{dt}\biggl|_{t=0}\ =\ h_{ij,ij}\ -\ h_{ii,jj}\ -\ R_{ij}h_{ij}
$$
where $R_{ij}$ is the Ricci tensor for $g_0${\rm .}
\endproclaim

A proof of this fact can be found in [LP, \S 8], or 
in Section 3.4 of this paper. 
Now if $g(t)$ is a conformal variation of $g_0$, and
$g(t)= e^{2f(t)} g_0$, then 
$$
L\ =\ e^{-(\frac{n}2+1)f} L_0 e^{(\frac{n}2-1)f}.
$$
It was  shown in [PR] that in odd dimensions,
$$
\det L\ =\ \det L_0.
$$
Indeed, this can be easily deduced using the canonical trace:
$$\align
\frac{d}{dt}\log\det L&=  \TR \dot L L^{-1}\\
& =  \TR \left( (-\frac{n}2-1) \dot f L+(\frac{n}2-1)L\dot f\right)L^{-1}
\ =\ -2\TR \dot f\, =\, 0.\endalign
$$
Because of this, it is sufficient to consider variations  $g(t)$ which
preserve the volume element, so that  $\dot g(0)=A$ where $A_{ii}=0$.
We have
$$
\frac{d}{dt}\biggl|_{t=0}(\log\det  L)
\ =\  \TR \dot L(0) L_0^{-1}.
$$
Lemma 3.2.1 gives
$$
\dot S(0)\ =\  A_{ij,ij}\ -\ R_{ij}A_{ij},\qquad 
\dot\Delta(0)\ =\ D_i A_{ij} D_j\ =\ A_{ij} D_i D_j + A_{ij,i} D_j.
\tag 3.2.1
$$

\vglue9pt

\demo{Proof of Theorem {\rm 1.3.8}}
$$
\align
 \frac{d}{dt}\biggl|_{t=0}\log\det L& =\ \TR \dot L(0)L_0^{-1}  \\ \noalign{\vskip5pt}
& =\ \TR\left(  \  A_{ij} D_i D_j\ +\ A_{ij,i} D_j  
\ +\ \al \left( A_{ij,ij}- R_{ij}A_{ij}\right) \right) L_0^{-1} \\ \noalign{\vskip5pt}
&=\ \int_M \Bigg( A_{ij}D_{x,i} D_{x,j}\G(x,y)\bigl|_{y=x}
  +\ A_{ij,i} D_{x,j}\G(x,y)\bigl|_{y=x}
\\ \noalign{\vskip5pt} &\hskip1.6in  +\ \al \left( A_{ij,ij}-R_{ij}A_{ij} \right)\G(x,x) \Bigg)
\,d\mu_0  \\\noalign{\vskip5pt}
&=\ \int_M \Bigg( A_{ij}(D_{x,i} D_{x,j}-\al R_{ij})\G(x,y)\bigl|_{y=x}
 \\ \noalign{\vskip5pt}
&\hskip1.5in  +\ \left( \al-\frac12 \right) A_{ij} D_i D_j  (\G(x,x)) \Bigg)
\,d\mu_0  \\ \noalign{\vskip5pt}
&=\ \int_M A_{ij}(x)  Q_{ij}(x) \, d\mu_0(x)
\endalign
$$
where  
$$
Q_{ij}(x)\ =\  (D_{x,i}D_{x,j}-\al R_{ij}) \G(x,y)\bigl|_{y=x}\ -\ \frac{n}{4(n-1)} D_iD_j (\G(x,x)).
\tag 3.2.2
$$
So $g_0$ is critical for $\det L$ if 
$$
Q_{ij}(x) \ =\  C(x)\delta_{ij}
$$
for some function $C(x)$ which can be computed by contracting over $i$ and $j$ and
using the fact that $L_x G(x,y)|_{y=x}=0$.  We get
$$
C(x)\ =\ \frac1{4(n-1)} \Delta_0 (\G(x,x)),
$$
so $g_0$ is critical for $\det L$ if and only if 
\medbreak
\hfill ${\displaystyle
(D_{x,i}D_{x,j}-\al R_{ij}) \G(x,y)\bigl|_{y=x}
\ =\ \frac{1}{4(n-1)}(n D_iD_j+\delta_{ij}\Delta_0 )(\G(x,x)).}$\quad \hfill \enddemo

\phantom{snow}
  
 For later use we note that  if $g_0$ is critical for $\det L$ then for any $Z\in T_{g_0}\M$, 
$$
\TR  (D_iZ_{ij} D_j-\al Z_{ij} R_{ij}+\alpha Z_{ij,ij}) L_0^{-1}
\ =\ -\frac{1}{4(n-1)}\TR Z_{ii,jj} L_0^{-1}.
\tag 3.2.3
$$

\vglue6pt

\demo{{\rm 3.3.} Second variation of $\detp\Delta$} 
\demo{Proof of Theorem {\rm 1.3.10}}
We want to calculate the Hessian of $F$ at a critical 
metric $g_0$, and we do so in  three steps.  First we calculate the conformal case 
$\Hess F(2\phi g_0, 2\phi g_0)$.  Then we calculate the cross-term
$\Hess F(2\phi g_0, A)$ when $A_{ii}=0$,
and finally we calculate the volume element preserving the
case $\Hess(A, A)$.  

Recall that
$$
\bar\phi\ =\ \phi-\int_M \phi\, d\mu_0,\qquad\qquad \bar\Delta_0^{-1}\ =\ \Delta_0^{-1}-Z(1)\Pi_0.
$$
We will start by showing that
$$
\multline
 \Hess F(2\phi g_0,2\phi g_0)= \  
\frac{(n+2)(n-2)}2   \int_M \bar\phi ^2\, d\mu_0
  \\ \noalign{\vskip6pt}-\ \frac{(n-2)^2 }4\TR \left( (\Delta_0 \phi)\bar{\Delta}_0^{-1}\right)^2.
\endmultline\tag 3.3.1$$
Set  
$$
g(t)=e^{2t\phi}g_0.
$$
Now 
$$
\Hess F(2\phi g_0,2\phi g_0 )
\ =\ \frac{d^2F}{dt^2}\biggl|_{t=0}
$$
and as in (3.1.1), 
$$
\frac{dF}{dt}\biggl|_{t=0}\ =\  \TR \dot \Delta \Delta^{-1}
\ -\ 2(\Pi \phi).
$$
Now 
$$
\align
 \Delta& =  -e^{-nt\phi} D_i e^{(n-2)t\phi} D_i, \\ \noalign{\vskip4pt}
 \dot \Delta(0)& =   n\phi D_i D_i \ -\ (n-2)D_i \phi D_i\\  \noalign{\vskip4pt}
& =  -\frac{n+2}2  \phi \Delta_0\ +\ \frac{n-2}2 \Delta_0\phi\ -\ \frac{n-2}2  (\Delta_0\phi),
\endalign
$$
where we temporarily use the slightly unwieldy convention that
$(\Delta_0 \phi)$ is the operator $\Delta_0$ applied to the function
$\phi$, while $\Delta_0\phi$ is the operator $\Delta_0$ composed with multiplication by $\phi$.
Since there is nothing special about $t=0$, we get
$$
\dot \Delta\ =\ -\frac{n+2}2  \phi \Delta\ +\ \frac{n-2}2 \Delta\phi\ -\ \frac{n-2}2  (\Delta\phi).
$$
Then
$$
\align
 \TR \dot \Delta \Delta^{-1}
& = \ \TR \left(-\frac{n+2}2  \phi \Delta
\ +\ \frac{n-2}2 \Delta\phi\ -\ \frac{n-2}2  (\Delta\phi)    \right) \Delta^{-1}  \\
&  =\  -2\TR \phi(I-\Pi)\ -\ \frac{n-2}2\TR   (\Delta\phi)   \Delta^{-1} \\
& = \  2\TR \phi \Pi\ -\ \frac{n-2}2 \TR (\Delta\phi) \Delta^{-1}.
\endalign
$$
Thus
$$
\frac{dF}{dt}\ =\ -\frac{n-2}2 \ \TR ( \Delta \phi) \Delta^{-1}.
\tag 3.3.2
$$
A similar formula was obtained by Richardson in [Ri]  using
more elaborate calculations involving the heat kernel instead of the canonical trace.
Now 
$$
\TR (\Delta\phi)\Pi\ =\ \int_M \Delta\phi\, d\mu\ =\ 0.
$$
Hence for any choice of $\beta$,
$$
\frac{dF}{dt}\ =\  -\frac{n-2}2 \TR (\Delta \phi)(\Delta+\beta\Pi)^{-1}.
\tag 3.3.3
$$
Let $Z(s)$ be the zeta function for $\Delta_0$ set 
$$
\beta\ =\ -\frac1{Z(1)}, \qquad\qquad \bar\Delta\ =\ \Delta+\beta \Pi,
 \qquad\qquad\bar\Delta_0\ =\ \Delta_0+\beta \Pi_0.
\tag 3.3.4
$$
(We assume that $Z(s)\neq 0$. In the case $Z(s)=0$, the following argument can be modified
by taking $\beta$ to be a constant and letting $|\beta|\to \infty$ at the end.)
Since $g_0$ is critical for  $F$ we have $\G(x,x)=Z(1)$ for all $x$, so that
$$
\bar\Delta_0^{-1}\ =\ \Delta_0^{-1}-Z(1)\Pi_0, \qquad 
K( (\bar\Delta^{-1})_\reg , x, x)\ =\  \G(x,x)-{Z(1)}\ =\ 0,
$$
and $\TR \psi\bar\Delta_0^{-1}=0$ for any function $\psi$.
Now
$$
 \Pi_0 \bar\Delta^{-1}_0 \ =\ \bar\Delta^{-1}_0 \Pi_0\ =\ -Z(1)\Pi_0,\qquad 
\Delta_0 \bar\Delta^{-1}_0\ =\ \bar\Delta^{-1}_0\Delta_0 \ =\ I-\Pi_0.
$$
Hence
$$
\align
&\tag 3.3.5 \\
\noalign{\vskip-8pt}
 \frac{d^2F}{dt^2}\biggl|_{t=0}
\ & =\ -\frac{n-2}2\left(  \TR (\dot \Delta(0) \phi)\bar\Delta_0^{-1}
\ +\  \TR ( \Delta_0 \phi)   ((\Delta+\beta\Pi)^{-1})^\cdot(0)\right)    \\
\ &  =\ -\frac{n-2}2 \TR ( \Delta_0 \phi)   ((\Delta+\beta\Pi)^{-1})^\cdot(0)     \\
& =\ \frac{n-2}2  \TR ( \Delta_0 \phi)  \bar\Delta^{-1}_0  (\dot \Delta_0+\beta\dot \Pi_0)\bar\Delta^{-1}_0    \\
& =\ \frac{n-2}2\left( 
\beta \TR ( \Delta_0 \phi) \bar\Delta^{-1}_0 \dot \Pi_0  \bar\Delta^{-1}_0 
\ +\ \TR ( \Delta_0 \phi) \bar\Delta^{-1}_0 \dot \Delta_0 \bar\Delta^{-1}_0 \right).
 \\
\endalign
$$
But
$$
\align
\noalign{\vskip2pt}
&\tag 3.3.6\\
\noalign{\vskip-18pt}
 \qquad \beta \TR ( \Delta_0 \phi) \bar\Delta^{-1}_0 \dot \Pi_0  \bar\Delta^{-1}_0
&=\ \beta\beta^{-1}\TR ( \Delta_0 \phi) \Pi_0 \dot \Pi_0   \bar\Delta^{-1}_0
 \\
&   =\  n\TR   ( \Delta_0 \phi) \Pi_0  \bar\phi \bar\Delta^{-1}_0 
\ =\ n\int\bar\phi^2\, d\mu_0,
\endalign
$$
while
$$
\align
&\hskip-32pt \TR ( \Delta_0 \phi) \bar\Delta^{-1}_0 \dot \Delta_0  \bar\Delta^{-1}_0\\  
 =\ & \TR ( \Delta_0 \phi) \bar\Delta^{-1}_0
\left(  -\frac{n+2}2  \phi \Delta_0\ +\ \frac{n-2}2 \Delta_0\phi\ -\ \frac{n-2}2  (\Delta_0\phi)\right)
\bar\Delta^{-1}_0  \\
 =\ & -\frac{n+2}2\TR ( \Delta_0 \phi) \bar\Delta^{-1}_0    \phi (I-\Pi_0)  
\ +\ \frac{n-2}2\TR ( \Delta_0 \phi) (I-\Pi_0 )\phi \bar\Delta^{-1}_0   \\
& \hskip2.2in  -\ \frac{n-2}2\TR ( \Delta_0 \phi) \bar\Delta^{-1}_0   (\Delta_0\phi)\bar\Delta^{-1}_0 \\
=\ & -2\TR ( \Delta_0 \phi)\, (I-\Pi_0 )\,\phi\, \bar\Delta^{-1}_0
\ -\  \frac{n-2}2\TR ( \Delta_0 \phi) \bar\Delta^{-1}_0   (\Delta_0\phi)\bar\Delta^{-1}_0 \\
  =\ & 2\TR ( \Delta_0 \phi)\, \Pi_0\, \phi\,  \bar\Delta^{-1}_0
\ -\  \frac{n-2}2\TR ( \Delta_0 \phi) \bar\Delta^{-1}_0   (\Delta_0\phi)\bar\Delta^{-1}_0 \\
 =\ &  2\TR   \int_M \bar\phi^2\, d\mu_0
\ -\  \frac{n-2}2\TR ( \Delta_0 \phi) \bar\Delta^{-1}_0   (\Delta_0\phi)\bar\Delta^{-1}_0.  
\endalign
$$
Substitution of this and (3.3.6) into (3.3.5) gives (3.3.1). 
 
Now we show that if $A_{ii}=0$  then
$$
\Hess F(2\phi g_0, A)
\ =\ -\frac{n-2}2 \TR \phi_{kk}\bar \Delta_0^{-1} D_i A_{ij}  D_j \bar\Delta_0^{-1}.
\tag 3.3.7
$$
Take a deformation
$$
g(s,t)= e^{2t\phi} B(s), \qquad \det(B_{ij})=1, \qquad  B^\prime_{ij}(0)=A_{ij},
$$
where the prime denotes differentiation with respect to $s$.  With the notation of (3.3.4),
$$
\frac{\partial F}{\partial t}\ =\ -\frac{n-2}2 \TR (\Delta\phi) \bar\Delta^{-1}.
$$
Then
$$
\align
\qquad\quad \frac{\partial^2 F}{\partial s\partial t}\biggl|_{s=t=0}
  =\ & -\frac{n-2}2 \TR (\Delta^\prime(0,0) \phi) \bar \Delta_0^{-1}\\
& \hskip1in -\ \frac{n-2}2 \TR (\Delta_0\phi)((\Delta+\beta\Pi)^{-1})^\prime(0,0)  \\
  =\ &  -\frac{n-2}2 \TR \phi_{kk}\bar\Delta_0^{-1} \Delta^\prime(0)\bar\Delta_0^{-1}
\\ =\ & -\frac{n-2}2 \TR \phi_{kk}\bar\Delta_0^{-1} D_i A_{ij} D_j\bar\Delta_0^{-1}.
\endalign
$$

 Finally, we show that 
$$
\Hess F(A, A)\ =\ -\frac1n \int_M A_{ij}^2\, d\mu_0
\ -\ \TR  D_i A_{ij}  D_j \bar\Delta^{-1} D_k A_{k\ell} D_\ell \bar\Delta^{-1} .
\tag 3.3.8
$$
Taking a deformation $g(t)$ of $g_0$ with $\det( g_{ij}(t))=1$ for all $t$ and
$\dot g(0)=A$, we also set $Z_{ij}=\partial_t^2(g^{-1})_{ij}(0)$ and note that
since $\det(g_{ij})=1$, 
$$
Z_{ii}\ =\ A_{ij}A_{ij}.
$$
 Now with the notation of (3.3.4) we have
$$
\frac{dF}{ds}\ =\ \TR \Delta^\prime \bar\Delta^{-1}
$$
so that
$$
\frac{d^2 F}{ds^2}\biggl|_{s=0}\ =\ \TR \Delta^\pprime(0) \bar\Delta^{-1}_0
\ -\ \TR \Delta^\prime(0) \bar\Delta^{-1}_0\Delta^\prime(0)\bar \Delta^{-1}_0.
\tag 3.3.9
$$
Now 
$$
\Delta^\prime(0) \ =\ D_i A_{ij} D_j,\qquad\qquad \Delta^\pprime(0)\ =\ -D_i Z_{ij} D_j,
$$
so that the second term on the right of (3.3.9) equals the second term on the right in (3.3.8).
Now since $g_0$ is critical,
$$
\align
\TR \bar\Delta^\pprime(0) \Delta^{-1}_0&  =\  -\TR D_i Z_{ij} D_j \bar\Delta_0^{-1}
\ =\ \int_M  Z_{ij} D_{y,i} D_{x,j} \G(x,y)\biggl|_{y=x}\, d\mu_0 \\
&  =\ -\frac1n \int_M  Z_{ij}\delta_{ij}\,d\mu_0
\ =\ -\frac1n \int_M A_{ij}A_{ij} \, d\mu_0.  \\
\noalign{\vskip-24pt}
\endalign      
$$
  \enddemo
\phantom{hold on}

\demo{{\rm 3.4.} Second variation  of $\det L$}
In this section we prove Theorem 1.3.11. First we remark that 
this theorem enables us to compute $\Hess \log\det L$ on the space $T_{g_0}\M$.  Indeed,
$$
T_{g_0}\M\ =\ (\conf g_0+\diff g_0)^\perp\ \oplus\ (\conf g_0 + \diff g_0).
$$
 Furthermore, $A\in (\conf g_0+\diff g_0)^\perp$ if and only if $A_{ii}=0$, $A_{ij,j}=0$,
and Theorem 1.3.11 gives $\Hess\log\det L$ under these conditions.  We claim that if  
$X\in \diff(g_0)+\conf(g_0)$ and $ Y\in T_{g_0}\M$, then
$\Hess\log\det L(X,Y)=0$.  Indeed,
$$
\Hess(X,Y)\ =\ \partial_s \partial_t F (g(s,t))|_{s=t=0}
\tag 3.4.1
$$  
where $g(s,t)$ is a 2-parameter deformation of $g_0$ with $\partial_s g(0,0)=X$
and $\partial_t g(0,0)=Y$.  Write $X=2\phi\delta_{ij}+X_1$ where
$X_1\in\diff(g_0)$.  It is not hard to construct  a deformation
$h(s,t)$ of $g_0$ with 
$h(s,0)\in\Diff(g_0)$ for all $s$, $\partial_s h(0,0)=X_1$ and $\partial_t h(0,0)=Y$.
Set $g(s,t)=e^{2s\phi}h(s,t)$.  Then $\partial_s g(0,0)=X$ and $\partial_t g(0,0)=Y$.
But then $\log\det L(g(s,t)=\log\det L(h(s,t))$ for all $s$ and $t$,
and so we can replace $g(s,t)$ by $h(s,t)$ on the right-hand side of (3.4.1), which then
vanishes because $h(s,0)$ is critical
for $\log\det L$ so that  $\partial_t \log\det L(h(s,t))|_{t=0}=0$.
\enddemo

\demo{Proof of Theorem {\rm 1.3.11}}
We take a deformation of $g_0$ with $\dot g_{ij}(0)=A_{ij}$ where $A_{ii}=0$ and $A_{ij,j}=0$, and 
set $(g^{-1}_{ij})^{\cdot\cdot}(0)=Z_{ij}$. 
By Lemma 3.2.1,  
$$
\align
&\hskip-30pt \Hess \log \det L(A,A)\\ \noalign{\vskip4pt}
   =\ &\frac{d^2}{dt^2}\biggl|_{t=0} \log\det L_{g(t)}
\ =\ \TR\ddot L(0) L_0^{-1}\ -\ \TR L_0^{-1} \dot L(0) L_0^{-1} \dot L(0)  \\ \noalign{\vskip4pt}
\ =\ &-\TR  D_i Z_{ij} D_j L_0^{-1}\ +\  \al \TR \ddot S(0) L_0^{-1}  \\ \noalign{\vskip4pt}
&\hskip.75in  -\TR   (D_i A_{ij}D_j-\al R_{ij}A_{ij}) L_0^{-1}
 (D_k A_{k\ell} D_\ell-\al R_{k\ell}A_{k\ell}) L_0^{-1}.
\endalign
$$
The last term in this final expression is the same as the last term in  (1.3.11).  
To equate the other terms, we need  the following.

\nonumproclaim{Lemma 3.4.1}  If $A_{ii}=0$ and $A_{ij,j}=0${\rm ,} then
$$
 \frac{d^2S}{dt^2}\ =\ Z_{ij} R_{ij}\ -\ Z_{ij,ij}\ +\  A_{ij,k} A_{ik,j}\ -\  \tfrac12 A_{ij,k} A_{ij,k}.
$$
\endproclaim

Using this Lemma, equation (3.2.3), and the fact that $Z_{ii}=A_{ij}A_{ij}$,  we see that  
$$
\align
&\hskip-4pt -\TR  D_i Z_{ij}   D_j L_0^{-1}\ +\  \al \TR \ddot S(0) L_0^{-1}  
\\ \noalign{\vskip4pt}
&\enspace =\, \TR \left( - D_i Z_{ij} D_j+\al  Z_{ij} R_{ij}\ -\ \al Z_{ij,ij}
\ +\  \al A_{ij,k} A_{ik,j}\ -\  \tfrac\al{2} A_{ij,k} A_{ij,k}\right) L_0^{-1}  \\ \noalign{\vskip4pt}
&\enspace =\, \TR   \left( \tfrac1{4(n-1)} (A_{ij} A_{ij})_{kk} \ +\ \al A_{ij,k} A_{ik,j}
\ -\   \tfrac\al{2} A_{ij,k} A_{ij,k}\right) L_0^{-1}.\\
\noalign{\vskip-18pt}
\endalign
$$
\enddemo
\phantom{phew}

\demo{Proof of Lemma {\rm 3.4.1}} 
We work in  coordinates. The metric $g(t)$ is a deformation of $g_0$.
We denote the covariant derivative in the metric $g(t)$ by $\nabla=\nabla(t)$, 
so that $\nabla(0)=D$.   If $\omega$ is some tensor field on $M$ then
$\nabla \omega$ is not tensorial in $\omega$, but $(\nabla-D)\omega$
is. 
Let $\Gamma_{jk}^i=\Gamma_{jk}^i(t)$ be the Christoffel symbol for $g(t)$, so that
$$
 \Gamma_{jk}^i\ =\ \frac{g^{i\ell}}2 \left( \partial_j g_{k\ell}
\ +\ \partial_k g_{j\ell}\ -\  \partial_\ell g_{jk} \right),
\qquad 
\nabla_{\partial_j} {\partial_k}\ =\ \Gamma_{jk}^i \partial_i. 
$$
Set
$$
\bar \Gamma_{jk}^i\ =\  \Gamma_{jk}^i(t)-\Gamma_{jk}^i(0).
$$
Then
$$
\bar \Gamma_{jk}^i\ =\  \frac{g^{i\ell}}2 \left( D_j  g_{k\ell}
\ +\ D_k g_{j\ell}\ -\  D_\ell g_{jk}\right).
$$
This can be verified at each point by working in normal coordinates for the metric $g_0$.

Now the   curvature tensor is
$$
R_{ijk}^\ell\ =\ \partial_i \Gamma_{jk}^\ell\ -\ \partial_j \Gamma_{ik}^\ell
\ +\    \Gamma_{ip}^\ell \Gamma_{jk}^p\ -\ \Gamma_{jp}^\ell \Gamma_{ik}^p,
$$
and  
$$
\bar  R_{ijk}^\ell\ =\ R_{ijk}^\ell(t)\ -\ R_{ijk}^\ell(0).
$$
Again working in normal coordinates for $g_0$, we find that 
$$
\bar R_{ijk}^\ell\ =\ D_i \bar \Gamma_{jk}^\ell\ -\ D_j \bar\Gamma_{ik}^\ell
\ +\ \bar\Gamma_{ip}^\ell \bar\Gamma_{jk}^p \ -\ \bar\Gamma_{jp}^\ell \bar\Gamma_{ik}^p .
$$
This enables us to compute the first and second variations of the Ricci curvature.
Indeed, we write
$$
\dot g_{ij}(0)=h_{ij}, \qquad  \ddot g_{ij}(0)=w_{ij}, 
\qquad   (g^{ij})^{\cdot\cdot}(0)=z^{ij},
$$
so  that
$$
z^{ij}=-w^{ij}+2h^{ik}h_k^j.
$$
Then
$$
\align
 \dot \Gamma_{jk}^i(0)& =\  \tfrac12 \left( D_j  h_k^i
\ +\ D_k h_j^i \ -\  D^i h_{jk}\right),  \\ \noalign{\vskip4pt}
 \dot \Gamma_{ik}^i(0)& =\  \tfrac12 D_k h_i^i,  \\ \noalign{\vskip4pt}
 \ddot \Gamma_{jk}^i(0)& =\  \tfrac12 \left( D_j  w_k^i
\ +\ D_k w_j^i \ -\  D^i w_{jk}\right)\ -\   h^{i\ell} \left( D_j  h_{k\ell}
\ +\ D_k h_{j\ell}\ -\  D_\ell h_{jk}\right),  \\ \noalign{\vskip4pt}
 \ddot \Gamma_{ik}^i(0)& =\  \tfrac12  D_k w_i^i \ -\ h^{i\ell} D_k h_{i\ell},  \\ \noalign{\vskip4pt}
  \dot R_{ijk}^i(0)& =\ D_i \dot \Gamma_{jk}^i\ -\ D_j \dot\Gamma_{ik}^i,  \\ \noalign{\vskip4pt}
  \ddot R_{ijk}^i(0)& =\ D_i \ddot \Gamma_{jk}^i\ -\ D_j \ddot\Gamma_{ik}^i
\ +\ 2\left(\dot\Gamma_{i\ell}^i \dot\Gamma_{jk}^\ell\ -\ \dot\Gamma_{j\ell}^i\dot\Gamma_{ik}^\ell
\right) ,
\endalign
$$ 
so  that
$$
\dot S\ =\ \partial_t ( g^{jk} R_{jk} )\ =\  -h^{jk} R_{jk}\ +\ g^{jk} \dot R_{jk}
\ =\  -h^{jk} R_{jk}\ +\ D_j D_k h^{jk}\ -\ D_j D^j h_k^k.
$$
This proves Lemma 3.2.1.  To calculate the second variation of the scalar curvature
we restrict to the case $h_i^i=0$ and  $D_j h^j_k=0$.
In this case $\dot\Gamma_{ik}^i=0$ and $\ddot \Gamma_{ik}^i=0$, so that
$$
\align
 \ddot S\ =\ \partial_t^2 ( g^{jk} R_{jk} ) =\  &z^{jk}R_{jk}
\ -\ 2h^{jk}\dot R_{jk}\ +\ g^{jk}\ddot R_{jk} 
\tag 3.4.2 \\ \noalign{\vskip4pt}
 =\ & z^{jk}R_{jk}
\ -\ 2h^{jk}\left(D_i \dot \Gamma_{jk}^i\ -\ D_j \dot\Gamma_{ik}^i   \right)      \\ \noalign{\vskip4pt}
& +\ g^{jk}\left(  D_i \ddot \Gamma_{jk}^i\ -\ D_j \ddot\Gamma_{ik}^i 
\ +\ 2\dot\Gamma_{i\ell}^i \dot\Gamma_{jk}^\ell\ -\ 2\dot\Gamma_{j\ell}^i\dot\Gamma_{ik}^\ell \right)  \\
\noalign{\vskip4pt}
 =\ & z^{jk}R_{jk}
\ -\ 2h^{jk} D_i \dot \Gamma_{jk}^i   
\ +\ g^{jk}\left(  D_i \ddot \Gamma_{jk}^i\ -\ 2\dot\Gamma_{j\ell}^i\dot\Gamma_{ik}^\ell \right) .
\endalign
$$
Now
$$
\align
 -2h^{jk} D_i \dot \Gamma_{jk}^i 
 =\ & -2h^{jk}D_i D_j  h_k^i\ +\ h^{jk} D_i D^i h_{jk},  \\ \noalign{\vskip4pt}
 g^{jk} D_i \ddot \Gamma_{jk}^i 
 =\ & D_i D_j  w^{ij}\ -\  \frac12 D_i D^i w_j^j  \\ \noalign{\vskip4pt}
  =\ &-D_i D_j z^{ij} \ +\ 2(D_i D_j  h^{i\ell})h_\ell^j
 \\ \noalign{\vskip4pt}
& +\ 2( D^j  h^{i\ell})D_i h_\ell^j  \ -\  \frac12 D_i D^i (h_{j\ell} h^{j\ell}),  \\ \noalign{\vskip4pt}
 -2g^{jk} \dot\Gamma_{j\ell}^i\dot\Gamma_{ik}^\ell 
 =\ & -\tfrac12 \left( D_j  h_\ell^i\ +\ D_\ell h_j^i \ -\  D^i h_{j\ell}\right)
  \left( D_i  h^{j\ell}\ +\ D^j h_i^\ell \ -\  D^\ell h_i^j  \right)   \\ \noalign{\vskip4pt}
 =\  &\tfrac12 D^i h^{j\ell} D_i  h_{j\ell}\ -\  (D^j  h^{i\ell})  D_i h_{j\ell} .
\endalign
$$
Substituting this into (3.4.2) gives 
\medbreak
\hfill ${\displaystyle 
\ddot S\ =\ z^{jk} R_{jk}\ -\ D_i D_j z^{ij}\ +\ (D^j  h_{i\ell})  D_i h_{j\ell}
\ -\ \tfrac12 (D^i h^{j\ell}) D_i  h_{j\ell}.
}$\phantom{swim} \hfill
\enddemo

\section{Critical points have finite index} 

In this section we prove Theorem 1.3.1.  As usual
let $g_0$ be a metric on the $n$-dimensional closed compact manifold
$M$ with $n$ odd,  write tensors with respect to a $g_0$-orthonormal frame for $TM$
and write $d\mu_0$ for the volume element for $g_0$ normalized so that $\mu_0(M)=1$.
Write $S^2_0T^*M$ for the bundle of symmetric $(0,2)$-tensors
on $M$ which are trace-free relative to $g_0$.

\nonumproclaim{Lemma 4.1.1} 
{\rm (a)} The operators  $P,Q:C^\infty(M)\to C^\infty(M)$ defined by
$$ \align
 (P\phi)(x)\,d\mu_0(x)& =\ \KERD (  \bar\Delta_0^{-1}\phi \bar\Delta_0^{-1}, \, x),  \\
(Q\phi)(x)\,d\mu_0(x)& =\ \KERD (  L^{-1} \phi L^{-1}, \, x),\endalign
$$
are pseudodifferential operators of order $n-4$.  $P$ and $Q$ have
the same principal symbol $p${\rm ,} where  $(-1)^{(n+1)/2}p(x,\xi)>0${\rm .}
\smallbreak

 {\rm (b)}  
The operator $R_0:C^\infty( S^2_0T^*M)\to C^\infty(M)$ defined by
$$
( R_0 A)(x)\, d\mu_0(x)\ =\  \KERD( \bar\Delta_0^{-1} D^2_{k\ell} A_{k\ell} \bar\Delta_0^{-1},\, x)
$$
is a pseudodifferential operator of order $n-2${\rm ,} and  
$R_0$ has the same principal symbol
as the operator 
$$
R_1\ =\ C^{(n)}  \Delta_0^{(n-4)/2}D^2_{k\ell},
$$
where  $C^{(n)}$ is a constant{\rm }, and $D^2_{k\ell}$ denotes the
divergence operator  $A_{k\ell}\to D^2_{k\ell} A_{k\ell}${\rm .}

\smallbreak {\rm (c)} The operators $U,V:C^\infty( S^2_0T^*M)\ \to \ C^\infty(S^2_0T^*M)$ 
defined by
$$
\align
& (UA(x))_{k\ell}\,d\mu_0(x)\ =\  \KERD( \bar\Delta_0^{-1} D^2_{ij} A_{ij} \bar\Delta_0^{-1} D^2_{k\ell},\, x), \\
& (VA(x))_{k\ell}\,d\mu_0(x)\ =\  \KERD ( L_0^{-1} D^2_{ij} A_{ij} L_0^{-1} D^2_{k\ell},\, x),
\endalign
$$
 are elliptic of order $n$.  $U$ and $V$ have the same principal symbol $u${\rm ,}
where  $(-1)^{(n+1)/2}u(x,\xi)>0${\rm .}
\endproclaim

\demo{{R}emark}
In general, the  homogeneous terms in the symbol expansions of $P$,
$Q$, $R_0$, $U$, $V$  down to but not including the
term of order zero, can be explicitly computed from $g_0,\,\phi,\, A_{ij}$.  To obtain
terms of order  $|\xi|^0$ and lower involves computing
the  Taylor expansion at the diagonal 
of the {\it regular} part of the Green function for $\bar\Delta_{0}$ or $L$.
In Section 5, we will consider the case when $g_0$ is the
standard metric on $S^n$ and we will compute $U$ explicitly for $n=3$ and  $V$ explicitly for all odd $n$.
\enddemo

\demo{Proof of Theorem {\rm 1.3.1}}  We will  see that there exists
an order $n$ pseudodifferential operator 
$Y:C^\infty( S^2T^*M)\to C^\infty(S^2T^*M)$ such that whenever 
$$
h_{ij}=2\phi\delta_{ij}+A_{ij},  \qquad \text{with }\qquad A_{ii}=0,\qquad A_{ij,j}=0, 
\tag 4.1.1
$$
we have
$$
\Hess F(h, h)
\ =\ \int_M h_{ij}(Yh)_{ij}.
\tag 4.1.2
$$
Moreover, the  principal symbol $y$ of $Y$
satisfies $(-1)^{(n+1)/2}y(x,\xi)>0$, and hence by
elliptic theory  $Y$ has at most finitely many nonnegative eigenvalues
if $n=4m+3$ and finitely many nonpositive eigenvalues in $n=4m+1$.

To see that (4.1.2) holds, we must consider the terms in  (1.3.10). 
The two main terms are 
$$
-\left(\frac{n-2}2\right)^2 \TR \phi_{ii} \bar\Delta_0^{-1}
\phi_{kk} \bar\Delta_0^{-1}\ =\ \int_M \phi (\Delta_0 P \Delta_0 \phi)\, d\mu_0,
\tag 4.1.3
$$
and 
$$
-\TR D_i A_{ij} D_j \bar\Delta_0^{-1} D_k A_{k\ell} D_\ell \bar\Delta_0^{-1}
\ =\ \int_M A_{ij} (UA)_{ij}d\mu_0.
\tag 4.1.4
$$
Here we have used Lemma 4.1.1.  Now the map
$$
h\to (\phi, A), \qquad \text{ where }\enspace h_{ij}=2\phi\delta_{ij}+A_{ij}, \quad A_{ii}=0,
$$
gives rise to the  splitting
$$
C^\infty(S^2T^*M)\ =\ LM\ \oplus\ C^\infty(S^2_0T^*M)
$$
where $LM$ is the trivial line bundle over $M$.  We can use this splitting to extend 
$P$ and $U$ to operators on $C^\infty( S^2T^*M)$ in the natural way, so that  the 
sum of (4.1.3) and (4.1.4) is
$$
\int_M \phi (\Delta_0 P \Delta_0 \phi)\ +\ \int_M A_{ij} (UA)_{ij}
\ =\ \int_M h_{ij} (Y_0h)_{ij}d\mu_0,
$$
where $Y_0$ is an elliptic pseudodifferential operator on $C^\infty(S^2T^*M)$
of order $n$ and its principal symbol $y$ in normal coordinates satisfies
$$
y_{ijk\ell}(x,\xi)h_{ij} h_{k\ell} \ =\ 
-(n-2)^2 p(x,\xi)|\xi|^4 (\phi(x))^2\ -\ u_{ijk\ell}(x,\xi)A_{ij}A_{k\ell},
$$
and $(-1)^{(n-1)/2} y(x,\xi)>0$.
The other terms in (1.3.10) correspond to operators of lower order.
In particular  when $A_{ij,j}=0$ we see that  $R_1 A=0$ and   
so setting $R=R_0-R_1$ we have  $R_0 A=R A$.  But 
$R$ is a pseudodifferential operator of order $n-3$.
The  cross term in (1.3.10) is
$$
-(n-2) \TR \phi_{ii} \bar\Delta_0^{-1} D_k A_{k\ell} D_\ell \bar\Delta_0^{-1} 
\ =\ (n-2) \int_M \phi \Delta_0 RA\, d\mu_0.
$$
The operator $\Delta_0 R$ has order $n-1$ which is less than the order of $Y_0$.
We conclude that (4.1.2) holds where $Y$ has
the same  principal symbol  as $Y_0$. 

Similarly by Lemma 4.1.1 and Lemma 3.4.1, 
$$
\Hess \log\det L(A, A)
\ =\ \int_M A_{ij}(ZA)_{ij}
$$
where $Z$ is an elliptic  pseudodifferential operator of order $n$,
whose principal symbol $z=-u$, and so
$Z$ has at most finitely many nonnegative eigenvalues
if $n=3$ $\mod 4$ and finitely many nonpositive eigenvalues in $n=1$ $\mod 4$. 
\enddemo

\demo{Proof of Lemma {\rm 4.1.1}}
We apply Theorem 1.2.6.  Work in normal coordinates $w$ about the point 
$x\in M$ and extend the orthonormal frame $\partial_{w_1},\dots, \partial_{w_n}$ 
at $x$
smoothly to a local orthonormal frame for $TM$, thus giving a 
local trivialization of $S^2_0T^*M$. 
The principal symbol of
$D^2_{k\ell}$ at $x$ is $-\xi_k \xi_\ell$, and the principal
symbol of $\bar\Delta^{-1}$ and $L^{-1}$ at $x$ is $|\xi|^{-2}$. 

The Fourier transform of $|w|^s$ on $\Bbb R^n$ is 
$$
C(s)|\xi|^{-n-s},  \qquad\qquad  C(s)=C_n(s)\ =\ 
\frac{ \pi^{n/2} 2^{n+s} \Gamma
\left(\frac{n+s}2\right)}{\Gamma(-s/2)}
\tag 4.1.5
$$
when $s\notin \{-n,2-n,\dots\}$;  see [Ta].
Set
$$
C^\prime \ =\   \frac{\Gamma(n/2)}
{ 2^{n+3} \pi^{ (n/2)-1} (n+1) \Gamma(n) }\ >\ 0.
\tag 4.1.6
$$

\medbreak (a)  By Theorem 1.2.6, $P$ and $Q$ are 
pseudodifferential operators and the principal symbol of each is 
$$
\align
 p(x,\xi)& =\ \bigl( (|\xi|^{-2})^\vee (|\xi|^{-2})^\vee \bigr)^\wedge
\ =\ \left(\frac1{2\pi}\right)^{2n} (C(-2))^2 \left(\
|w|^{2-n} |w|^{2-n}\right)^{\wedge}  \\
& =\ \left(\frac1{2\pi}\right)^{2n} (C(-2))^2 C(4-2n) |\xi|^{n-4}\\
& =\ (-1)^{(n+1)/2} 16(n-1)(n+1)C^\prime |\xi|^{n-4}.
\endalign
$$

\smallbreak (b)  By Theorem 1.2.6, $R_0$ is a pseudodifferential
operator.  Its principal symbol, 
$$
r(x,\xi): S^2_{0}T^*_xM\ \to\ \Bbb C
$$
can be written as a symmetric matrix $r_{k\ell}(x,\xi)$ so that
$$
r(x,\xi) A\ =\ r_{k\ell}(x,\xi) A_{k\ell}.
$$
Since $A$ is trace-free, we can calculate $r_{k\ell}(x,\xi)$ modulo terms of the form 
$u(\xi)\delta_{k\ell}$.
$$
\align
 r_{k\ell}(x,\xi)& =\ 
-\left(\  (|\xi|^{-2})^{\vee}\, (\xi_k \xi_\ell |\xi|^{-2})^{\vee} \right)^{\wedge}  \\
\allowdisplaybreak
& =\ \left(\frac1{2\pi}\right)^{2n} (C(-2))^2 \left(\
|w|^{2-n} \frac{\partial^2}{\partial w_k \partial w_{\ell}}|w|^{2-n}
\right)^{\wedge}  \\
\allowdisplaybreak
& =\ \left(\frac1{2\pi}\right)^{2n} (C(-2))^2  n(n-2)  \left(\
|w|^{-2n} w_k w_\ell \right)^{\wedge}  \\
\allowdisplaybreak
&=\ -\left(\frac1{2\pi}\right)^{2n} (C(-2))^2 C(-2n) n(n-2) 
\frac{\partial^2}{\partial \xi_k \partial \xi_{\ell}}|\xi|^{n}  \\
\allowdisplaybreak
&=\  -\left(\frac1{2\pi}\right)^{2n} (C(-2))^2 C(-2n) n^2  (n-2)^2\  |\xi|^{n-4}
\xi_k \xi_{\ell}.  
\endalign
$$
Clearly the operator
$$
R_1\ =\ -\left(\frac1{2\pi}\right)^{2n} (C(-2))^2 C(-2n) n^2  (n-2)^2
\  \Delta_0^{(n-4)/2}D_i D_j
$$
has the same principal symbol as $R_0$. 

\smallbreak (c)
 Finally, by Theorem 1.2.6, $U$ and $V$ are pseudodifferential operators
and    their principal symbol $u$  can be written as a tensor $u_{ijk\ell}$
symmetric in $i$ and $j$ and in $k$ and $\ell$,  with
$$
(u(x,\xi)A)_{ij}\ =\ u_{ijk\ell}(x,\xi) A_{k\ell},
$$
and this matrix can be calculated 
modulo terms of the form $\delta_{ij} v(\xi)+ w(\xi)\delta_{k\ell}$. 
$$
\align
 u_{ij k\ell}& =\ 
\left(\  (|\xi|^{-2}\xi_i\xi_j)^{\vee} \,
 (\xi_k \xi_\ell |\xi|^{-2})^{\vee} \right)^{\wedge}  \\ \noalign{\vskip4pt}
& =\ \left(\frac1{2\pi}\right)^{2n} (C(-2))^2 \left(\ 
\left( \frac{\partial^2}{\partial w_i \partial w_{j}}|w|^{2-n}\right)
 \left( \frac{\partial^2}{\partial w_k \partial w_{\ell}}|w|^{2-n}\right)
\ \right)^{\wedge}  \\ \noalign{\vskip4pt}
& =\ \left(\frac1{2\pi}\right)^{2n} (C(-2))^2 n^2 (n-2)^2 \left(\  |w|^{-4-2n}
  w_i  w_j w_k  w_{\ell} 
\ \right)^{\wedge}  \\ \noalign{\vskip4pt}
& =\ \left(\frac1{2\pi}\right)^{2n} (C(-2))^2 C(-2n-4) n^2 (n-2)^2 
 \frac{\partial^4}{\partial \xi_i\partial \xi_j\partial \xi_k \partial \xi_{\ell}}
|\xi|^{n+4}  \\ \noalign{\vskip4pt}
&=\   \frac{(-1)^{(n+1)/2}C^\prime}{(n+2)(n+4)}
\frac{\partial^4}{\partial \xi_i\partial \xi_j\partial \xi_k \partial \xi_{\ell}}
|\xi|^{n+4}  \ =\   (-1)^{(n+1)/2}C^\prime \tilde u_{ijk\ell}
\endalign
$$
where 
$$
\multline
\tilde u_{ijk\ell}\ =\ 
n(n-2) \, |\xi|^{n-4} \xi_i \xi_j \xi_k \xi_{\ell} \\   \noalign{\vskip4pt}  +\ 
n \, |\xi|^{n-2} 
\left( \xi_i  \xi_k \delta_{j\ell}  +  \xi_i  \xi_\ell \delta_{jk} 
+  \xi_j  \xi_\ell \delta_{ik}  +  \xi_j  \xi_k \delta_{i\ell} \right) \\  \noalign{\vskip4pt}   +\ |\xi|^{n} 
\left( \delta_{ik}   \delta_{j\ell}  +    \delta_{jk} \delta_{i\ell}\right)
 .
\endmultline
$$
The symbol $\tilde u$ is clearly elliptic and positive.  Indeed dropping the summation convention, we have
$$
\align
\sum_{i,j,k,\ell} A_{ij}\tilde u_{ijk\ell}A_{k\ell}    
 =\  &n(n-2)\left(\sum_{i,j} A_{ij} \xi_i \xi_j\right)^2  |\xi|^{n-4}  \\ \noalign{\vskip4pt}
& +\ 4n \sum_j 
\left(\sum_i A_{ij}\xi_i\right)^2|\xi|^{n-2}    
\ +\  2 \sum_{i,j} (A_{ij})^2   |\xi|^{n}.  \\
\noalign{\vskip-28pt} 
\endalign
$$
\enddemo

 \phantom{wind}

\section{Standard spheres which are extremal} 

 Formulas for the Hessian of the functionals $F=\log\detp\Delta-(2/n)\log V$ and $\log\det L$ at 
general critical points were proved in Sections 3.3 and 3.4.  In Section 4 it was shown that the
critical points always have finite index.  In this section we evaluate explicitly the formula for $\Hess F$ on
the $3$-sphere and for $\Hess \log\det L$ on all odd dimensional spheres, 
to show that these spheres  are local extremals.
By contrast, in Section 6 we will see that the $(4m+3)$-sphere with $m=1,2,\dots$,
is a saddle point for  $F$.  

In Section 5.1, we develop the general theory to  compute the Hessian of the
log determinant functionals on the standard spheres,
and in Sections 5.2 and 5.3 we apply this theory to $F$ and $\log\det L$  respectively.
Throughout we write $SS^n$ for the  unit sphere bundle of $S^n$, $g_0$ for the standard
metric on $S^n$,   $d\s_{n}$ for the standard volume form,
$V_n$ for the volume $\s_n(S^n)$ and $d\mu_n$ for the normalized volume form $d\s_n/V_n$. 
We write $S^2_0 T^*M$ for the bundle of symmetric $(0,2)$ tensors on $M$ which are trace-free 
relative to $g_0$, 
and  the components of all tensors are given with respect to a $g_0$-orthonormal frame.

\demo{{\rm 5.1.} Fourier series for singular kernels} \enddemo

\nonumproclaim{Lemma 5.1.1} Let $P, Q:C^\infty(S^n)\to C^\infty (S^n)$ be odd class pseudodifferential operators
which are functions of the  Laplacian on the standard $n$\/{\rm -}\/sphere{\rm ,} and 
suppose $A\in C^\infty(S^2_0T^*S^n)$  and $A_{ij,j}=0${\rm .} Then
$$
\KERD( P A_{ij} D^2_{ij}  Q,x)\ =\ 0.
$$
\endproclaim

\demo{Proof} Write $\Bbb T$ for the circle $\Bbb R/(2\pi\Bbb Z)$. 
It is convenient to introduce some spaces of functions on $\Bbb T\setminus 0$.
$$
\align
\F\ &=\ \{\ f\in C^\infty(\Bbb T\setminus 0)\ :\  f(x)=f(-x)\ \}.  
\tag 5.1.1  \\
\F_0\ &=\ \{\ f\in \F\ :\  f(\pi)=0\ \}. 
\endalign
$$
Define $T:\F\to\F_0$ by 
$$
Tf\ =\  f^{\pprime}(r)-\frac{\cos r}{\sin r} f^\prime(r)
\ =\ \sin r\frac{d}{dr} \frac1{\sin r} \frac{df}{dr}.
\tag 5.1.2
$$
The right inverse $T^{-1}$  is given by 
$$
T^{-1} f\ := \   I\left( \sin r\  I\left(\frac{f}{\sin r}\right)\right),
\qquad\text{where }  If(x)\ =\ \int_{\pi}^x f(t)\, dt,\quad
x\in (0, 2\pi).
$$
Write $r=r(x,y)$ for the distance  between points $x$ and $y$
measured on the standard $n$-sphere.
If $P$ and $Q$ are functions of the Laplacian 
then we can define $\Psi,\Phi\in\Cal F$ by 
$$
\Phi(r)=K(P,x,y), \qquad\qquad \Psi(r)=K(Q,x,y),  \qquad 0<r\leq \pi,
$$
and extend $\Phi$ and $\Psi$ to functions in $\Cal F$.
On the standard sphere, $D^2_{y,ij} \cos r=-\delta_{ij}\cos r$.  Here, the subscript $y$
indicates that we are taking the partial covariant derivative with respect to $y$.  We have
$$
D^2_{y,ij} \Psi(r)\ =\ \left( \Psi^{\pprime}(r)-\frac{\cos r}{\sin r} \Psi^\prime(r)
\right) r_{y,i} r_{y,j}\ +\ \frac{\cos r}{\sin r} \Psi^\prime(r)\delta_{ij}.
$$
Write   $a(x,y)= r_{y,i} r_{y,j} A_{ij}(y)$.       Since $A_{ii}=0$, we get 
$$ \align
\tag 5.1.3\\
\noalign{\vskip-16pt}
\qquad K(A_{ij}D^2_{ij} Q,y, x)& =\ A_{ij}(y) D^2_{y,ij} K(Q,y,x)\\
\noalign{\vskip4pt}
& =\  A_{ij}(y)( T\Psi (r) ) r_{y,i} r_{y,j} \ =\ T\Psi(r) a(x,y).
\endalign
$$ 
 First we consider the case when $P$ and $Q$ are pseudodifferential
operators which are functions of the standard Laplacian,   $\ord P<-n$ and $\ord Q<-n-2$.
 For the following argument we do not need to assume  that $P$ and $Q$ are 
odd class or integral order.     
$$
\align
&\hskip-8pt K( P A_{ij} D^2_{ij} Q, x,x)\\ \noalign{\vskip4pt}
& \ =\  \int_{S^n} K(P,x,y) \, K(A_{ij}D^2_{ij} Q,y,x)\, d\mu_n(y) \\ \noalign{\vskip4pt}
& \ =\ \frac1{V_n}\int_{S^n} \Phi(r) a(x,y)( T\Psi)(r) d\s_n(y)  \\ \noalign{\vskip4pt}  
& \ =\ \frac1{V_n}\int_{r=0}^\pi\int_{v\in SS^n_x}  a(x,\exp_x(rv))\Phi(r)( T\Psi)(r)\sin^{n-1}r 
 \, drd\s_{n-1}(v) \\ \noalign{\vskip4pt}
& \ =\ \frac1{V_n}\int_{r=0}^\pi\int_{v\in SS^n_x}  a(x,\exp_x(rv))\,T T^{-1}(\Phi(r)( T\Psi)(r)\sin^{n-1}r)  
\, drd\s_{n-1}(v)  \\ \noalign{\vskip4pt}
& \ =\ \int_{S^n}  A_{ij}(y) D^2_{y,ij} T^{-1}(\Phi(r)( T\Psi)(r)\sin^{n-1}r)  \, d\mu_n(y)  \\ \noalign{\vskip4pt}
& \ =\ \int_{S^n} (D^2_{y,ij} A_{ij}(y)) T^{-1}(\Phi(r)( T\Psi)(r)\sin^{n-1}r)  \, d\mu_n(y)\ =\ 0.  \\
\noalign{\vskip6pt}
\endalign
$$
The orders of $P$ and $Q$ were chosen low enough so that all distributions involved are
bounded functions and   integration by parts is valid.  The calculation fails when
$\ord P$ and $\ord Q$ are large, but the lemma is easily obtained by setting
$$
P_z= P (I+\Delta)^{z/4}, \qquad\qquad Q_z=Q (I+\Delta)^{z/4} ;
$$
thus applying the lemma when  $\Re z\ll 0$ gives 
$K(  P_z A_{ij} D^2_{ij} Q_z, x,x)  = 0$. Hence by analytic continuation,
$\KERD(P_0 A_{ij} D^2_{ij} Q_0,x)=0$. 
\enddemo

\phantom{run}

We will now derive a formula for $\TR( A_{ij} D^2_{ij} P B_{k\ell}D^2_{k\ell} Q)$
when $A, B\in C^\infty(S^2_0 T^*S^n)$.

 For $(x,v)\in SS^n$ and $r\in \Bbb R$,  let  $(x_r,v_r)$ be the
image of $(x,v)$ after  flowing distance $r$ along the geodesic
with initial point $x$ and initial direction $v$, so that
 $(x_r, v_r)= (\exp|_x(rv), \exp_*|_{(x,rv)}(v))$.
 For a  function $g$ on $\Bbb T$, define the  Fourier 
coefficients of $g$:
$$
\hat g_k\ =\ \frac1{2\pi} \int_0^{2\pi}  g(r)  e^{ikr} \, dr.
$$
This definition extends to distributions $g$ by duality.
 For a function $a(x,v)$ on $SS^3$,  define 
$$
\hat a_k (x,v)\ =\ 
\frac1{2\pi } \int_0^{2\pi}   a(x_s, v_s) e^{iks} \, dr.
\tag 5.1.4
$$

\nonumproclaim{Lemma 5.1.2} Suppose that $n$ is odd and $P$ and $Q$
are odd class pseudodifferential operators on $S^n$  which are functions of the
standard Laplacian{\rm .} Suppose that $P${\rm ,} $Q$ have order at most $-2${\rm .}
Let $A,B\in C^\infty(S^2_0 T^*S^n)$   and 
for $v\in T^*_x S^n$ write $a(x,v)=A_{ij}(x)v_i v_j${\rm ,}
$b(x,v)=B_{ij}(x)v_i v_j${\rm .}  Then 
$$\align
&\tag 5.1.5\\
\noalign{\vskip-30pt}
\phantom{ run }&\TR(   A_{ij} D^2_{ij} P B_{k\ell}  D^2_{k\ell} Q, x)\\ 
\noalign{\vskip6pt}
&\qquad\qquad =\ \frac\pi{V_n^2} \int_{SS^n}   \sum_{k=-\infty}^\infty  \hat f_k \hat a_k(x,v) \hat b_{-k}(x,v) 
 d\s_{n-1}(v) d\s_n(x)
\endalign
$$
where $f(r)$ is the even distribution  on $\Bbb T$ defined by
$$
 f(r)\ =\ T\Phi(r)\, T\Psi(r) \sin^{n-1} r, \qquad\qquad 0<r<2\pi,
\tag 5.1.6
$$
and regularized so that for $r$ close to zero{\rm ,}
$f(r)$ is the sum of a polyhomogeneous distribution with
regular parity  and a bounded function{\rm .}  There is an asymptotic expansion for $\hat f_k${\rm :}
$$
\hat f_k\ \sim\ c_n |k|^n \ +\ c_{n-2} |k|^{n-2}\ +\ \dots\ +\ c_1 |k|\ +\ o(1),
\qquad\qquad |k|\to \infty.
\tag 5.1.7
$$
\endproclaim

\demo{Proof} First, we consider the case when $P$ and $Q$ are pseudodifferential
operators which are functions of the standard Laplacian and $\ord P,\,\ord Q< -2$, 
$\ord P+ \ord Q<-n-4$.  For the following argument we do not need to assume  that $P$ and $Q$ are 
odd class or integral order.    The function $f(r)$ in (5.1.6) is integrable and
$$
\align
&\tag 5.1.8 \\
\noalign{\vskip-32pt}
\qquad\quad &  K(  A_{ij} D^2_{ij} P B_{k\ell} D^2_{k\ell} Q,x,x)  \\ \noalign{\vskip6pt}
&\hskip12pt = \ \frac1{V_n} \int_{S^n}  A_{ij}(x) (D^2_{x,ij}\Phi(r)) 
B_{k\ell}(y) (D^2_{y,k\ell} \Psi(r))  \, d\s_n(y)   \\ \noalign{\vskip6pt}
&\hskip12pt =\ \frac1{V_n}\int_{S^n}  A_{ij}(x) r_{x;i}r_{x;j}  B_{k\ell}(y) 
 r_{y;k} r_{y;l} T\Phi(r) T\Psi(r) \,  d\s_n(y)  \\ \noalign{\vskip6pt}
&\hskip12pt =\  \frac1{2V_n} \int_{SS^3_x}\int_0^{2\pi} a(x,v) b(x_r, v_r)
 T\Phi(r) T\Psi(r)   \sin^{n-1} r
\, dr d\s_{n-1}(v)  \\ \noalign{\vskip6pt}
&\hskip12pt = \  \frac1{2V_n} \int_{SS^3_x}\int_0^{2\pi} f(r) a(x,v) b(x_r, v_r)
\, dr d\s_{n-1}(v). 
\endalign
$$
Integrating over $x\in S^n$, we get  
$$
\align
&\tag 5.1.9 \\
\noalign{\vskip-12pt}
&\hskip-12pt \trace(    A_{ij} D^2_{ij} P B_{k\ell} D^2_{k\ell} Q,x,x)
 \\  \noalign{\vskip4pt}
&\quad =\  \frac1{2V_n^2} \int_{SS^3}\int_0^{2\pi} f(r) a(x,v) b(x_r, v_r)\, dr d\s_{n-1}(v)d\s_n(x) \\  \noalign{\vskip4pt}
&\quad= \ \frac1{4\pi V_n^2} \int_{SS^3}\int_0^{2\pi}\int_0^{2\pi} a(x_s,v_s) b(x_{s+r}, v_{s+r}) f(r)
\, ds dr d\s_{n-1}(v)d\s_n(x)  \\  \noalign{\vskip4pt}
&\quad = \   \frac{\pi}{V_n^2} \int_{SS^3} 
\sum_{k=-\infty}^\infty  \hat f_k \hat a_k(x,v) \hat b_{-k}(x,v)   d\s_{n-1}(v)d\s_n(x). 
\endalign
$$ 
Here, we have used the fact that the measure $ d\s_{n-1}(v)d\s_n(x)$
is invariant under the geodesic flow $(x,v)\to (x_s,v_s)$. 

To establish the general case we may suppose without loss of generality
that   $P$ and $Q$ have order $-2$, 
form the analytic family $P_z=P(I+\Delta)^{z/2}$ and set $\Phi_z(r)=K(P_z,x,y)$
and $f_z(r)=T\Phi_z(r)T\Psi(r)\sin^{n-1}r $.
Then   there exist  constants $p_j(z)$,  $u_j(z)$, 
$c_j(z)$, $ \tilde c_j(z)$ depending
analytically on $z$ when $\Re z<2$, such that as $r\downarrow 0$, 
$$
\multline
\phantom{T} \Phi_z(r)\sim\  \left( p_{2-n}(z) r^{2-n-z}\,+\, p_{4-n}(z) r^{4-n-z}\,+\,\dots\right)
\\     +\  \left( p_{0}(z)\,+\, p_{2}(z) r^{2}\,+\,\dots\right),  
 \endmultline
$$
$$\multline T\Phi_z(r) \sim\  \left( t_{-n}(z)
r^{-n-z}\,+\, t_{2-n}(z) r^{2-n-z}\,+\,\dots\right)
\\    +\  \left( t_{2}(z)r^2\,+\, t_{4}(z) r^{4}\,+\,\dots\right),   \endmultline
$$
$$ 
 T\Psi(r) \sim\  \left( u_{-n} r^{-n}\,+\, u_{2-n} r^{2-n}\,+\,\dots\right)
\ +\  \left( u_{2}r^2\,+\, u_{4} r^{4}\,+\,\dots\right),  \qquad\enspace
$$
 so that
$$\multline
 f_z(r) \sim\  \left( c_{-1-n}(z) r^{-1-n-z}\,+\, c_{1-n}(z) r^{1-n-z}\,+\,\dots\right)
\\    +\ \left( c_{1}(z) r^{1-z}\,+\,  c_{3}(z) r^{3-z}\,+\,\dots\right) \phantom{ r^{1-n-z}\,T\Phi(r)\Phi_z(r)}  \\ 
      +\  \left(\tilde c_0(z)\,+\,\tilde c_1(z) r\,+\, \tilde c_2(z) r^{2}\,+\,\dots\right).
\endmultline
$$  
The terms  $c_{j}(z) r^{j-z}$ with even $j<0$  analytically continue  to $z=0$ to give
regular parity distributions, and 
$$
 f_z(r)\ -\ \left( c_{-1-n}(z) r^{-1-n-z}\,+\, c_{1-n}(z) r^{1-n-z}\,+\,\dots \,+\, c_{-1}(z) r^{-1-z}\right)
$$
analytically continues to $z=0$ to give a bounded function of $r$. Hence as a distribution, $f_z(r)$ 
analytically continues to $z=0$ to give a regular parity distribution in $r$  plus an 
integrable function of $r$.  The Fourier coefficients of the function $f_0$ clearly satisfy (5.1.7).
Equation (5.1.8) holds when $z<\!<0$, and extends 
in the distributional sense to $z=0$.  \enddemo

 Finally we outline the procedure to calculate $\hat f_k$ from  the function $f$.

\demo{Definition {\rm 5.1.3}}
If $f\in \Cal F$,
we say $\al_k$ are  {\it Fourier coefficients}
for $f$ if the sequence $\al_k$ is {\it tempered\/}; that is, for some $N$,
$$
|\al_k|\ \leq \  C(1+|k|^N),
$$
and the distribution 
$$
\sum_k \al_k e^{-ikr}
$$
is equal to the  function $f(r)$ on $\Bbb T^n\setminus 0$.
\enddemo
  
\nonumproclaim{Lemma 5.1.4}
Suppose  $f$ is a function in $\Cal F$ and  $\al_k${\rm ,} $\gamma_k$
are two sets of Fourier coefficients  for $f${\rm .}  Then there
exists a polynomial $p(x)$ such that 
$$
\al_k\ =\ \gamma_k\ +\ p(k).
$$
\endproclaim

Our strategy to compute the Fourier coefficients of the distribution $f$
in (5.1.6) will be to compute  Fourier coefficients $\gamma_k$ for the  restriction of  
$f$ to $\Bbb T\setminus 0$, 
and then to adjust these 
by adding on a polynomial $p(k)$ to get coefficients $\al_k$ with the 
asymptotics given by  (5.1.7).
These are then necessarily  the Fourier coefficients of the distribution $f$. 

To calculate  Fourier coefficients of the function $f\in\Cal F$, we look for a 
harmonic function  $h(z)$ on the unit disc
$\D=\{z\in\Bbb C:|z|<1\}$ which extends continuously to $\overline\D\setminus 1$ 
and equals $f$ on the circle in the sense that
$$
f(r)=h(e^{ir}), \qquad\qquad r\in\Bbb T\setminus 0.
\tag 5.1.10
$$

\nonumproclaim{Lemma  5.1.5} If {\rm (5.1.10)} holds and 
$$
h(z)\ =\ \al_0\ +\ \sum_{k=1}^\infty \left( \al_k z^k
+\al_{-k} \overline z^k\right),\qquad\qquad |z|<1,
\tag 5.1.11
$$
where $\al_k$ are tempered{\rm ,} 
then $\al_k$ are {\it Fourier coefficients} for $f${\rm .}
\endproclaim

We omit the proofs of  Lemmas 5.1.4 and 5.1.5.  It is convenient 
to write 
$$
f\ \sim\ \al_0\ +\ \sum_{k=1}^\infty \left( \al_k z^k
+\al_{-k} \overline z^k\right) 
\tag 5.1.12
$$
when (5.1.10) and (5.1.11) hold.

\demo{{\rm 5.2.} The Hessian of $F$ at $S^3$}
Theorem 1.3.10 gives the Hessian of the function $F$ at a critical metric $g_0$.
In dimension $3$, 
$$
\align \noalign{\vskip4pt}
& \tag{5.2.1} \\
\noalign{\vskip-32pt} \qquad \quad &
  \Hess F( 2\phi g_0+A  ,2\phi g_0+A)
 \\&\hskip16pt =\  -\ \frac13\int_{S^3} A_{ij} A_{ij}\,d\mu_0 
\ +\ \frac52  \int_{S^3} \overline\phi^2\, d\mu_0
\ -\ \frac14 \TR \phi_{ii} \bar \Delta_0^{-1}  
\phi_{jj} \bar \Delta_0^{-1}
  \\
&\hskip32pt   -\  \TR   A_{ij} D^2_{ij}    \bar \Delta_0^{-1}  
\phi_{jj} \bar \Delta_0^{-1}  
\ -\  \TR A_{ij} D^2_{ij}  \bar \Delta_0^{-1}  
A_{k\ell} D^2_{k\ell} \bar \Delta_0^{-1} ,
\endalign
$$
where $A\in S^2_0T^*S^n$ with $A_{ij,j}=0$.
Our aim is to show that this is strictly negative unless $A_{ij}=0$ and $\phi\delta_{ij}\in\diff(g_0)$.
The term 
$$
\TR   A_{ij} D^2_{ij}    \bar \Delta_0^{-1}  
\phi_{jj} \bar \Delta_0^{-1} 
$$
vanishes by Lemma 5.1.1.  
Richardson [Ri] studied the case of conformal deformations and his results
show that when $A_{ij}=0$,  $\Hess F\leq 0$   with strict inequality when
$\phi\delta_{ij}\notin\diff(g_0)$. Hence we only need to understand the case $\phi=0$. 
 For $n>0$, set
$$ 
\align
I(k)\ & =\ \cases 1 \quad & k\text{ even,}
 \\
0  & k\text{ odd.} \endcases\tag 5.2.2 \\
E(k)\  &  =\ \sum_{j=1}^k \frac{I(k-j)}{j^2}
\ +\ \frac{\pi^2 I(k)}{12}.
\endalign
$$

 \enddemo

\nonumproclaim{Lemma 5.2.1} For the standard metric on $S^3${\rm ,} suppose that $$A\in  C^\infty(S^2_0T^*S^2).$$
Write $a(x,v)=A_{ij}(x)v_i v_j${\rm .} Then
$$
\TR   A_{ij} D^2_{ij}  \bar \Delta_0^{-1} 
A_{k\ell}  D^2_{k\ell} \bar \Delta_0^{-1} 
=\ \frac1{4\pi^3} \int_{SS^3}   \sum_{k=-\infty}^\infty  \al_k |\hat a_k(x,v)|^2 
 d\s_2(v) d\s_3(x)
\tag 5.2.3
$$
where $\al_k>0$ for all $k${\rm ,} so that  {\rm (5.2.3)} is positive unless
$A=0${\rm .}
\endproclaim

\demo{Proof} Let $G(r)$ be  Green's function where $r$ denotes the distance between $x$ and $y$; that is
$$
G(r)\ =\ K(\bar \Delta_0^{-1},x,y) \ =\ \frac{(\pi-r)\cos r}{2\sin r} +\ \frac12.
\tag 5.2.4
$$
By Lemma 5.1.2, (5.2.3) holds where $\al_k$ are the Fourier coefficients
of  $f(r)= (TG(r))^2 \sin^2 r$.  We will follow the strategy outlined
below Lemma 5.1.4 to compute these Fourier coefficients. Now
$$
G^\prime(r)  =  \frac{ -\cos r}{2\sin r}
+\frac{(r-\pi)}{2\sin^2 r}, 
$$
and since
$$
G^\pprime(r)  +\frac{2\cos r}{\sin r}G^\prime(r)\ =\ 1,
$$
we see that
$$
\align
TG(r)\ =\ & G^\pprime(r)  -\frac{\cos r}{\sin r}G^\prime(r) =
1  -\frac{3\cos r}{\sin r}G^\prime(r) \\
=\ & 1 +  \frac{ 3\cos^2 r}{2\sin^2 r}
 - \frac{3(r-\pi)\cos r}{2\sin^3 r}  
\  =\  \frac{-3(r-\pi)\cos r}{2\sin^3 r}
 + \frac{ 3}{2\sin^2 r}
 -  \frac{1}2.
\endalign
$$
Hence
$$
\align
\quad f(r)& =\  (TG(r))^2\sin^2 r \ =\ \biggl( \frac{-3(r-\pi) \cos r}{2\sin^2 r}
+ \frac{ 3}{2\sin r}- \frac{\sin r}2 \biggr)^2   \tag 5.2.5
\\ \noalign{\vskip4pt}
&=\  \frac{9(r-\pi)^2 \cos^2 r}{4\sin^4 r}
 -\frac{9(r-\pi) \cos r}{2\sin^3 r} +\frac{ 9}{4\sin^2 r}\\ \noalign{\vskip4pt}
&\hskip2.5in+\frac{3(r-\pi)\cos r}{2\sin r} 
-\frac32+ \frac{\sin^2 r}4    \\ \noalign{\vskip4pt}
&=\ \frac38 \left( \frac{d^2}{dr^2}-2\right) \left( \frac{(r-\pi)^2}{\sin^2 r}
 - 2 \frac{(r-\pi)\cos r}{\sin r} \right) 
\ -\ \frac32\ +\  \frac{\sin^2 r}4 \\ \noalign{\vskip4pt}
& =\ -\frac38 \left( \frac{d^2}{dr^2}-2\right) \frac{d}{dr} \frac{(r-\pi)^2\cos r}{\sin r}
\ -\ \frac32\ +\  \frac{\sin^2 r}4.
\endalign
$$
\medskip
We are now going to compute  Fourier 
coefficients for $f$ restricted to $\Bbb T\setminus 0$.   To do this we will find a harmonic function 
on the disc having $f$ as its boundary values.   We  rely on the simple
fact that although products of harmonic functions are
not in general harmonic, the product of two
analytic functions or of two conjugate analytic functions
is harmonic.  Define
$$
p(z)\ =\ -i\frac{1+z^2}{1-z^2}\ =\ -i\left( 1+ 2\sum_{j=1}^\infty  z^{2j}\right),\qquad 
h(z)\ =\ \frac{\pi^2}6+\sum_{k=1}^\infty\frac2{k^2}z^k.
$$
Setting $z=e^{ir}$ we have
$$
\frac{\cos r}{\sin r}\ =\  i \frac{z+\frac1z}{z-\frac1z}\ =\ -i\frac{1+z^2}{1-z^2}
\ =\ p(z); \qquad \text{ also, }\qquad\frac{\cos r}{\sin r}\ =\  \overline p(z).
$$
Calculation  of Fourier coefficients gives
$$
(r-\pi)^2\ =\ \frac{\pi^2}3\ +\ \sum_{k=1}^\infty
\frac2{k^2}(z^k+\overline z^k)\ =\ h(z)+\overline h(z).
$$
Now
$$
p(z) h(z)\ =\  -i\left( 1+ 2\sum_{j=2}^\infty I(j) z^{j}\right)
\left( \frac{\pi^2}6\ +\ \sum_{k=1}^\infty\frac2{k^2}z^k\right) 
\ =\ \frac{-i\pi^2}6 \ +\ \sum_{k=1}^\infty \beta_k z^k 
$$
where
$$
\beta_k\ =\ \cases -2i \qquad & k=1, \\
-i \left(   4E(k-2)\ +\  \frac{2}{k^2}  \right)  \qquad & k>1.  \endcases
$$
Hence
$$
\frac{(r-\pi)^2\cos r}{\sin r}\ =\ 
p(z) h(z)\ +\ \bar p(z) \bar h(z)\ =\ 
\frac{-i\pi^2}{3}  \ +\ \sum_{k=1}^\infty \beta_k (z^k-\bar z^k).
$$
 From (5.2.5), we see that
$$
f(r)\ \sim\ \gamma_0\ +\ \sum_{k=1}^\infty \gamma_k (z^k+\bar z^k)
$$
where
$$ \align
\gamma_k& =\ \cases -\frac32+\frac18 & k=0 \\
\frac{3}{8}(k^2+2)ik\beta_k   & k=1 \text{ (or } k>2) \\
\frac{18i\beta_2 }4  \ -\ \frac1{16}            & k=2  
\endcases
\\ \noalign{\vskip5pt}
& =\  \cases  -\frac{11}8   & k=0 \\
\frac{9}{4}    \ \  & k=1,  \\
 \frac{3\pi^2}2 \ +\  \frac94   \ -\ \frac1{16}            & k=2 \\
\frac{3}{2}(k^3+2k)E(k-2)\ +\  \frac{3k}{4}\ +\ \frac{3}{2k}     \ \  &  k>2.
\endcases \endalign
$$
Now
$$
E(k-2) \ =\  \frac{\pi^2}8
  - \frac1{2k}-\frac1{2k^2}-   \frac1{3k^3}
+ O(k^{-4})
$$
as $k\to\infty$, so  that
$$
\gamma_k\ =\ \frac{3\pi^2}{16} k^3\ -\ \frac34 k^2
\ +\ \frac{3\pi^2}8 k\ -\ 2\ +\ O(\frac1k)
$$
as $k\to\infty$.
The Fourier coefficients $\al_k=\hat f_k$ are therefore given
by 
$\al_k=\gamma_k+(3/4)k^2+2$ and so
$$
\al_k\ =\ \frac{3\pi^2}{16} k^3
\ +\ \frac{3\pi^2}8 k\ +\ O(\frac1k)
$$
as $k\to\infty$.  This gives 
$$
\al_k\ =\ \cases  \frac58 & k=0 \\
        5   &  k=1  \\
                \frac{3\pi^2}2+7+\frac{3}{16}  &  k=2 \\ 
\frac32 (k^3+2k)E(k-2)
 + \frac34\left(k^2+k+\frac2k\right)+2  &  k>2 
            \endcases
$$
so that $\al_k>0$ for all $k>1$.   \enddemo

\phantom{lunch}

 5.3. {\it The Hessian of $\log\det L$ at $S^n$}.
Let $g_0$ be the standard metric on $S^n$ where $n$ is odd.
As usual $r$ denotes the distance between points $x,y\in S^n$.  Green's function for $L_0$ is
$$
G(r)\ =\ K(L_0^{-1},x,y)\ =\ \frac{C_n }{\sin^{n-2}(r/2)},\qquad\qquad C_n\ =\ \frac{1}{(n-2)V_{n-1}}.
$$
Since this has singular parity as $r\downarrow 0$,  $(L_0^{-1})_\reg=0$. 
 Furthermore, $R_{ij}=2\delta_{ij}$. From this we see that (1.3.11) reduces to 
$$
\Hess \log\det L(A,A)\ =\ 
\ - \TR  A_{ij} D^2_{ij} L_0^{-1} A_{k\ell} D^2_{k\ell} L_0^{-1} 
\tag 5.3.1
$$
when $A\in S^2_0 T^*S^n$ with $A_{ij,j}=0$.

\nonumproclaim{Lemma 5.3.1} Suppose that $A\in C^\infty(S^2_0T^*S^2)${\rm ,}
where the trace is with respect to the
standard metric{\rm .}  Then
when $a(x,v)=A_{ij}(x)v_i v_j${\rm ,} 
$$ \align
\noalign{\vskip2pt}
&\tag 5.3.2\\
\noalign{\vskip-31pt}
&\qquad \TR   A_{ij} D^2_{ij} L_0^{-1} A_{k\ell} D^2_{k\ell} L_0^{-1} \\
&\hskip60pt  =\ (-1)^{(n+1)/2}\frac{\pi}{V_n^2} \int_{SS^3}   \sum_{k=0}^\infty  
\al_k |\hat a_k(x,v)|^2 \,d\s_{n-1}(v) d\s_n(x), 
\endalign
$$
where  $\al_k>0$ for all $k${\rm ,} so that {\rm (5.3.1)} is nonzero unless
$A_{k\ell}=0${\rm ,} and has the same sign as $(-1)^{(n-1)/2}${\rm .}
\endproclaim

\demo{Proof} 
By Lemma 5.1.2, (5.3.2) holds where $\al_k$ are the Fourier coefficients
of  $(-1)^{(n+1)/2}f(r)= (-1)^{(n+1)/2}(TG(r))^2 \sin^2 r$.  We will follow the strategy outlined
below Lemma 5.1.4 to compute these Fourier coefficients.
$$
TG(r)\ =\ \sin r\frac{d}{dr}\frac1{\sin r}\frac{d}{dr}G(r)\ =\ 
C_n^\prime \frac{\cos^2 (r/2)}{\sin^n(r/2)} , 
$$
and
$$
f(r)\ =\ \sin^{n-1} r (TG(r))^2 
\ =\ C_n^\pprime\cos^2(r/2)\left( \frac{\cos (r/2)}{\sin (r/2)}\right)^{n+1},  \qquad C_n^\pprime>0.
\tag 5.3.3
$$
Now if $z=e^{ir}$,
$$
(-1)^{(n+1)/2}\left( \frac{\cos (r/2)}{\sin (r/2)}\right)^{n+1}
\ =\  \left( \frac{1+z}{1-z}\right)^{n+1}.
$$

\nonumproclaim{Lemma 5.3.2}  For fixed $n${\rm ,} 
let $p(k)$ be the coefficients defined by 
$$
\left( \frac{1+z}{1-z}\right)^{n+1}\ =\ \sum_{k=0}^\infty p(k) z^k.
$$
Then $p(k)>0$ for all $k${\rm .} Moreover for $k>n+1${\rm ,} 
$$
p(k)=a_n k^n + a_{n-2} k^{n-2} +\dots + a_1 k, 
 \qquad \text{ where } a_j>0\text{ for all }j.
\tag 5.3.4
$$
\endproclaim

\demo{Proof}  
$$
\align
 \left( \frac{1}{1-z}\right)^{n+1}& =\ \sum_{k=0}^\infty
\frac{(k+1)(k+2)\dots(k+n)}{n!} z^k,  \\
 (1+z)^{n+1}& =\ \sum_{\ell=0}^{n+1} \left( \matrix n+1 \\ \ell\endmatrix\right) z^\ell.
\endalign
$$
Since the coefficients of both of these power series are positive, so are the coefficients $p(k)$
in the product.  Moreover when  $k>n+1$, $p(k)$ is equal to
$$
\align
&\hskip-36pt \frac1{n!} \sum_{\ell=0}^{n+1} \left( \matrix n+1 \\ \ell \endmatrix \right)
\prod_{j= 1}^{n} (k+j-\ell)  \\
& \ =\ \frac1{n!} \sum_{\ell=0}^{(n+1)/2} \left( \matrix n+1 \\ \ell \endmatrix \right)
\left( \prod_{j= 1-\ell}^{n-\ell} (k+j)
\ +\    \prod_{j= 1-\ell}^{n-\ell} (k-j)      \right),  \\
\endalign
$$
but since $n$ is odd,
$$
\prod_{j= 1-\ell}^{n-\ell} (k+j)
\ +\    \prod_{j= 1-\ell}^{n-\ell} (k-j)   
$$
is a polynomial having the form given in (5.3.4). 
\enddemo 

Now we see that when $z=e^{ir}$, 
$$
\cos ^2(r/2)\ =\  \left(\frac{z^{1/2}+z^{-1/2}}{2}\right)^2
\ =\ \frac14 (z+2+z^{-1}),
$$
so  that
$$
\cos ^2(r/2)\ =\  \frac14 (z+1)\ +\ \frac14 (1+\bar z)
$$
and
$$
\align
&\hskip-40pt (-1)^{(n+1)/2} \cos ^2(r/2)  \left( \frac{\cos (r/2)}{\sin (r/2)}\right)^{n+1} \\  \sim\ &
 \frac14 (z+1)\sum_{k=0}^\infty p(k) z^k
\ +\  \frac14 (1+\bar z)\sum_{k=0}^\infty p(k) \bar z^k  \\
 =\ &\frac14 \left( p(0)\ +\ \sum_{k=1}^\infty (p(k)+p(k-1)) z^k
\right)
\\&\hskip1.5in  +\ \frac14 \left( p(0)\ +\ \sum_{k=1}^\infty (p(k)+p(k-1)) \bar z^k
\right)  \\
  =\ & c(0)\ +\ \sum_{k=1}^\infty c(k) (z^k + \bar z^k).
\endalign
$$
Now let us examine the coefficient $c(k)$ for large $k$.
$$
\align
c(k)\ & =\ \frac14\left( p(k)\ +\ p(k-1) \right)\\ 
&=\ \frac14 \Big( a_n k^n + a_{n-2} k^{n-2} +\dots + a_1 k
\\
&\qquad\quad\qquad +\ a_n(k-1)^n + a_{n-2} (k-1)^{n-2} +\dots + a_1 (k-1) \Big) \\
&=\  b_n k^n + b_{n-1} k^{n-1} + \dots + b_0
\endalign
$$
where 
$$
\cases b_n >0 \qquad & n \text{ odd,} \\
b_n <0 \qquad & n \text{ even.}
\endcases
$$
Now for every $k\geq 0$, set
$$
e(k)\ =\ - b_{n-1} k^{n-1} - b_{n-3} k^{n-3} - \dots - b_0
$$
so that  $e(k)$ is positive.  Then
$$
\align
(-1)^{(n+1)/2} \cos ^2(\th/2)\left( \frac{\cos (r/2)}{\sin (r/2)}\right)^{n+1} 
 \sim\  & c(0)+e(0)\\
&+\sum_{k=1}^\infty (c(k)+e(k)) (z^k + \bar z^k)
\endalign
$$
and $(-1)^{(n+1)/2}\hat f_k\ =\ C_n^\pprime(c(k)+e(k))>0$.  \hfill\qed 
\enddemo

\section{Saddle points for $\detp\Delta$}

In this section we  show that the $(4m+3)$-sphere with $m=1,2,\dots$,
is a saddle point for  $F$ under  conformal deformations, by proving Lemma 1.3.12 which 
states that   if $n>3$ is odd and $\phi\in\Cal H_{2}$, then 
$$
\multline
\frac{(n+2)(n-2)}{2}  \int_{S^n} \bar\phi^2\, d\mu_n 
\ -\ \frac{(n-2)^2}4 \TR  \phi_{ii} \bar{\Delta}_0^{-1}  \phi_{kk}
\bar{\Delta}_0^{-1}\\ =\ \frac{ 2(n-2)\left(-2 + (n-2)(n+1)\sum_{j=2}^{n-1}\frac1j\right)}
{(n-3)n}\int_{S^n} \phi^2\, d\mu_n.
\endmultline
\tag 6.1.1
$$
The left side of (6.1.1) is $\Hess F(2\phi g_0,2\phi g_0)$.
The space $\Cal H_k$ has dimension
$$
d_{k}\ =\ \frac{ (k+1)(k+2)\dots (k+n-2) (2k+n-1)}{(n-1)!}
$$
and the eigenvalue of $\Delta$ on $\Cal H_k$ is
$$
\la_{k}=k(k+n-1).
$$
Write $Z_n(s)$ for the zeta function 
$$
Z(s)\ =\ \sum_{j=1}^\infty d_{j} \la_{j}^{-s}
$$
defined via analytic continuation for $\Re s\leq n/2$.

Define $\tilde \la_k^{-1}$ to be the eigenvalues of $\bar\Delta_0^{-1}$, so  that
$$
\frac1{\tilde\la_{k}}\ =\ \cases -Z(1), & \ \ k=0, \\
1/\la_{k}, & \ \ k\geq 1.
\endcases
$$
(We do not rule out the possibility that $Z(1)=0$.) 
Define the measure $d\nu_n$ on $[-1,1]$ by  
$$
d\nu_n(t)\ =\  \frac1{v} (1-t^2)^{(n/2)-1}\, dt,  \qquad v=\frac{(n-1)!\pi}{2^{n-1}(((n-1)/2)!)^2 } .
$$
The normalization ensures that $\nu_n([-1,1])=1$.  Fix $n$, and 
for $k=0,1,2,\dots,$ let $P_{k}(t)$ be the orthogonal polynomials on $[-1,1]$ with respect to this measure
normalized so that
$$
\int_{-1}^1 P_{k}^2(t)\, d\nu_n(t)\ =\ d_{k}.
$$
These are eigenfunctions for
$$
\Delta_n\ =\ -(1-t^2)^{1-(n/2)}\frac{d}{dt} (1-t^2)^{n/2} \frac{d}{dt}.
$$
The orthogonal projection $\Pi_{k}$ of $L^2(S^n)$
onto $\Cal H_{k}$ has kernel $K(x,y)=P_{k}(x\cdot y)$.
To prove Lemma 1.3.12, we just have to check that (6.1.1) holds for some $\phi\in \Cal H_2$.
 Fix a point $x_0\in S^n$ and set 
$$ \align
\phi(x)&=  \frac1{\la_2}  P_{2}(x,x_0)\qquad\text{ so that }\tag 6.1.2
\\
 \phi_{kk}(x)&= P_{n,2}(x,x_0),
\qquad \int \bar\phi^2\, d\mu_n =\frac{d_2}{\la_2^2}= \frac{n(n+3)}{8(n+1)^2}.
\endalign $$
 From this we can calculate the first term in (6.1.1), and Lemma 1.3.12 follows if
we show that 
$$
\TR \phi_{kk} \bar\Delta^{-1} \phi_{jj} \bar\Delta^{-1}\ =\ 
\frac{ (n+3)\left( n^2+n-4 - 4(n+1) \sum_{j=2}^{n-1} \frac1j\right)}{4(n-3)(n+1)^2}.
\tag 6.1.3
$$

Now
$$
\align
\trace \phi_{ii} \Pi_k \phi_{jj} \Pi_\ell
& =\ \int_{S^n\times S^n} \phi_{ii}(x) P_{k}(x\cdot y) \phi_{jj}(y) P_{\ell}(x\cdot y)\, d\mu_n  d\mu_n    \\
& =\ \int_{S^n\times S^n} P_2(x\cdot x_0) P_{k}(x\cdot y) P_2(y,\cdot x_0) P_{\ell}(x\cdot y)\, d\mu_n  d\mu_n  .
\endalign
$$
By symmetry this does not depend on $x_0$ and we can average over $x_0\in S^n$ to obtain
$$
\int_{S^n\times S^n} P_2(x\cdot y) P_{k}(x\cdot y) P_{\ell}(x\cdot y)\, d\mu_n  d\mu_n  
\ =\ \int_{-1}^1 P_2 P_k P_\ell \, d\nu_n  .
$$
Set  $p=(n-1)/2$, and $B=\Delta^2+p^2$ which has eigenvalue $(j+p)^2$ on $\Cal H_j$,
and set
$$
C(j,m,\ell)\ =\ \int_{-1}^1 P_j P_{j+m} P_\ell \, d\nu_n, \qquad\qquad
 F(j,m,\ell)\ =\ \frac{C(j,m,\ell)}{\la_j \la_{j+m}}.
\tag 6.1.4
$$
Then
$$
\align
&\hskip-12pt \TR \phi_{kk} \bar\Delta^{-1} \phi_{jj} \bar\Delta^{-1}\ =\ 
\trace \phi_{kk} \bar\Delta^{-1}B^{z/2} \phi_{jj} \bar\Delta^{-1}\bigl|^\mer_{z=0}\tag 6.1.5  
 \\
=\ &\sum_{j,k=0}^\infty \frac{ \trace \phi_{ii} \Pi_j \phi_{\ell\ell} \Pi_k}{\tilde\la_j \tilde\la_k}
(j+p)^{z}\biggl|^\mer_{z=0}\ =\ \sum_{j,k=0}^\infty \frac{ \int P_{j} P_{k} P_{2}\, d\nu_n   }{\tilde\la_j \tilde\la_k}
(j+p)^{z}\biggl|^\mer_{z=0}\\
   =\ &\frac2{\tilde\la_0}\frac{ d_{2} }{ \la_2} 
\ +\ \sum_{j,k=1}^\infty \frac{ \int P_{j} P_{k} P_{2}\, d\nu_n   }{\la_j \la_k}
(j+p)^{z}\biggl|^\mer_{z=0} \\
 =\ & \frac2{\tilde\la_0}\frac{ d_{2} }{ \la_2} 
\ +\ \Bigg( \sum_{j=1}^\infty F(j,0,2)  (j+p)^{z} 
\ +\ \sum_{j=1}^\infty  F(j,2,2)  (j+p)^{z}   \\
&\hskip2.75in  +\ \sum_{j=3}^\infty F(j,-2,2) (j+p)^{z}    
\Bigg)   \biggl|^\mer_{z=0}  \\
 =\ &\frac2{\tilde\la_0}\frac{ d_{2} }{ \la_2} 
\ +\  \sum_{j=1}^2 F(j,0,2) \ +\ \sum_{j=1}^2  F(j,2,2)    
\\
& \hskip1.25in+ \sum_{j=3}^\infty (F(j,0,2)+F(j,2,2)+ F(j,-2,2)) (j+p)^{z}  \biggl|^\mer_{z=0} . 
\endalign
$$ 
The final sum in (6.1.5) and the quantity $1/\tilde\la_0=Z(1)$  both involve a 
zeta regularization.  In order to evaluate these two quantities we make use of the following lemma.

\nonumproclaim{Lemma 6.1.1}
Suppose that $G$ is an even polynomial\/{\rm ;} that is{\rm ,} $G(-j)=G(j)${\rm ,} and $\ell>1${\rm .}  Then
$$
 \sum\Sb j\neq \ell \\ 1\leq j<\infty \endSb \frac{ G(j)}{j^2-\ell^2}j^{z}   \biggl|^\mer_{z=0}
\ =\ \frac{G(0)}{2\ell^2}\ +\ \frac{G(\ell)}{4\ell^2}\ -\ \frac{G^\prime(\ell)}{2\ell}.
$$
\endproclaim

\demo{Proof}
Set 
$$
H(j)\ =\ \cases \frac{ G(j)-G(\ell)}{j^2-\ell^2}  \qquad & j^2\neq \ell  \\
\frac{G^\prime(j)}{2\ell} \qquad & j^2=\ell^2.\endcases
$$
Then $H$ is an even polynomial and since 
$\zeta(-2k)=0$ for $k=1,2,\dots$, and $\zeta(0)=-1/2$, we see that 
$$
\align
 \sum_{j=1}^\infty  H(j) j^{z}\biggl|^\mer_{z=0}
& =\ \sum_{j=1}^\infty \sum_{k=0}^{\deg H} \frac{H^{(2k)}(0)}{(2k)!}j^{2k+z}\biggl|^\mer_{z=0} \\
&   =\  \sum_{k=0}^{\deg H}\frac{H^{(2k)}(0)}{(2k)!}\zeta(-2k)
\ =\  -\frac{H(0)}{2}\ =\ \frac{G(0)}{2\ell^2}-\frac{G(\ell)}{2\ell^2}.
\endalign
$$
Hence
$$
\align
  \sum_{j=1}^{\ell-1}   \frac{ G(j)}{j^2-\ell^2}& +\ \frac{G^\prime(\ell)}{2\ell}
\ +\ \sum_{j=\ell+1}^\infty \frac{ G(j)}{j^2-\ell^2} j^{z}   \biggl|^\mer_{z=0} \tag 6.1.6
 \\ \noalign{\vskip4pt}
& =\ \sum_{j=1}^\infty  H(j)  j^{z}   \biggl|^\mer_{z=0} 
\ +\ G(\ell) \left( \sum_{j=1}^{\ell-1} \frac{ 1}{j^2-\ell^2}  
\ +\ \sum_{j=\ell+1}^\infty \frac{ 1}{j^2-\ell^2} \right)  \\ \noalign{\vskip4pt}
&=\ \frac{G(0)}{2\ell^2}-\frac{G(\ell)}{2\ell^2}\ +\ \frac{3G(\ell)}{4\ell^2}.\\
\noalign{\vskip-24pt}
\endalign
$$
\enddemo
\phantom{not yet}
\vglue-6pt

 From Lemma 6.1.1, we get the following.

\phantom{not yet}

\nonumproclaim{Lemma 6.1.2}
$$
Z(1)\ =\ -\frac1{n-1} \sum_{j=1}^{n-1} \frac1j.
$$
\endproclaim

\phantom{not yet}

\demo{Proof}
As before, set $p=(n-1)/2$, and $B=\Delta^2+p^2$ which has eigenvalue $(j+p)^2$ on $\Cal H_j$.
Using the canonical trace, we have
$$
\align
 Z(1)& \ =\ \TR \Delta^{-1}B^{z/2}\bigl|^\mer_{z=0}
\ =\ \sum_{j=1}^\infty \frac{d_{j}}{\la_{j}}(j+p)^{z}\biggl|^\mer_{z=0}  \\ \noalign{\vskip4pt}
& \ =\ \frac{2}{(n-1)!} \sum_{j=1}^\infty \frac{ \prod_{k=0}^{p-1}( (j+p)^2 -k^2)}{(j+p)^2 -p^2}
(j+p)^{z}\biggl|^\mer_{z=0}   \\ \noalign{\vskip4pt}
&\ =\ \frac{2}{(n-1)!} \sum_{j=p+1}^\infty \frac{ \prod_{k=0}^{p-1}( j^2 -k^2)}{j^2 -p^2}
j^{z}\biggl|^\mer_{z=0}.   \\
\endalign
$$
Now, applying  Lemma 6.1.1 with $\ell=p$ and 
$$
G(j)\ =\ \frac{2}{(n-1)!} \prod_{k=0}^{p-1}( j^2 -k^2)\ =\ \frac{2j(j-p+1)\dots(j+p-1)}{(n-1)!} ,
$$
we find that
$$
\align
& G(j)\ =\ 0\qquad\qquad \text{ for }|j|<p,  \\ \noalign{\vskip4pt}
& G(p)\ =\ 1, \\  \noalign{\vskip4pt}
& G^\prime(p)\ =\ \sum_{j=1}^{n-1}\frac1j \ +\ \frac1{n-1},
\endalign
$$
and so
\vglue19pt
 \hfill ${\displaystyle
\sum_{j=p+1}^\infty \frac{ G(j)}{j^2-\ell^2}j^{z}   \biggl|^\mer_{z=0}\ =\ 
\frac{G(p)}{4p^2}\ -\ \frac{G^\prime(p)}{2p}\ =\ -\frac1{n-1}\sum_{j=1}^{n-1}\frac1j.
}$ \hfill
\enddemo

\phantom{I wanna}

 From Lemma 6.1.2, we can  compute the first term in (6.1.5) and get
$$
\frac2{\tilde\la_0}\frac{ d_{2} }{ \la_2} \ =\ \frac{n(n+3)}{2(n-1)(n+1)} \sum_{j=1}^{n-1} \frac1j.
\tag 6.1.7
$$
To compute the other terms in (6.1.5), 
standard techniques give the recurrence relation
$$
xP _{k}(x)\ =\ \frac{k+1}{2k+n+1}P _{k+1}\ +\ 
\frac{k+n-2}{2k+n-3}P _{k-1}.
$$
 From this it can be deduced that 
$$
P_2(t)\ =\ \frac12 (n+3)\, ((n+1)x^2-1),
$$
and the linearization coefficients defined in  (6.1.4), are  computed as in [V, 9.11.6]:
$$
\matrix
\noalign{\vskip3pt}
\hfill C(j,0,2)&\hskip-6pt  = \, 
{  \frac{n+3}{2(n-2)!} }(j)\cdots(j+p-2)(j+p)^2 (j+2)\cdots (j+n-1)\hfill& j\geq 0, \hfill\\ \noalign{\vskip4pt}
 C(j,2,2)&\hskip-6pt  = \, {  \frac{(n+1)(n+3)}{4(n-1)!}} (j+1)\cdots (j+p) (j+p+2)\cdots(j+n),\hfill& j\geq 0,\hfill\\
\noalign{\vskip4pt}
 C(j,-2,2)&\hskip-6pt  = \, {  \frac{(n+1)(n+3)}{4(n-1)!} }(j-1)\cdots (j+p-2) (j+p)\cdots(j+n-2), \hfill& j\geq\hfill
2. 
\endmatrix \tag 6.1.8
$$
Hence
$$ \align
\sum_{j=1}^2 F(j,0,2) \ +\ \sum_{j=1}^2  F(j,2,2)
 =\ &\frac{n+1}n \ +\ \frac{(n-1)n(n+3)^2}{4(n+1)^2(n+5)}\tag 6.1.9 \\ 
& +\ \frac{(n+1)^2}{6(n+2)}
\ +\ \frac{n(n+1)(n+2)}{32(n+5)}.
\endalign
$$
We use Lemma 6.1.1 to compute the final sum in (6.1.5). 
 First we write this term in the form 
$$
 \sum_{j=p+3}^\infty (F(j-p,0,2)+F(j-p,2,2)+ F(j-p,-2,2)) j^{z}  \biggl|^\mer_{z=0}.
\tag 6.1.10
$$
Using the formulas for the linearization coefficients we get
$$
\align
\noalign{\vskip-8pt}
&\tag 6.1.11\\
\noalign{\vskip-14pt}
&\hskip-18pt F(j-p,0,2)+F(j-p,2,2)+F(j-p,-2,2)
\\ \noalign{\vskip4pt}
=\ & \frac{n+3}{2(n-1)!}\biggl( \frac{ 2p(j-p)\dots(j-2) j^2 (j+2)\dots(j+p)}{(j-p)^2(j+p)^2}
\\  \noalign{\vskip4pt}
&\hskip1.5in+\ \frac{(p+1)(j-p+1)\dots j(j+2)\dots(j+p+1)}{(j-p)(j+p)(j-p+2)(j+p+2)}  \\  \noalign{\vskip4pt}
&\hskip1.5in   +\ \frac{(p+1)(j-p-1)\dots(j-2)j\dots(j+p-1)}{(j-p)(j+p)(j-p-2)(j+p-2)} \biggr) \\  \noalign{\vskip4pt}
=\ &\frac{n+3}{(n-1)!} \frac{ (j-p+1)\dots(j+p-1)}{(j-p)(j+p)(j-p-2)(j+p-2)(j-p+2)(j+p+2)} j Q(j)
\endalign
$$
where  
$$
 Q(j)\ =\ nj^2\ -\ \frac{n(n+3)^2}4 +8.
$$
Writing
$$
G(j)\ =\ \frac{n+3}{(n-1)!} (j-p+1) (j-p+3)\dots(j+p-3)(j+p-1)j Q(j),
$$
we see that (6.1.10) equals 
$$ \multline
  \sum_{j=p+3}^\infty \frac{G(j)}{(j^2-p^2)(j^2-(p+2)^2)}\\
\ =\ \frac1{(p+2)^2-p^2}
\biggl( \sum_{j=p+3}^\infty \frac{G(j)}{(j^2-(p+2)^2)} 
 -\sum_{j=p+3}^\infty \frac{G(j)}{(j^2-p^2)}\biggr) ,  \endmultline
$$
to which we can apply Lemma 6.1.1. 
Using the fact that $G(j)=0$ for $0\leq j\leq p-3$ and $j=p-1$,  we find that (6.1.10) equals
$$
\multline
\frac{1}{4(p+1)}
\biggl(-\frac{G(p-2)}{(p-2)^2-(p+2)^2}\ -\ \frac{G(p)}{(p^2-(p+2)^2)}
\ -\ \frac{G(p+1)}{(p+1)^2-(p+2)^2}\\  \noalign{\vskip4pt}
  +\ \frac{G(p+2)}{4(p+2)^2}  \ -\ \frac{G^\prime(p+2)}{2(p+2)}  
 \ +\ \frac{G(p-2)}{(p-2)^2-p^2}\ +\ \frac{G(p+1)}{(p+1)^2-p^2}  
\\  \noalign{\vskip4pt}
  +\ \frac{G(p+2)}{(p+2)^2-p^2}\ -\ \frac{G(p)}{4 p^2}
\ +\ \frac{G^\prime(p)}{2p} \biggr) .\\
\endmultline
$$
We evaluate this  using the values
$$
\align
  G(p-2)  =\ &\frac{ 2(n+1)(n+3)}{(n-1)(n-3)},   \qquad G(p)\ =\ -\frac{(n+3)(n^2+n-4)}{2(n-3)}, \\  \noalign{\vskip4pt}
 G(p+1)  =\ & -\frac{(n+1)(n+3)(n+4)}{6},   \qquad G(p+2)\ =\ \frac{n(n+3)^2}{2(n-1)}, \\  
 G^\prime(p)   =\ &\frac{(n-1)n(n+3)}{4(n-3)}\\   
& -\ \frac{(n+3)(n^2+n-4)}{2(n-3)}
\left( \sum_{j=2}^{n-1}\frac1j\ +\ \frac12\ -\ \frac1{n-3}\  +\ \frac1{n-1}\right),\\   
  G^\prime(p+2)  =\ &\frac{n(n+3)^2}{2(n-1)}\left( \sum_{j=3}^{n-2}\frac1j\ -\ \frac14
\ +\ \frac1n \ +\ \frac2{n+3}\ +\ \frac{n(n+3)}8.\right),
\endalign
$$
and we find that (6.1.7), (6.1.9) and  (6.1.10) sum to  the right-hand side of~(6.1.3).  
\enddemo

\AuthorRefNames [BFK1]
\references

[Br]   
\name{T. Branson},
Sharp inequalities, the functional determinant, and the
complementary series,  
{\it Trans\. A.M.S.}
{\bf 347} (1995), 3671--3742.

[BFK1]
\name{D.  Burghelea, L. Friedlander}, and \name{T. Kappeler},  
On the determinant of elliptic differential and finite difference 
operators in vector 
bundles over $S^1$,
{\it Comm.\ Math Phys\/}.\ {\bf 138} (1991), 1--18.

[BFK2]
\bibline,  
Meyer-Vietoris type formula for determinants of elliptic differential
operators,
{\it J. Funct. Anal\/}.\ {\bf 107} (1992), 34--65.

[BO1]
\name{T. Branson} and \name{B. {\O}rstead}, 
Conformal geometry and global invariants,
{\it Differential Geom.\  and Appl\/}.\ {\bf 1}  
 (1991), 279--308.

[BO2]
\bibline, 
Explicit functional determinants in four dimensions,
{\it Proc. A.M.S.} {\bf 113}  
(1991), 669--682.

[BCY]
\name{T. Branson, S.-Y.\ A.\ Chang}, and \name{P. Yang},
Estimates  and extremals for zeta function determinants on four-manifolds,
{\it Commun.\ Math.\ Phys.\/} {\bf 149} (1992), 241--262.

[CQ1]
\name{S.-Y.\ A. Chang}  and \name{J. Qing},
Zeta functional  determinants on manifolds with boundary,
{\it Math.\ Res.\ Lett.\/} {\bf 3} (1996), 1--17.

[CQ2]
\bibline, The zeta functional determinants on manifolds with boundary. 
I. The formula. 
{\it J. Funct.\ Anal.\/} {\bf 147} (1997), 
327--362.

[CQ3]
\bibline,
The zeta functional determinants on manifolds with boundary. II. Extremal metrics and compactness of
isospectral set, 
{\it J. Funct.\ Anal.\/}
{\bf 147} (1997),  363--399.

[CY]
\name{S.-Y.\ A.\ Chang} and \name{P. Yang},
Extremal metrics of zeta function determinants on $4$-manifolds,
{\it Ann.\ of Math.\/} {\bf 142} (1995), 171--212.

[Chi]
\name{P. Chiu},
Height of flat tori,
{\it Proc.\ A.M.S.} {\bf 125} (1997), 723--730.

[Cho]
\name{B. Chow},
The Ricci flow on the $2$-sphere,
{\it J. Differential Geom.\/} {\bf 33} (1991), 325--334.

[Eb]
\name{D. Ebin},
On the space of Riemannian metrics,
{\it Bull.\ A.M.S.}
{\bf 75} (1968), 1001--1003.

[Fo1]
\name{R. Forman},
Functional determinants and geometry,
{\it Invent. Math.} {\bf 88}  (1987),  447--493; Erratum
{\it ibid\/}.\ {\bf 108} (1992), 453--454.

[Fo2]
\bibline,
Determinant of the Laplacian on graphs,
{\it Topology} {\bf 32}  (1993), 35--46.

[Fo3]
\bibline,
Determinants, finite-difference operators and
boundary value problems,
{\it Commun.\ Math.\ Phys.\/} {\bf 147} (1992), 485--526.

[Ha]
\name{R. Hamilton},
 The Ricci flow on surfaces,  in {\it Mathematics and General Relativity\/}
 (Santa Cruz, 1986), 237--262, 
{\it Contemp.\ Math.\/} {\bf 71}, A.M.S., Providence, RI  (1988).

[Ka]
\name{C. Kassel},
Le r\'esidu non commutatif (d'apres M.\ Wodzicki),
{\it S{\rm \'{\rm e}}m.\  Bourbaki}, 
Exp.\ No.\ 708 (1988-89),  199--229.

[Gui]
\name{V. Guillemin},
A new proof of Weyl's formula on the asymptotic distribution of eigenvalues,
{\it Adv.\ in Math.\/}
{\bf 55}  (1985), 131--160.

[Gur]
\name{M. Gursky},
Uniqueness of the functional determinant,
{\it Commun.\ Math.\ Phys.\/}
{\bf 189} (1997), 655--665.

[KV]
\name{M. Kontseivich} and \name{S. Vishik},
Geometry of determinants of elliptic operators, in 
{\it Funct.\ Anal.\ on the Eve of the $21^{\rm st}$ Century}, {\it Vol}.\  1, 
Birkh\"auser  Boston, MA, {\it Progr.\ Math\/}.\ {\bf 131} (1995),
173--197.

[LP]
\name{J. Lee} and \name{T. Parker},
The Yamabe problem,
{\it Bull.\ A.M.S.} {\bf 17} (1987), 37--91.

[Ok1]
\name{K. Okikiolu},
The Campbell-Hausdorff theorem for elliptic 
operators and a related trace formula, 
{\it  Duke Math.\ J.}  {\bf 79} (1995), 687--722.

[Ok2]
\bibline,
 The multiplicative anomaly for determinants  
of  elliptic operators, 
{\it  Duke Math.\ J.} {\bf 79} (1995), 723--750.

[Ok3]
\bibline,
Critical metrics for   spectral
zeta functions, 
 preprint.

[Ok4]
\bibline,
A note on the space of metrics modulo diffeomorphisms,
in preparation.

[On]
\name{E.\ Onofri},
On the positivity of the effective action in a theory of random surfaces,
{\it Commun.\ Math.\ Phys.\/} {\bf 86} (1982), 321--326.

[OPS1]
\name{B. Osgood, R. Phillips}, and \name{P. Sarnak}, 
Extremals of determinants of Laplacians, 
{\it  J.  Funct.\ Anal.\/} {\bf  80} (1988), 148--211. 

[OPS2]
\bibline, 
Compact isospectral sets of surfaces, 
{\it  J.  Funct.\ Anal.\/}  {\bf 80} (1988), 212--234.

[OPS3]
\bibline, 
Moduli space, heights and isospectral sets of plane domains, 
{\it  Ann.\  of  Math.\/}  {\bf 129} (1989), 293--362.

[PR]
\name{T. Parker} and \name{S. Rosenberg},
Invariants of conformal Laplacians,
{\it J. Differential Geom.\/}  {\bf 25} (1987), 199--222.

[Po]
\name{A. Polyakov},
Quantum geometry of bosonic strings,
{\it Phys.\ Lett.\ B} {\bf 103} (1981), 207--210.

[RS]
\name{D. Ray} and \name{I. Singer},  
$R$-torsion and the Laplacian on  Riemannian manifolds,
{\it Adv.\ in Math.\/} {\bf 7} (1971), 145--210.

[Ri]
\name{K. Richardson}, 
{Critical points of the determinant of the Laplace operator},
{\it J.  Funct.\ Anal.\/} {\bf 122} (1994), 52--83.

[Sp]
\name{M. Spivak},
{\it A Comprehensive Introduction to Differential Geometry},
{\it Vol\/}.\ 1,  Publish or Perish, Wilmington, Del.\  (1979).

[Ta]
\name{M. Taylor},
{\it Partial Differential Equations\/} I., {\it Basic Theory\/},
{\it Appl.\ Math.\ Sci.\/} {\bf 115}, Springer-Verlag,
New York,  1996.

[Tr]
 \name{F. Tr\`eves},
{\it Basic Linear Partial Differential Equations},
{\it Pure and Appl.\ Math\/}.\ {\bf 62}, Academic Press, New York (1975).

[Vi]
\name{N.  Vilenkin}, 
{\it Special Functions and the Theory of Group Representations},
 {\it Transl.\  Math. Mono\/}.\ {\bf 22}, A.M.S., Providence,
 RI (1968).

[Wo]
\name{M. Wodzicki},  
Noncommutative residue.\ I.\  Fundamentals,
in {\it $K$-Theory, Arithmetic and Geometry} (Moscow, 1984--1986), 
 {\it Lecture Notes in Math.\/} {\bf 1289} (1987), 320--399,
 Springer-Verlag, New York.

\endreferences

\phantom{oh my}
\vglue-30pt

\centerline{\ninerm (Received September 20, 1999)}
\bye
\end